\documentclass[a4paper,10pt]{article}
\usepackage{amssymb}
\usepackage{latexsym}
\usepackage{amsmath}
\usepackage{amsthm}
\usepackage{amsfonts}
\usepackage{amssymb}
\usepackage{graphicx}

\newtheorem{Th}{Theorem}
\newtheorem{Prop}{Proposition}
\newtheorem{Co}{Corollary}
\newtheorem{Lm}{Lemma}
\newtheorem{Lma}{Lemma}[section]
\newtheorem{Dfi}{Definition}
\newtheorem{Rm}{Remark}

\newcommand{\be}{\begin{equation}}
\newcommand{\ee}{\end{equation}}
\newcommand{\R}{\mathbb{R}}
\newcommand{\N}{\mathbb{N}}
\newcommand{\C}{\mathbb{C}}
\newcommand{\Z}{\mathbb{Z}}

\newcommand\res{\mathop{\hbox{\vrule height 7pt width .5pt depth 0pt
\vrule height .5pt width 6pt depth 0pt}}\nolimits}

\def\theequation{\thesection.\arabic{equation}}
\def\theTh{\Roman{section}.\arabic{Th}}
\def\theProp{\Roman{section}.\arabic{Prop}}
\def\theCo{\Roman{section}.\arabic{Co}}

\def\theRm{\Roman{section}.\arabic{Rm}}
\newcommand{\reset}{\setcounter{equation}{0}\setcounter{Th}{0}\setcounter{Prop}{0}\setcounter{Co}{0}
\setcounter{Lm}{0}\setcounter{Rm}{0}}
\def\La{\Lambda}
\def\La{\Lambda}
\def\ti{\tilde}
\def\lf{\left}
\def\rg{\right}

\def\al{\alpha}
\def\la{\lambda}

\def\ep{\varepsilon}
\def\ds{\displaystyle}
\def\ov{\overline}
\def\Om{\Omega}
\def\om{\omega}
\def\p{\partial}

\def\vP{\vec{\Phi}}

\def\vq{\vec{q}}
\def\gP{g_{\vec{\Phi}}}
\def\gI{|\vec{\mathbb I}_{\vP}|^2_{\gP}}

\def\res{\mathop{\hbox{\vrule height 7pt width .5pt 
depth 0pt\vrule height .5pt width 6pt depth 0pt}}\nolimits}

\begin{document}
\title{A Viscosity Method in the Min-max Theory of Minimal Surfaces.}
\author{ Tristan Rivi\`ere\footnote{Department of Mathematics, ETH Zentrum,
CH-8093 Z\"urich, Switzerland.}}
\date{ }
\maketitle

{\bf Abstract :}{\it We present the min-max construction of critical points of the area using penalization arguments. Precisely, for any immersion of a closed surface $\Sigma$ into a given closed manifold, we add to the area Lagrangian a  term equal to the $L^q$
norm of the second fundamental form of the immersion times a ``viscosity'' parameter.  This relaxation of the area functional satisfies the Palais-Smale condition  for $q>2$. This permits to construct critical points of the relaxed Lagrangian using classical min-max arguments such as the mountain pass lemma. The goal of this work is to describe the passage to the limit when the ``viscosity'' parameter tends to zero. Under some natural entropy condition, we establish a varifold convergence of these critical points towards a parametrized integer stationary varifold realizing the min-max value. It is proved in \cite{PiRi} that parametrized integer stationary varifold  are given by smooth maps exclusively.  As a consequence we conclude that every surface area minmax  is realized by a smooth possibly branched minimal immersion.}

\medskip

\noindent{\bf Math. Class. 49Q05, 53A10, 49Q15, 58E12, 58E20}
\section{Introduction}
The study of minimal surfaces, critical points of the area, has stimulated the development of entire fields in analysis and in geometry. The calculus of variations is one of them. The origin of the field is very much linked to the question of proving the existence of minimal
2-dimensional discs bounding a given curve in the euclidian 3-dimensional space and minimizing the area. This question, known as {\it Plateau Problem}, has been posed since the XVIIIth century by Joseph-Louis Lagrange,
the founder of the Calculus of Variation after Leonhard Euler. This question has been ultimately solved independently by Jesse Douglas and Tibor Rad\'o around 1930. In brief the main strategy of the proofs was to minimize the Dirichlet energy instead of the area, which is lacking coercivity properties, the two lagrangians being identical on conformal maps. After these proofs, successful attempts  have been made to solve the Plateau Problem
in much more general frameworks. This has been in particular at the origin of the field of {\it Geometric Measure Theory} during the 50's, where the notions of rectifiable current which were proved to be the {\it ad-hoc} objects
for the minimization process of the area (or the mass in general) in the most general setting.

The search of absolute or even local minimizers is of course the first step in the study of the variations of a given lagrangians but is far from being exhaustive while studying the whole set of critical points. In many problems there is even no minimizer at all, this is for instance the case of closed surfaces in simply connected manifolds with also trivial two dimensional homotopy groups. This problem is already present in the 1-dimensional counter-part of minimal surfaces, the study of closed geodesics. For instance in a sub-manifold of ${\R}^3$ diffeomorphic to $S^2$  there is obviously no closed geodesic minimizing the length. In order to construct closed geodesics in such manifold, Birkhoff around 1915 introduced a technic called ''min-max'' which permits to generate  critical points of the length with non trivial index. In two words this technic consists in considering the space of paths of closed curves within a non-trivial homotopy classes of paths in the sub-manifold (called ``sweep-out'') and to minimize, out of all such paths or ''sweep-outs'', the maximal length of the curves realizing each ''sweep-out''. In order to do so, one is facing the difficulty posed by a lack of coercivity of the length with respect to this minimization process within this ''huge space'' of sweep-outs. In order to ''project'' the problem to a much smaller space of ''sweep-outs'' in which the length would become more coercive, George Birkhoff replaced each path by a more regular one made of very particular closed curves joining finitely many points with portions of geodesics minimizing the length between these points. This replacement method also called nowadays ``curve shortening process'' has been generalized in many situations in order to perform min-max arguments.

Back to minimal surfaces, in a series of two works (see \cite{CM1} and \cite{CM2}), Tobias Colding and Bill Minicozzi, construct by min-max methods minimal 2 dimensional spheres in homotopy 3-spheres (The analysis carries over to general target riemannian manifolds). The main strategy
of the proof combines the original approach of Douglas and Rad\'o, consisting in replacing the area functional by the Dirichlet energy, with a ''Birkhoff type'' argument of optimal replacements. Locally to any map from a given
''sweep-out'' one performs a surgery, replacing the map itself by an harmonic extension minimizing the Dirichlet energy. The convergence of such a ''harmonic replacement'' procedure, corresponding in some sense to Birkhoff ``curve shortening procedure'' in one dimension, is ensured by a fundamental result regarding the local convexity
of the Dirichlet energy into a manifold under small energy assumption and a unique continuation type property. What makes possible the use of the Dirichlet energy instead of the area functional, as in \cite{SU}, is the fact that the domain $S^2$ posses
only one conformal structure and modulo a re-parametrization any $W^{1,2}$ map can be made $\ep-$conformal (due to a fundamental result of Charles Morrey see theorem I.2 \cite{Mo}). This is not anymore the case if one wants to
extend Colding-Minicozzi's approach to general surfaces. This has been done however successfully by Zhou Xin in \cite{Xi1} and \cite{Xi2} following the original Colding-Minicozzi approach. These papers are based on an involved argument in which to any ''sweep-out'' 
of $W^{1,2}-$maps a path of smooth conformal structures together with a  path of re-parametrization  are assigned in order to be as close as possible to paths of conformal maps.

Because of the finite dimensional nature of the moduli space of conformal structures in 2-D, and the ``optimal properties'' of the Dirichlet energy,  Colding-Minicozzi's min-max method is intrinsically linked to two dimensions as Douglas-Rad\'o's resolution of the Plateau problem was too. The field of Geometric Measure Theory, which was originally designed to remedy to this limitation and to
solve the Plateau Problem for arbitrary dimensions in various homology classes, has been initially developed with a minimization perspective and the framework of rectifiable currents as well as the lower semicontinuity of the mass for weakly converging sequences was matching perfectly this goal. In order to solve min-max problems in the general framework of Geometric Measure Theory, the notion of {\it varifold} has been successfully introduced by  William Allard 
and by Fred Almgren. A complete GMT min-max procedure has been finally set up by Jon Pitts in \cite{Pit} who introduced the notion of {\it almost minimizing varifolds} and developed their regularity theory in co-dimension 1. Constructive comparison arguments as well as combinatorial type arguments are also needed in this rather involved and general procedure (The reader is invited also to consult \cite{CD} and \cite{MN1} for thorough presentations of the GMT approaches to min-max
procedures).

\medskip

The aim of the present work is to present a direct min-max approach for constructing minimal surfaces  in a given closed sub-manifold $N^n$ of ${\R}^m$.  The general scheme is simple : one works with a special subspace of $C^1$ immersions of a given surface $\Sigma$, one adds to the area of each of such an immersion $\vec{\Phi}$ a relaxing ``curvature type'' functional multiplied by a small viscous parameter $\sigma^2$
\be
\label{0-0}
A^\sigma(\vec{\Phi}):=\mbox{Area}(\vec{\Phi})+\sigma^2\int_\Sigma\mbox{ curvature terms} \ dvol_{g_{\vec{\Phi}}} 
\ee
where $dvol_{g_{\vec{\Phi}}}$ is the volume form on $\Sigma$ induced by the immersion $\vec{\Phi}$. The ``curvature terms'' is chosen in order to ensure that $A^\sigma$ satisfies the {\it Palais-Smale property} on the {\it ad-hoc} 
corresponding Finsler manifold of $C^1-$immersions. This offers the suitable framework in which Palais deformation theory can be applied to produce critical points realizing an arbitrary minmax value. Once a min-max critical point of $A^\sigma$ is produced one passes to the limit
$\sigma\rightarrow 0\cdots$

\medskip

More precisely, we introduce the space ${\mathcal E}_{\Sigma,p}$ of $W^{2,2p}-$immersions $\vec{\Phi}$
of a given closed surface $\Sigma$ for $p>1$ \footnote{The condition $p>1$ ensures that $\vec{\Phi}$ is $C^1$. This last fact permits to use the classical definition of an immersion. The case $p=1$ was considered in previous works by the author
where the notion of immersion had to be weakened.} into $N^n\subset{\R}^m$. It is proved below that this space has a nice structure of {\it Banach Manifold}. For such immersions we consider the relaxed energy
\[
A^\sigma(\vec{\Phi}):=\mbox{Area}(\vec{\Phi})+\sigma^2\int_{\Sigma}\lf[1+|\vec{\mathbb I}_{\vec{\Phi}}|^{2}\rg]^p\ dvol_{g_{\vec{\Phi}}}
\]
where $g_{\vec{\Phi}}$ and $\vec{\mathbb I}_{\vec{\Phi}}$ are respectively the first and second fundamental forms of $\vec{\Phi}(\Sigma)$ in $N^n$. Unlike previous existing viscous relaxations for min-max problems in the literature, the energy $A^\sigma$ is intrinsic in the sense that it is invariant under re-parametrization of $\vec{\Phi}$ : $A^\sigma(\vec{\Phi})=A^\sigma(\vec{\Phi}\circ\Psi)$ for any smooth diffeomorphism $\Psi$ of $\Sigma$. Modulo a choice of parametrization it is proved in \cite{La} and \cite{KLL} that for a fixed $\sigma\ne 0$ the Lagrangian $A^\sigma$ satisfies the {\it Palais-Smale condition}.
Hence we can consider applying the mountain path lemma to this Lagrangian. 

We introduce now the following definition

\begin{Dfi}
\label{df-target-harmonic}
Let $\Sigma$ be a closed Riemann surface and $M^m$ be a closed sub-manifold $M^m\subset{\R}^Q$. A map $\vec{\Phi}\in W^{1,2}(\Sigma,M^m)$ together with an $L^\infty$ bounded integer multiplicity $N_x$ is called `` integer target harmonic" if for almost every\footnote{The notion of {\it almost every domain} means for every smooth domain $\Om$ and any smooth function $f$ such that $f^{-1}(0)=\p\Om$ and $\nabla f\ne 0$ on $\p\Om$ then for almost every $t$ close enough to zero and regular value for $f$ one considers the domains contained in $\Omega$ or containing $\Omega$ and bounded by $f^{-1}(\{t\})$.}
domain $\Om\subset \Sigma$ and any smooth function $F$ supported in the complement of an open neighborhood of $\vec{\Phi}(\p \Om)$ we have
\be
\label{targ}
\int_{\Om} N_x\ \lf<d(F(\vec{\Phi})), d\vec{\Phi}\rg>_{g_0}- N_x\ F(\vec{\Phi})\ A(\vec{\Phi})(d\vec{\Phi},d\vec{\Phi})_{g_0}\ dvol_{g_0}
\ee
where $g_0$ is an arbitrary metric whose conformal structure is the one given by the Riemann surface $\Sigma$, $<\cdot,\cdot>_{g_0}$ denotes the scalar product in $T^\ast\Sigma$ issued from $g_0$,
$A(\vec{y})$ denotes the second fundamental form of $M^m\subset {\R}^m$ at the point $\vec{y}$. When the function $N$ is constant on $\Sigma$ we simply speak about ``target harmonic'' maps.
\hfill $\Box$
\end{Dfi}
Our main result in the present work is the following convergence theorem.
\begin{Th}
\label{th-main}
Let $N^n$ be a closed $n-$dimensional sub-manifold of ${\R}^m$ with $3\le n\le m-1$ being arbitrary. Let $\Sigma$ be an arbitrary closed Riemanian 2-dimensional manifold. Let $\sigma_k\rightarrow 0$ and let $\vec{\Phi}_k$ be a sequence of critical points
of 
\[
A^{\sigma_k}(\vec{\Phi}):=\mbox{Area}(\vec{\Phi})+\sigma_k^2\int_{\Sigma}\lf[1+|\vec{\mathbb I}_{\vec{\Phi}}|^{2}\rg]^p\ dvol_{g_{\vec{\Phi}}}
\]
 in the space of $W^{2,2p}-$immersions  of $\Sigma$ and satisfying the entropy condition
\be
\label{iii}
\sigma_k^2\int_{\Sigma}\lf[1+|\vec{\mathbb I}_{\vec{\Phi}_k}|^{2}\rg]^p\ dvol_{g_{\vec{\Phi}_k}}=o\lf( \frac{1}{\log\, \sigma_k^{-1}}  \rg)
\ee
Then, modulo extraction of a subsequence, there exists a closed Riemann surface $(S,h_0)$ with
 \[
\mbox{genus}(S)\le\mbox{genus}(\Sigma) 
\]
and a conformal integer target harmonic
map $(\vec{\Phi}_\infty,N)$ from $S$ into $N^n$ such that
\[
\lim_{k\rightarrow +\infty} A^{\sigma_k}(\vec{\Phi}_k)=\frac{1}{2}\int_SN\,|d\vec{\Phi}_\infty|_{h_0}^2 \ dvol_{h_0}
\]
Moreover, the oriented varifold associated to $\vec{\Phi}_k$ converges in the sense of Radon measures towards the stationary integer varifold associated to $(\vec{\Phi}_\infty,N)$. \hfill $\Box$
\end{Th}
The regularity of target harmonic maps is established in \cite{Ri5} and \cite{PiRi}. According to the main results in these works the limit $(\vec{\Phi}_\infty,N)$ is a smooth minimal branched immersion equipped with a smooth integer valued multiplicity.

\medskip

{\bf Open problem}: Assuming $\vec{\Phi}_k$ has a uniformly bounded {\it Morse index} for the Lagrangian $A^{\sigma_k}$ one  expects that the convergence is a strong $W^{1,2}-$ ``bubble tree'' convergence (i.e. strong away from finitely many points) which is equivalent to $N\equiv 1$ on $S$.

\medskip

The main difficulty in proving theorem~\ref{th-main} in contrast with existing non intrinsic viscous approximations of min-max procedures in the literature is that there is a-priori no $\epsilon-$regularity property \underbar{independent} of the viscosity $\sigma$ available. Indeed the following result
is proved in \cite{MR}.
\begin{Prop}
\label{pr-no-ep-reg}
There exists $\vec{\Phi}_k\in C^\infty(T^2,S^3)$ and $\sigma_k\rightarrow 0$ such that $\vec{\Phi}_k$ is a sequence of  immersions, critical points of $A^{\sigma_k}$, which is conformal into $S^3$ from a converging sequence of flat tori ${\R}^2/{\Z}+(a_k+i\,b_k){\Z}$  towards ${\R}^2/{\Z}+(a_\infty+i\,b_\infty){\Z}$, for which
\[
\limsup_{k\rightarrow +\infty}A^{\sigma_k}(\vec{\Phi}_k)<+\infty
\]
such that also $\vec{\Phi}_k$  weakly converges to a limiting map $\vec{\Phi}_\infty$ in $W^{1,2}({\R}^2/{\Z}+(a_\infty+i\,b_\infty){\Z},S^3)$ but $\vec{\Phi}_k$ nowhere strongly converges : precisely
\[
\forall\ U\ \mbox{open set in }{\R}^2/{\Z}+(a_\infty+i\,b_\infty){\Z} \quad\quad\int_{U} |\nabla \vec{\Phi}_\infty|^2\ dx^2<\liminf_{k\rightarrow +\infty}\int_{U} |\nabla \vec{\Phi}_k|^2\ dx^2\quad.
\]
\hfill $\Box$
\end{Prop}
In order to overcome this major difficulty in the passage to the limit $\sigma_k\rightarrow 0$ we prove a quantization result, lemma~\ref{lm-quanti}, which roughly says that there is a positive number $Q_0$, depending only on the target $N^n\subset {\R}^m$, below which for $k$ large enough, under the entropy condition assumption, there is no critical point of $A^{\sigma_k}$. This result is used at several stages in the proof. The main strategy goes as follows. We first establish the stationarity of the limiting varifold. The proof is based on an {\it almost divergence form} of the Euler Lagrange equation associated to $A^\sigma$ following the approach introduced in \cite{Ri2} for the Willmore Lagrangian in ${\R}^m$. The existence of such an {\it almost divergence form} is due to the symmetry group associated to the same Lagrangian in flat space and the application of Noether theorem (see \cite{Be}). As in \cite{MR}, the exact divergence form in Euclidian space is just an {\it almost-divergence form} in manifold. Next  we choose a conformal  parametrization of $\vec{\Phi}_k$ on a possibly degenerating sequence of Riemann surfaces $(\Sigma,h_k)$ (where $h_k$ denotes the  constant curvature metric of volume 1 conformally equivalent to $\vec{\Phi}_k^\ast g_{N^n}$). 
We use Deligne Mumford compactification in order to make converge $(\Sigma,h_k)$ towards a nodal Riemann surface with punctures (see for instance \cite{Hum}). We then use the monotonicity formula, deduced from the stationarity, in order to prove that, away from a so called {\it oscillation set}, the limiting volume density measure on the thick parts of the limiting nodal surface is absolutely continuous with respect to the Lebesgue measure. We then use the monotonicity formula again in order to prove the quantization result lemma~\ref{lm-quanti}. This quantization result is used in order to show that the limiting volume density measure  restricted to the oscillation set is equal to finitely many Dirac masses.
The quantization result is again used in order to prove that for the weakly converging sequence $\vec{\Phi}_k$   there is no energy loss neither in the necks in
each thick parts of the limiting nodal surface, nor in the collars regions separating possible bubbles, which are possibly formed (see lemma~\ref{lm-no-neck}). The previous results are proved to show the rectifiability of the limiting varifold (see lemma~\ref{lm-rect}).
We then prove that there is no measure concentrated on the set of points where the rank of the weak limit $\vec{\Phi}_\infty$ on each thick part and on each bubble is not equal to 2. Finally we use all the previous results to prove  a ``bubble tree convergence'' of the sequence $\vec{\Phi}_k$ on each thick part (lemma~\ref{lm-strong}) which gives in particular that the limiting rectifiable stationary varifold is integer. The last lemma, lemma~\ref{lm-integ}, establishes that the limiting map is a {\it conformal target harmonic map} on each thick part of the nodal surface and on each bubble.  

Theorem~\ref{th-main} can be used to prove various existence results of optimal varifolds realizing a min-max energy level. We first define the following notion.
\begin{Dfi}
\label{df-admissible}
A family of subsets ${\mathcal A}\subset{\mathcal P}({\mathcal M})$ of a Banach manifold ${\mathcal M}$ is called {\bf admissible family} if for every
homeomorphism $\Xi$ of ${\mathcal M}$ isotopic to the identity  we have
\[
\forall\, A\in{\mathcal A}\quad\quad\Xi(A)\in {\mathcal A}
\]
\hfill $\Box$
\end{Dfi}

\medskip

\noindent{\bf Example.} Consider ${\mathcal M}:=W^{2,{2p}}_{imm}(\Sigma,N^n)$  for some closed oriented surface $\Sigma$ and some closed sub-manifold $N^n$ of ${\R}^m$ and take for any $q\in {\N}$ and  $c\in \pi_q(\mbox{Imm}(\Sigma),N^n)$ then
the following family is admissible
\[
{\mathcal A}:=\lf\{
\ds \vec{\Phi}\in C^0(S^q,W^{2,2p}_{imm}(S^2,N^n))\ ;\ \quad\mbox{ s.t. } [\vec{\Phi}]=c\rg\}\quad.
\]
\hfill $\Box$

Our second main result is the following.
\begin{Th}
\label{th-min-max}
Let ${\mathcal A}$ be an admissible family in the space of $W^{2,2p}-$immersions into a closed sub-manifold of an euclidian space $N^n$.   Assume
\be
\label{0-01}
\inf_{A\in {\mathcal A}}\ \max_{\vec{\Phi}\in A}\mbox{Area}(\vec{\Phi})=\beta^0>0\quad,
\ee
then there exists a closed riemann surface $(S,h_0)$ with genus$(S)\le$genus$(\Sigma)$ and a  conformal  integer target harmonic
map $(\vec{\Phi}_\infty,N)$ from $S$ into $N^n$ such that
\[
\frac{1}{2}\int_SN\,|d\vec{\Phi}_\infty|_{h_0}^2 \ dvol_{h_0}=\beta^0\quad.
\]
\hfill $\Box$
\end{Th}
This general existence result has to be put in perspective with the previous min-max existence results partly discussed above either in {\it GMT} (see \cite{Pit}, \cite{SS}, \cite{CD}, \cite{MN1}, \cite{MN2}$\cdots$)  in {\it harmonic map theory} (see \cite{CM1},\cite{CM2},\cite{Xi1},\cite{Xi2}) or  using {\it level set-PDE} approaches (see \cite{HuTo}, \cite{ToWi}, \cite{Gua}, \cite{GaGua}, \cite{Ste-1}, \cite{Ste-2}).  Combined with the main regularity results in \cite{Ri5} and \cite{PiRi} theorem~\ref{th-min-max} implies in particular all known
results for the realization of arbitrary minmax by minimal surfaces. One technical advantage of the present work over the previous existing literature on minmax theory for surfaces in GMT or harmonic map theory, is that our proof of theorem~\ref{th-min-max}  does not require any ''replacement argument''.   The viscosity approach  gives moreover, without any additional work, an upper bound of the genus of the optimal surface. Such lower semicontinuity of the genus has been established in the GMT approach in  \cite{DP} in co-dimension 1 and was not given by the min-max procedure itself.
As in the geodesic case studied recently in \cite{MR} and where a passage to the limit in the second derivative is proved, the viscosity approach  gives under the multiplicity one assumption\footnote{The multiplicity one condition $N\equiv 1$ is expected to hold for finite index minmax problems in general. See the open problem in the first part of the introduction.} ($N= 1$ a.e on $S$) informations on the limiting index (see \cite{Riv-index}). This fact was left  open in  the {\it GMT} , the {\it harmonic map}  as well as in the {\it level set-PDE} approaches in it's full generality (see however partial important results in this direction for the PDE approach in \cite{MN3}).

\medskip

The second, and possibly main advantage, of the viscosity method resides in the fact that one can explore min-max within the space of \underbar{immersions} of fixed closed surfaces. The spaces Imm$(\Sigma,N^n)$ offers a  richer topology than the space of {\it integer rectifiable 2-cycles} ${\mathcal Z}_2(N^n)$ considered by Almgren whose homotopy type is  more coarse.  
The author has recently taken advantage of the full strength of theorem~\ref{th-min-max} for introducing new families of minmax problems at the level of immersions called {\it minmax hierarchies} (see \cite{Riv-index}).

\medskip

In order to simplify the presentation and in particular the computations of the Euler Lagrange equation to $A^\sigma$ we are presenting the proof of theorem~\ref{th-main}, in the special case $N^n=S^3$. There is however no
argument below which is specific to that case and the proof in the general case follows each step word for word of the $S^3$ case. Indeed, the almost conservation law in general target
manifold is perturbed by lower order terms (see for instance the explicit expression for $p=1$ and general target in \cite{MoR}). In arbitrary co-dimension each tensor has it's  counterpart 
which are possibly geometrically more involved but can be treated analytically identically as in the codimension 1 case. As soon as the strong $W^{1,2}-$bubble tree convergence is established, the passage to the limit in the non-linearity of the {\it harmonic map equation}\footnote{Recall that the mean-curvature vector of an immersion $\vec{\Phi}$ into a closed sub-manifold $N^n$ of an euclidian space is given by
\[
2\,\vec{H}_{\vec{\Phi}}=\Delta_{g_{\vec{\Phi}}}\vec{\Phi}+A(\vec{\Phi})(d\vec{\Phi},d\vec{\Phi})_{g_{\vP}}
\]
where $\Delta_{g_{\vec{\Phi}}}$ is the {\it negative Laplace Beltrami} operator with respect to the metric $g_{\vec{\Phi}}$. In conformal coordinates this becomes
\[
2\,\vec{H}_{\vec{\Phi}}=e^{-2\la}\,\lf[\Delta\vec{\Phi}+A(\vec{\Phi})(\nabla\vec{\Phi},\nabla\vec{\Phi})\rg]
\]
 where $e^\la:=|\p_{x_i}\vec{\Phi}|$.} is totally independent of the type of non linearity the target is producing. We keep from this non linearity, usually denoted $A(\vec{\Phi})(\nabla\vec{\Phi},\nabla\vec{\Phi})$ where $A$ is the second fundamental form of the target $N^n$, exclusively the quadratic dependence in the gradient. The conformal nature of the maps makes 
moreover the manipulation of the harmonic map equations straightforward independently of the existence or not of symmetries in the target. We took this point of view in order to ease a bit the reading of the proof.

\medskip

\noindent{\bf Acknowledgements} {\it The author is very grateful to Alexis Michelat and Alessandro Pigati as well as to the referees for their very careful reading of the paper and for having pointed out imprecisions and incorrectnesses in the preliminary version of the present work.}

\section{The viscous relaxation of the area for surfaces.}
\reset
\subsection{The Finsler manifold of immersions into the spheres with $L^q$ bounded second fundamental form.}

 For $k\in {\N}$ and $1\le q\le +\infty$, we recall the definition of $W^{k,q}$ Sobolev function on a closed smooth surface $\Sigma$ (i.e. $\Sigma$ is compact without boundary). To that aim we take
some reference smooth metric $g_0$ on $\Sigma$ and we set
\[
W^{k,q}(\Sigma,{\R}):=\lf\{f\ \mbox{ measurable } \; s.t \quad \nabla_{g_0}^kf\in L^q(\Sigma,g_0)\rg\}\quad,
\]
where $\nabla^k_{g_0}$ denotes the $k-$th iteration of the Levi-Civita connection associated to $\Sigma$. Since the surface is closed the space defined in this way is independent of $g_0$.
Let $N^n$ be a closed  $n-$dimensional sub-manifold of ${\R}^m$ with $3\le n\le m-1$ being arbitrary. The Space of $W^{k,q}$ into $N^n$ is defined as follows
\[
W^{k,q}(\Sigma, N^n):=\lf\{\vec{\Phi}\in W^{k,q}(\Sigma,{\R}^{m})\quad;\quad\vec{\Phi}\in N^n\ \mbox{ almost everywhere}\rg\}\quad.
\]
We have the following well known proposition
\begin{Prop}
\label{pr-prelim-1}
Assuming $kq>2$, the space  $W^{k,q}(\Sigma,N^n)$ defines a Banach Manifold. \hfill $\Box$
\end{Prop}
\noindent{\bf Proof of proposition~\ref{pr-prelim-1}.}

  This comes mainly from the fact that, under our assumptions,
\be
\label{prelim-II.1}
W^{k,q}(\Sigma,{\R}^m)\quad\hookrightarrow\quad C^0(\Sigma,{\R}^m)\quad.
\ee
The Banach manifold structure is then defined as follows. Choose $\delta>0$ such that each geodesic ball $B_\delta^{N^n}(z)$ for any $z\in N^n$ is strictly convex and the exponential map
\[
\exp_z\ :\ V_z\subset T_zN^n\ \longrightarrow\ B_\delta^{N^n}(z)
\]
realizes a $C^\infty$ diffeomorphism for some open neighborhood of the origin in $T_zN^n$ into the geodesic ball $B_\delta^{N^n}(z)$. Because of the embedding (\ref{prelim-II.1}) there exists $\ep_0>0$ such that
\[
\begin{array}{l}
\ds\forall\ \vec{u}\ ,\ \vec{v}\, \in W^{k,q}(\Sigma,{N}^n)\quad\|\vec{u}-\vec{v}\|_{W^{k,q}}<\ep_0\quad\\[3mm]
\ds\quad\quad\Longrightarrow\quad\|\mbox{dist}_N(\vec{u}(x),\vec{v}(x))\|_{L^\infty(\Sigma)}<\delta\quad.
\end{array}
\]
We equip now the  space $W^{k,q}(\Sigma,N^n)$ with the distance issued from the $W^{k,q}$ norm  and for any $\vec{u}\in {\mathcal M}=W^{k,q}(\Sigma,N^n)$ we denote by $B^{\mathcal M}_{\ep_0}(\vec{u})$ the open ball in ${\mathcal M}$ of center $\vec{u}$ and radius $\ep_0$. 

As a covering of ${\mathcal M}$ we take $(B^{\mathcal M}_{\ep_0}(\vec{u}))_{\vec{u}\in {\mathcal M}}$. We denote by 
\[
E^{\vec{u}}:=\Gamma_{W^{k,q}}\lf(\vec{u}^{-1}TN\rg):=\lf\{ \vec{w}\in W^{k,q}(\Sigma,{\R}^m)\ ;\ \vec{w}(x)\in T_{\vec{u}(x)}N^n\ \forall\ x\in\Sigma \rg\}
\] 
this is the Banach space of $W^{k,q}-$sections of the bundle $\vec{u}^{-1}TN$ and for any $\vec{u}\in {\mathcal M}$ and $\vec{v}\in B^{\mathcal M}_{\ep_0}(\vec{u})$ we define $\vec{w}^{\,\vec{u}}(\vec{v})$ to be the following element of $E^{\vec{u}}$
\[
\forall\ x\in \Sigma\quad\quad\vec{w}^{\,\vec{u}}(\vec{v})(x):=\exp^{-1}_{\vec{u}(x)}(\vec{v}(x))
\]
It is not difficult to see that
\[
\vec{w}^{\,\vec{v}}\circ\,(\vec{w}^{\,\vec{u}})^{-1}\ :\ \vec{w}^{\,u}\lf(B^{\mathcal M}_{\ep_0}(\vec{u})\cap B^{\mathcal M}_{\ep_0}(\vec{v})\rg)\ \longrightarrow\  \vec{w}^{\,v}\lf(B^{\mathcal M}_{\ep_0}(\vec{u})\cap B^{\mathcal M}_{\ep_0}(\vec{v})\rg)
\]
defines a $C^\infty$ diffeomorphism. \hfill $\Box$

\medskip

 For $p>1$ we define
\[
{\mathcal E}_{\Sigma,p}=W^{2,2p}_{imm}(\Sigma^2,{N}^n):=\lf\{ \vec{\Phi}\in W^{2,2p}(\Sigma^2,{N}^n)\ ;\  \mbox{rank}\,(d\Phi_x)=2\quad\forall x\in\Sigma^2 \rg\}\quad.
\]
The set $W^{2,2p}_{imm}(\Sigma^2,{N}^n)$ as an open subset of the normal Banach Manifold $W^{2,2p}(\Sigma^2,{N}^n)$ inherits a Banach Manifold structure.

\medskip

We equip now the space $W^{2,2p}_{imm}(\Sigma,N^n)$ with a {\it Finsler manifold structure} on it's  tangent bundle (see the definition of Banach bundle space and Tangent bundle to a Banach manifold in \cite{Lan}). For the convenience of the reader we recall the notion of Finsler structure.
\begin{Dfi}
\label{df-finsler}
Let ${\mathcal M}$ be a normal\footnote{The assumption to be {\it normal} is a relatively strong separation axiom which ensures  that the defined {\it Finsler structure} generates a distance which makes the topology of the Banach manifold metrizable (see \cite{Pa} pages 201-202). This assumption can be weakened to {\it regular} but not to {\it Hausdorff} only.} and let ${\mathcal V}$ be a Banach bundle space over ${\mathcal M}$. A {\bf Finsler structure} on ${\mathcal V}$ is a continuous function
\[
\|\cdot\|\ :\ {\mathcal V}\ \longrightarrow\ {\R}
\]
such that for any $x\in {\mathcal M}$
\[
\|\cdot\|_x:=\lf.\|\cdot\|\rg|_{\pi^{-1}(\{x\})}\quad\mbox{ is a norm on }{\mathcal V}_x\quad.
\]
Moreover for any local trivialization $\tau_i$ over $U_i$ and for any $x_0\in U_i$ we define on ${\mathcal V}_x$ the following norm
\[
\forall\ \vec{w}\in \pi^{-1}(\{x\})\quad\quad\|\vec{w}\|_{x_0}:=\|\tau_i^{-1}\lf(x_0,\rho(\tau_i(\vec{w}))\rg)\|_{x_0}\quad,
\]
where $\rho$ is the canonical projection $\rho\  :\ U_i\times E\rightarrow E$ and there exists $C_{x_0}>1$ such that
\[
\forall\ x\in U_i\quad\quad C_{x_0}^{-1}\ \|\cdot\|_x\le\|\cdot\|_{x_0}\le C_{x_0}\ \|\cdot\|_x\quad.
\] 
In a $C^q$ Banach bundle, the Finsler structure is said to be $C^l$ for $l\le q$ if, in local charts, the dependence of $\|\cdot\|_{x}$ is $C^l$ with respect to $x$.
\hfill $\Box$
\end{Dfi}
\begin{Dfi}
\label{df-finsler-manifold} Let ${\mathcal M}$ be a normal $C^p$ Banach manifold. $T{\mathcal M}$ equipped with a Finsler structure is called a {\bf Finsler Manifold}.
\hfill $\Box$
\end{Dfi}
\begin{Rm}
\label{rm-prelim-II.2}
A Finsler structure on $T{\mathcal M}$ defines in a canonical way a dual Finsler structure on $T^\ast{\mathcal M}$.\hfill $\Box$
\end{Rm}
  The tangent
space to ${\mathcal E}_{\Sigma,p}$ at a point $\vec{\Phi}$ is the space $\Gamma_{W^{2,2p}}(\vec{\Phi}^{\,-1} TN^n)$ of $W^{2,2p}-$sections of the bundle $\vec{\Phi}^{\,-1} TN^n$, i.e.
\[
T_{\vec{\Phi}}{\mathcal E}_{\Sigma,p}=\lf\{\vec{w}\in W^{2,2p}(\Sigma^2,{\R}^m)\ ;\ \vec{w}(x)\in T_{\vec{\Phi}(x)}N^n\quad\forall\, x\in\, \Sigma^2\rg\}\quad.
\]
We equip $T_{\vec{\Phi}}{\mathcal E}_{\Sigma,p}$ with the following norm
\[
 \|\vec{v}\|_{\vec{\Phi}}:=\lf[ \int_{\Sigma}\lf[|\nabla^2\vec{v}|^2_{g_{\vec{\Phi}}}+|\nabla\vec{v}|^2_{g_{\vec{\Phi}}}+|\vec{v}|^2\rg]^{p}\ dvol_{g_{\vec{\Phi}}}\rg]^{1/2p}+\|\,|\nabla\vec{v}|_{g_{\vec{\Phi}}}\,\|_{L^\infty(\Sigma)}\quad,
 \]
 where we keep denoting, for any $j\in {\N}$, $\nabla$ to be the connection on $(T^\ast\Sigma)^{\otimes^j}\otimes\vec{\Phi}^{-1}TN$ over $\Sigma$ defined by $\nabla:=\nabla^{g_{\vec{\Phi}}}\otimes\vec{\Phi}^\ast \nabla^h$ and $\nabla^{g_{\vec{\Phi}}}$ is the Levi Civita connection on $(\Sigma,g_{\vec{\Phi}})$ and $\nabla^h$ is the Levi-Civita connection on $N^n$. 
 
 We check for instance that $\nabla\vec{v}$ , resp. $\nabla^2\vec{v}$ defines a $C^0$, resp. $L^{2p}$, section of $(T^\ast\Sigma)\otimes \vec{\Phi}^{-1}TN$, resp. $(T^\ast\Sigma)^2\otimes \vec{\Phi}^{-1}TN$.
 
The fact that we are adding to the $W^{2,2p}$ norm of $\vec{v}$ 
 with respect to $g_{\vec{\Phi}}$ the $L^\infty$ norm of $|\nabla\vec{v}|_{g_{\vec{\Phi}}}$ could look redundant since $W^{2,2p}$ embeds in $W^{1,\infty}$. We are doing it in order to ease the proof of the completeness of the Finsler Space equipped with the Palais distance below.

 Observe that, using Sobolev embedding and in particular due to the fact $W^{2,q}(\Sigma,{\R}^m)\hookrightarrow C^1(\Sigma,{\R}^m)$ for $q>2$, the norm $\|\cdot\|_{\vec{\Phi}}$ as a function on the Banach tangent bundle $T{\mathcal E}_{\Sigma,p}$ is obviously continuous.
 
 \begin{Prop}
 \label{pr-prelim-II.2}
 The norms $\|\cdot\|_{\vec{\Phi}}$  defines a $C^2-$Finsler structure on the space ${\mathcal E}_{\Sigma,p}$.\hfill $\Box$
 \end{Prop}
\noindent{\bf Proof of proposition~\ref{pr-prelim-II.2}.}
We introduce the following trivialization of the Banach bundle. For any $\vec{\Phi}\in {\mathcal E}_{\Sigma,p}$ we denote $P_{\vec{\Phi}(x)}$ the orthonormal
projection in ${\R}^m$ onto the $n-$dimensional vector subspace of ${\R}^m$ given by $T_{\vec{\Phi}(x)}N^n$ and for any $\vec{\xi}$ in the ball $B^{{\mathcal E}_{\Sigma,p}}_{\ep_1}(\vec{\Phi})$ for some $\ep_1>0$ and any $\vec{v}\in T_{\vec{\xi}}{\mathcal E}_{\Sigma,p}=\Gamma_{W^{2,2p}}(\vec{\xi}^{-1}TN)$ we assign the map $\vec{w}(x):=P_{\vec{\Phi}(x)}\vec{v}(x)$. It is straightforward to check
that for $\ep_1>0$ chosen small enough the map which to $\vec{v}$ assigns $\vec{w}$ is an isomorphism from $T_{\vec{\xi} }{{\mathcal E}_{\Sigma,p}}$ into $T_{\vec{\Phi}}{\mathcal E}_{\Sigma,p}$ and that there exists
$k_{\vec{\Phi}}>1$ such that $\forall \vec{v}\in TB^{{\mathcal E}_{\Sigma,p}}_{\ep_1}(\vec{\Phi})$
\[
k_{\vec{\Phi}}^{-1}\,\|\vec{v}\|_{\vec{\xi}}\le \|\vec{w}\|_{\vec{\Phi}}\le k_{\vec{\Phi}}\,\|\vec{v}\|_{\vec{\xi}}\quad.
\]
The $C^2-$dependence of $\|\cdot\|_{\vec{\xi}}$  with respect to $\vec{\xi}$ in the chart above is left to the reader. This concludes the proof of proposition~\ref{pr-prelim-II.2}.\hfill $\Box$

\medskip

\subsection{Palais deformation theory applied to the space of  $W^{2,2p}-$immersions.}

\begin{Th}
\label{th-distance-palais}{\bf [Palais 1970]} Let $({\mathcal M},\|\cdot\|)$ be a Finsler Manifold. Define on ${\mathcal M}\times{\mathcal M}$ 
\[
d(p,q):=\inf_{{\omega}\in \Omega_{p,q}}\int_0^1\lf\|\frac{d\omega}{dt}\rg\|_{\omega(t)}\ dt\quad,
\]
where
\[
\Omega_{p,q}:=\lf\{\om\in C^1([0,1],{\mathcal M})\ ;\ \om(0)=p\quad\om(1)=q\rg\}\quad.
\]
Then $d$ defines a distance on ${\mathcal M}$ and $({\mathcal M},d)$ defines the same topology as the one of the Banach Manifold. $d$ is called {\bf Palais distance}
of the Finsler manifold $({\mathcal M},\|\cdot\|)$.
\hfill $\Box$
\end{Th}
Contrary to the first appearance the non degeneracy of $d$ is not straightforward and requires a proof (see \cite{Pa}). This last result combined with the famous result of Stones (see \cite{Sto}) on the paracompactness of metric spaces gives the following corollary.
\begin{Co}
\label{co-finsler-para} Let $({\mathcal M},\| \cdot\|)$ be a Finsler Manifold then ${\mathcal M}$ is paracompact.\hfill $\Box$
\end{Co} 
The following result\footnote{As a matter of fact the proof of the completeness with respect to the Palais distance is  skipped in various applications
of Palais deformation theory in the literature.} is going to play a central role in adapting Palais deformation theory to our framework of $W^{2,2p}-$immersions.
\begin{Prop}
\label{th-complete}
Let $p>1$ and ${\mathcal M}:={\mathcal E}_{\Sigma,p}$ be the space of $W^{2,2p}-$immersions of a closed oriented surface $\Sigma$ into a closed
sub-manifold $N^n$ of ${\R}^m$
\[
{\mathcal E}_{\Sigma,p}=W^{2,2p}_{imm}(\Sigma^2,{N}^n):=\lf\{ \vec{\Phi}\in W^{2,2p}(\Sigma^2,{N}^n)\ ;\  \mbox{rank}\,(d\Phi_x)=2\quad\forall x\in\Sigma^2 \rg\}\quad.
\]
The Finsler Manifold given by the structure
\[
 \|\vec{v}\|_{\vec{\Phi}}:=\lf[ \int_{\Sigma}\lf[|\nabla^2\vec{v}|^2_{g_{\vec{\Phi}}}+|\nabla\vec{v}|^2_{g_{\vec{\Phi}}}+|\vec{v}|^2\rg]^{p}\ dvol_{g_{\vec{\Phi}}}\rg]^{1/2p}+\|\,|\nabla\vec{v}|_{g_{\vec{\Phi}}}\,\|_{L^\infty(\Sigma)}
 \]
is \underbar{complete} for the Palais distance.
\hfill $\Box$
\end{Prop}
\noindent{\bf Proof of proposition~\ref{th-complete}.}
For any $\vec{\Phi}\in {\mathcal M}$ and $\vec{v}\in T_{\vec{\Phi}}{\mathcal M}$ we introduce the tensor in $(T^\ast\Sigma)^{\otimes^2}$ given in coordinates by
\[
\begin{array}{l}
\ds\nabla\vec{v}\,\dot{\otimes}\,d\vec{\Phi}+d\vec{\Phi}\,\dot{\otimes}\,\nabla\vec{v}=\sum_{i,j=1}^2\lf[\nabla_{\p_{x_i}}\vec{v}\cdot{\p_{x_j}}\vec{\Phi}+ \p_{x_i}\vec{\Phi}\cdot\nabla_{\p_{x_j}}\vec{v} \rg]\ dx_i\otimes dx_j\\[5mm]
\ds\quad\quad\quad=\sum_{i,j=1}^2\lf[\nabla^h_{\p_{x_i}\vec{\Phi}}\vec{v}\cdot{\p_{x_j}}\vec{\Phi}+ \p_{x_i}\vec{\Phi}\cdot\nabla^h_{\p_{x_j}\vec{\Phi}}\vec{v} \rg]\ dx_i\otimes dx_j
\end{array}
\]
where $\cdot$ denotes the scalar product in ${\R}^m$. Observe that we have
\[
\lf|\nabla\vec{v}\,\dot{\otimes}\,d\vec{\Phi}+d\vec{\Phi}\,\dot{\otimes}\,\nabla\vec{v}\rg|_{g_{\vec{\Phi}}}\le 2\ |\nabla\vec{v}|_{g_{\vec{\Phi}}}\quad.
\]
Hence, taking a $C^1$ path  $\vec{\Phi}_s$ in $\mathcal M$ one has for $\vec{v}:=\p_s\vec{\Phi}$
\be
\label{x-I.1}
\begin{array}{l}
\ds\| |d\vec{v}\dot{\otimes}d\vec{\Phi}+d\vec{\Phi}\dot{\otimes}d\vec{v}|^2_{g_{\vec{\Phi}}} \|_{L^\infty(\Sigma)}=\lf\|\sum_{i,j,k,l=1}^2\ g^{ij}_{\vec{\Phi}}\ g^{kl}_{\vec{\Phi}}\ \p_s(g_{\vec{\Phi}})_{ik}\ \p_s(g_{\vec{\Phi}})_{jl}\rg\|_{L^\infty(\Sigma)}\\[5mm]
\ds\quad=\lf\|\ |\p_s(g_{ij} dx_i\otimes dx_j)|^2_{g_{\vec{\Phi}}}\rg\|_{L^\infty(\Sigma)}=\lf\| \ |\p_sg_{\vec{\Phi}}|^2_{g_{\vec{\Phi}}}\rg\|_{L^\infty(\Sigma)}\quad.
\end{array}
\ee
Hence
\be
\label{x-I.2}
\int_0^1\lf\| \ |\p_sg_{\vec{\Phi}}|_{g_{\vec{\Phi}}}\rg\|_{L^\infty(\Sigma)}\ ds\le2\ \int_0^1\|\p_s\vec{\Phi}\|_{\vec{\Phi}_s}\, ds\quad.
\ee
We now use the following lemma
\begin{Lm}
\label{lm-I.1}
Let $M_s$ be a $C^1$ path into the space of positive $n$ by $n$ symmetric matrix then the following inequality holds
\[
\mbox{Tr}\,(M^{-2}(\p_sM)^2)= \|\p_s\log M\|^2=\mbox{Tr}\,((\p_s\log M)^2)\quad
\]
\hfill $\Box$
\end{Lm}
\noindent{\bf Proof of lemma~\ref{lm-I.1}.} We write $M=\exp A$ and we observe that
\[
\mbox{Tr}\,(\exp(-2 A)(\p_s\exp A)^2)=\mbox{Tr}\,(\p_sA)^2\quad.
\]
Then the lemma follows.\hfill $\Box$

Combining the previous lemma with (\ref{x-I.1}) and (\ref{x-I.2}) we obtain in a given chart
\be
\label{x-I.2-a}
\int_0^1\|\p_s\log (g_{ij})\|\ ds=\int_0^1 \sqrt{\mbox{Tr}\,((\p_s\log g_{ij})^2)}\, ds\le 2\ \int_0^1\|\p_s\vec{\Phi}\|_{\vec{\Phi}_s}\, ds\quad.
\ee
This implies that in the given chart the $\log$ of the matrix $(g_{ij}(s))$ is uniformly bounded for $s\in [0,1]$ and hence $\vec{\Phi}_1$ is an immersion. It remains to show
that it has a controlled $W^{2,q}$ norm. We  denote
\[
\mbox{Hess}_p(\vec{\Phi}):=\int_\Sigma[1+|\nabla d\vec{\Phi}|_{g_{\vec{\Phi}}}^2]^p\ dvol_{g_{\vec{\Phi}}}
\] 
and we compute
\be
\label{x-I.3}
\begin{array}{l}
\ds\frac{d}{ds}(\mbox{Hess}_p(\vec{\Phi}))=p\,\int_\Sigma\p_s|\nabla d\vec{\Phi}|_{g_{\vec{\Phi}}}^2\, [1+|\nabla d\vec{\Phi}|_{g_{\vec{\Phi}}}^2]^{p-1}\ dvol_{g_{\vec{\Phi}}}\\[5mm]
\ds\quad\quad+\int_\Sigma [1+|\nabla d\vec{\Phi}|_{g_{\vec{\Phi}}}^2]^{p}\ \p_s(dvol_{g_{\vec{\Phi}}})\quad.
\end{array}
\ee
Classical computations give
\[
\p_s(dvol_{g_{\vec{\Phi}}})=\lf<\nabla\p_s\vec{\Phi},d\vec{\Phi}\rg>_{g_{\vec{\Phi}}}\ dvol_{g_{\vec{\Phi}}}\quad.
\]
So we have
\be
\label{x-I.4}
\begin{array}{l}
\ds\lf|\int_\Sigma [1+|\nabla d\vec{\Phi}|_{g_{\vec{\Phi}}}^2]^{p}\ \p_s(dvol_{g_{\vec{\Phi}}})\rg|\le \|\,|\nabla\p_s\vec{\Phi}|_{g_{\vec{\Phi}}}\,\|_{L^\infty(\Sigma)}\int_\Sigma [1+|\nabla d\vec{\Phi}|_{g_{\vec{\Phi}}}^2]^{p}\ dvol_{g_{\vec{\Phi}}}\\[5mm]
\ds \quad\quad\quad\quad\quad\le \|\p_s\vec{\Phi}\|_{\vec{\Phi}}\ \int_\Sigma [1+|\nabla d\vec{\Phi}|_{g_{\vec{\Phi}}}^2]^{p}\ dvol_{g_{\vec{\Phi}}}\quad.
\end{array}
\ee
In local charts we have
\[
|\nabla d\vec{\Phi}|_{g_{\vec{\Phi}}}^2=\sum_{i,j,k,l=1}^2g^{ij}_{\vec{\Phi}}\, g^{kl}_{\vec{\Phi}}\, \lf<\nabla^h_{\p_{x_i}\vec{\Phi}}\p_{x_k}\vec{\Phi},\nabla^h_{\p_{x_j}\vec{\Phi}}\p_{x_l}\vec{\Phi}\rg>_h\quad.
\]
Thus in bounding $\int_\Sigma\p_s|\nabla d\vec{\Phi}|_{g_{\vec{\Phi}}}^2\, [1+|\nabla d\vec{\Phi}|_{g_{\vec{\Phi}}}^2]^{p-1}\ dvol_{g_{\vec{\Phi}}}$ we first have to control terms of the form
\be
\label{x-I.5}
\lf|\int_\Sigma\sum_{i,j,k,l=1}^2\p_sg^{ij}_{\vec{\Phi}}\, g^{kl}_{\vec{\Phi}}\, \lf<\nabla^h_{\p_{x_i}\vec{\Phi}}\p_{x_k}\vec{\Phi},\nabla^h_{\p_{x_j}\vec{\Phi}}\p_{x_l}\vec{\Phi}\rg>_h\ [1+|\nabla d\vec{\Phi}|_{g_{\vec{\Phi}}}^2]^{p-1}\  dvol_{g_{\vec{\Phi}}}\rg|\quad.
\ee
We write
\[
\begin{array}{l}
\ds\sum_{i,j,k,l=1}^2\p_sg^{ij}_{\vec{\Phi}}\, g^{kl}_{\vec{\Phi}}\, \lf<\nabla^h_{\p_{x_i}\vec{\Phi}}\p_{x_k}\vec{\Phi},\nabla^h_{\p_{x_j}\vec{\Phi}}\p_{x_l}\vec{\Phi}\rg>_h\\[5mm]
\ds\quad=\sum_{i,j,k,l,t,r=1}^2\p_sg^{ij}_{\vec{\Phi}}\, g_{jt}\, g^{tr} g^{kl}_{\vec{\Phi}}\, \lf<\nabla^h_{\p_{x_i}\vec{\Phi}}\p_{x_k}\vec{\Phi},\nabla^h_{\p_{x_j}\vec{\Phi}}\p_{x_l}\vec{\Phi}\rg>_h\\[5mm]
\ds\quad=-\sum_{i,j,k,l,=1}^2\, \lf(\sum_{t,r=1}^2\p_sg_{jt}\, g^{tr}\rg)\ g^{ij}_{\vec{\Phi}}\, g^{kl}_{\vec{\Phi}}\, \lf<\nabla^h_{\p_{x_i}\vec{\Phi}}\p_{x_k}\vec{\Phi},\nabla^h_{\p_{x_j}\vec{\Phi}}\p_{x_l}\vec{\Phi}\rg>_h\quad.
\end{array}
\]
Hence, using (\ref{x-I.1}),
\be
\label{x-I.6}
\begin{array}{l}
\ds\lf|\int_\Sigma\sum_{i,j,k,l=1}^2\p_sg^{ij}_{\vec{\Phi}}\, g^{kl}_{\vec{\Phi}}\, \lf<\nabla^h_{\p_{x_i}\vec{\Phi}}\p_{x_k}\vec{\Phi},\nabla^h_{\p_{x_j}\vec{\Phi}}\p_{x_l}\vec{\Phi}\rg>_h\ [1+|\nabla d\vec{\Phi}|_{g_{\vec{\Phi}}}^2]^{p-1}\  dvol_{g_{\vec{\Phi}}}\rg|\\[5mm]
\ds\quad\le\ \|\,|\p_s g_{\vec{\Phi}}|_{g_{\vec{\Phi}}}\|_{L^\infty(\Sigma)}\ \int_{\Sigma}\ [1+|\nabla d\vec{\Phi}|_{g_{\vec{\Phi}}}^2]^{p}\ dvol_{g_{\vec{\Phi}}}\\[5mm]
\ds\quad\le\ 2\,\|\p_s\vec{\Phi}\|_{\vec{\Phi}_s}\ \int_{\Sigma}\ [1+|\nabla d\vec{\Phi}|_{g_{\vec{\Phi}}}^2]^{p}\ dvol_{g_{\vec{\Phi}}}\quad.
\end{array}
\ee
We have also
\[
\begin{array}{l}
\ds\p_s\lf<\nabla^h_{\p_{x_i}\vec{\Phi}}\p_{x_k}\vec{\Phi},\nabla^h_{\p_{x_j}\vec{\Phi}}\p_{x_l}\vec{\Phi}\rg>_h\\[5mm]
\ds\quad=\lf<\nabla^h_{\p_s\vec{\Phi}}\lf(\nabla^h_{\p_{x_i}\vec{\Phi}}\p_{x_k}\vec{\Phi}\rg),\nabla^h_{\p_{x_j}\vec{\Phi}}\p_{x_l}\vec{\Phi}\rg>_h+\lf<\nabla^h_{\p_{x_i}\vec{\Phi}}\p_{x_k}\vec{\Phi},\nabla^h_{\p_s\vec{\Phi}}\lf(\nabla^h_{\p_{x_j}\vec{\Phi}}\p_{x_l}\vec{\Phi}\rg)\rg>_h\quad.
\end{array}
\]
By definition we have
\[
\nabla^h_{\p_s\vec{\Phi}}\lf(\nabla^h_{\p_{x_i}\vec{\Phi}}\p_{x_k}\vec{\Phi}\rg)=\nabla^h_{\p_{x_i}\vec{\Phi}}\lf(\nabla^h_{\p_{s}\vec{\Phi}}\p_{x_k}\vec{\Phi}\rg)+R^h(\p_{x_i}\vec{\Phi},\p_s\vec{\Phi})\p_{x_k}\vec{\Phi}\quad,
\]
where we have used the fact that $[\p_s\vec{\Phi},\p_{x_i}\vec{\Phi}]=\vec{\Phi}_\ast[\p_s,\p_{x_i}]=0$. Using also that $[\p_s\vec{\Phi},\p_{x_k}\vec{\Phi}]=0$, since $\nabla^h$ is torsion free,
we have finally
\be
\label{x-I.7}
\nabla^h_{\p_s\vec{\Phi}}\lf(\nabla^h_{\p_{x_i}\vec{\Phi}}\p_{x_k}\vec{\Phi}\rg)=\nabla^h_{\p_{x_i}\vec{\Phi}}\lf(\nabla^h_{\p_{x_k}\vec{\Phi}}\p_{s}\vec{\Phi}\rg)+R^h(\p_{x_i}\vec{\Phi},\p_s\vec{\Phi})\p_{x_k}\vec{\Phi}\quad,
\ee
where $R^h$ is the Riemann tensor associated to the Levi-Civita connection $\nabla^h$. We have
\be
\label{x-I.8}
\nabla^h_{\p_{x_i}\vec{\Phi}}\lf(\nabla^h_{\p_{x_k}\vec{\Phi}}\p_{s}\vec{\Phi}\rg)=(\nabla^h)^2_{\p_{x_i}\vec{\Phi}\p_{x_k}\vec{\Phi}}\p_s\vec{\Phi}+\nabla^h_{\nabla^h_{\p_{x_i}\vec{\Phi}}\p_{x_k}\vec{\Phi}}\p_s\vec{\Phi}
\ee
Hence
\be
\label{x-I.9}
\begin{array}{l}
\lf<\nabla^h_{\p_s\vec{\Phi}}\lf(\nabla^h_{\p_{x_i}\vec{\Phi}}\p_{x_k}\vec{\Phi}\rg),\nabla^h_{\p_{x_j}\vec{\Phi}}\p_{x_l}\vec{\Phi}\rg>_h=\lf<(\nabla^h)^2_{\p_{x_i}\vec{\Phi}\p_{x_k}\vec{\Phi}}\p_s\vec{\Phi},\nabla^h_{\p_{x_j}\vec{\Phi}}\p_{x_l}\vec{\Phi}\rg>_h\\[5mm]
\ds\quad+\lf<\nabla^h_{\nabla^h_{\p_{x_i}\vec{\Phi}}\p_{x_k}\vec{\Phi}}\p_s\vec{\Phi},\nabla^h_{\p_{x_j}\vec{\Phi}}\p_{x_l}\vec{\Phi}\rg>_h+ \lf< R^h(\p_{x_i}\vec{\Phi},\p_s\vec{\Phi})\p_{x_k}\vec{\Phi}
 ,\nabla^h_{\p_{x_j}\vec{\Phi}}\p_{x_l}\vec{\Phi}\rg>_h\quad.
\end{array}
\ee
Combining all the previous and observing that
\be
\label{x-I.9-a}
\lf|\sum_{i,j,k,l=1}^2 g^{ij}_{\vec{\Phi}}\, g^{kl}_{\vec{\Phi}}\, \lf<\nabla^h_{\nabla^h_{\p_{x_i}\vec{\Phi}}\p_{x_k}\vec{\Phi}}\p_s\vec{\Phi},\nabla^h_{\p_{x_j}\vec{\Phi}}\p_{x_l}\vec{\Phi}\rg>_h \rg|\le C\,|\nabla\p_s\vec{\Phi}|_{g_{\vec{\Phi}}}\ |\nabla d\vec{\Phi}|^2_{g_{\vec{\Phi}}}\ 
\ee
gives then
\be
\label{x-I.10}
\begin{array}{l}
\ds\lf|\int_\Sigma\sum_{i,j,k,l=1}^2 g^{ij}_{\vec{\Phi}}\, g^{kl}_{\vec{\Phi}}\, \p_s\lf<\nabla^h_{\p_{x_i}\vec{\Phi}}\p_{x_k}\vec{\Phi},\nabla^h_{\p_{x_j}\vec{\Phi}}\p_{x_l}\vec{\Phi}\rg>_h\ \ [1+|\nabla d\vec{\Phi}|_{g_{\vec{\Phi}}}^2]^{p-1}\  dvol_{g_{\vec{\Phi}}}\rg|\\[5mm]
\ds\quad\le C\, \int_\Sigma\ \lf|\lf<\nabla^2\p_s\vec{\Phi},\nabla d\vec{\Phi}\rg>_{g_{\vec{\Phi}}}\rg|\ [1+|\nabla d\vec{\Phi}|_{g_{\vec{\Phi}}}^2]^{p-1}\ dvol_{g_{\vec{\Phi}}}\\[5mm]
\ds\quad\quad+\, C\,\int_\Sigma\ |\nabla\p_s\vec{\Phi}|_{g_{\vec{\Phi}}}\ |\nabla d\vec{\Phi}|^2_{g_{\vec{\Phi}}}\ [1+|\nabla d\vec{\Phi}|_{g_{\vec{\Phi}}}^2]^{p-1}\ dvol_{g_{\vec{\Phi}}}\\[5mm]
\ds\quad\quad+\, C\, \|R^h\|_{L^\infty(N^n)}\ \int_\Sigma |\p_s\vec{\Phi}|_h\ |\nabla d\vec{\Phi}|_{g_{\vec{\Phi}}}\ [1+|\nabla d\vec{\Phi}|_{g_{\vec{\Phi}}}^2]^{p-1}\ dvol_{g_{\vec{\Phi}}}\quad.
\end{array}
\ee
Combining all the above we finally obtain that
\be
\label{x-I.11}
\lf|\p_s\mbox{Hess}_p(\vec{\Phi})\rg|\le C\, \|\p_s\vec{\Phi}\|_{\vec{\Phi}}\ \lf[\mbox{Hess}_p(\vec{\Phi})+\mbox{Hess}_p(\vec{\Phi})^{1-1/2p}\rg]\quad.
\ee
Combining (\ref{x-I.2-a}) and (\ref{x-I.11}) we deduce using Gromwall lemma that if we take a $C^1$ path from $[0,1)$ into ${\mathcal E}_{\Sigma,p}$ with finite length for the Palais distance $d$, the limiting map $\vec{\Phi}_1$
is still a $W^{2,2p}-$immersion of $\Sigma$ into $N^n$, which proves the completeness of $({\mathcal E}_{\Sigma,p},d)$.\hfill $\Box$

\medskip
The following definition is central in Palais deformation theory.
\begin{Dfi}
\label{df-palais-smale} Let $E$ be a $C^1$ function on a Finsler manifold $({\mathcal M},\|\cdot\|)$ and $\beta\in E({\mathcal M})$. On says that $E$ fulfills the {\bf Palais-Smale condition}
at the level $\beta$ if for any sequence $u_n$ staisfying
\[
E(u_n)\longrightarrow \beta\quad\mbox{ and }\quad\|D E_{u_n}\|_{u_n}\longrightarrow 0\quad,
\]
then there exists a subsequence $u_{n'}$ and $u_\infty\in {\mathcal M}$ such that
\[
d(u_{n'},u_\infty)\longrightarrow 0\quad.
\]
and hence $E(u_\infty)=\beta$ and $DE_{u_\infty}=0$. \hfill $\Box$
\end{Dfi}
The following result is the Palais Smale condition for the functional 
$$
A_p^\sigma(\vec{\Phi}):=\mbox{Area}(\vec{\Phi})+\sigma^2\int_{\Sigma}\lf[1+|\vec{\mathbb I}_{\vec{\Phi}}|^2\rg]^{p}\ dvol_{g_{\vec{\Phi}}}\quad.
$$
\begin{Th}
\label{th-palais}
Let $p>1$ and $\vec{\Phi}_k$ such that
\[
\limsup_{k\rightarrow +\infty}A_p^\sigma(\vec{\Phi}_k)<+\infty\quad,
\]
and satisfying
\be
\label{V.1-a}
\lim_{k\rightarrow+\infty}\sup_{\|\vec{w}\|_{\vec{\Phi}_k}\le 1}D A_p^\sigma(\vec{\Phi}_k)\cdot \vec{w}=0\quad.
\ee
Then, modulo extraction of a subsequence, there exists a sequence of $W^{2,2p}-$diffeomorphisms $\Psi_k$ such that $\vec{\Phi}_k\circ\Psi_k$ converges strongly in ${\mathcal E}_{\Sigma,p}$ for the Palais distance to a critical point
of $A_p^\sigma$. Moreover, if one assume that $\vec{\Phi}_k$ stays inside a fixed ball of the Palais distance one can take $\Psi_k(x)=x$.  \hfill $\Box$
\end{Th}
\begin{Rm}
\label{rm-palais}
The first part of this theorem has been proved in \cite{KLL}  (theorem 5.1) in the flat framework which does not differ much from our case of $W^{2,2p}-$immersions into $N^n$. See also
\cite{RiB} for a proof making use of the underlying conservation laws. The second part is a direct consequence of the proof of proposition~\ref{th-complete} above and is being used below since in the main Palais theorem~\ref{th-II-main} the flow issued by the pseudo-gradient maintains the image at a finite Palais distance. \hfill $\Box$
\end{Rm}

\medskip

\begin{Dfi}
\label{df-admissible}
A family of subsets ${\mathcal A}\subset{\mathcal P}({\mathcal M})$ of a Banach manifold ${\mathcal M}$ is called {\bf admissible family} if for every
homeomorphism $\Xi$ of ${\mathcal M}$ isotopic to the identity  we have
\[
\forall\, A\in{\mathcal A}\quad\quad\Xi(A)\in {\mathcal A}\quad.
\]
\hfill $\Box$
\end{Dfi}

\medskip

\noindent{\bf Example.} Consider ${\mathcal M}:=W^{2,q}_{imm}(S^2,{\R}^3)$ and take\footnote{It is proved in \cite{Sma} and \cite{Hir} that 
\[
\mbox{Imm}(S^2,{\R}^3)\simeq_{homotop.} SO(3)\times\Om^2(SO(3))
\]} $c\in \pi_1(\mbox{Imm}(S^2,{\R}^3))={\Z}_2\times{\Z}$ then
the following family is admissible
\[
{\mathcal A}:=\lf\{
\ds \vec{\Phi}\in C^0([0,1],W^{2,q}_{imm}(S^2,{\R}^3))\ ;\ \vec{\Phi}(0,\cdot)=\vec{\Phi}(1,\cdot)\quad\mbox{ and } [\vec{\Phi}]=c\rg\}\quad.
\]
\hfill $\Box$

\medskip

We recall the main theorem of Palais deformation theory.
\begin{Th}
\label{th-II-main}{\bf[Palais 1970]}
Let $({\mathcal M},\|\cdot\|)$ be a Banach manifold together with a $C^{1,1}-$Finsler structure. Assume ${\mathcal M}$ is \underbar{complete} for the induced Palais distance $d$ and
let $E\in C^1({\mathcal M})$ satisfying the Palais-Smale condition $(PS)_\beta$ for the level set $\beta$. Let ${\mathcal A}$ be an admissible family in ${\mathcal P}({\mathcal M})$
such that
\[
\inf_{A\in{\mathcal A}}\ \sup_{u\in A} E(u)=\beta\quad,
\]
then there exists $u\in {\mathcal M}$ satisfying
\be
\label{II.01}
\lf\{
\begin{array}{l}
D E_u=0\\[5mm]
E(u)=\beta
\end{array}
\rg.
\ee
\hfill $\Box$
\end{Th}

\subsection{Struwe's monotonicity trick.}

Because of theorem~\ref{th-palais}, theorem~\ref{th-II-main} can be applied to each of the lagrangian $A_p^\sigma$ for any admissible family $\mathcal A$ of ${\mathcal E}_{\Sigma,p}$ satisfying 
\be
\label{0-01-rep}
\inf_{A\in {\mathcal A}}\ \max_{\vec{\Phi}\in A}\mbox{Area}(\vec{\Phi})=\beta^0>0\quad.
\ee
However, beside the difficulty of establishing a convergence of any nature to the corresponding sequence of critical points $\vec{\Phi}_\sigma$ given by theorem~\ref{th-II-main}, although it is clear that
\[
\lim_{\sigma\rightarrow 0}\inf_{A\in {\mathcal A}}\ \max_{\vec{\Phi}\in A}A^\sigma_p(\vec{\Phi}_\sigma)=\beta^0\quad,
\]
nothing excludes {\it a-priori} that
 \[
\lim_{\sigma\rightarrow 0}\inf_{A\in {\mathcal A}}\ \max_{\vec{\Phi}\in A}\mbox{Area}(\vec{\Phi}_\sigma)<\beta^0\quad,
\]
and it could be that the smoothing part of the lagrangian $\sigma^2\int_{\Sigma}\lf[1+|\vec{\mathbb I}_{\vec{\Phi}_\sigma}|^2\rg]^{p}\ dvol_{g_{\vec{\Phi}_\sigma}}$ does not go to zero. In order to prevent this unpleasant situation where the smoothed min-max
procedure is not approximating properly the limiting min-max procedure, M.Struwe invented a technic - called sometimes ``Struwe's monotonicity trick'' - consisting in localizing the action 
of the pseudo-gradient close to the level set Area$(\vec{\Phi})=\beta_0$ exclusively . Precisely we have the following result.
 \begin{Th}
\label{th-III.3-reg}
Let $({\mathcal M},\|\cdot\|)$ be a complete Finsler manifold. Let $E^\sigma$ be a family of $C^1$ functions for $\sigma\in [0,1]$ on ${\mathcal M}$ such that
for every $\vec{\gamma}\in {\mathcal M}$ 
\be
\label{III.4-reg}
\sigma\longrightarrow E^\sigma(\vec{\gamma})\quad\mbox{ and }\quad\sigma\longrightarrow \partial_\sigma E^\sigma(\vec{\gamma})\quad,
\ee
are increasing and continuous functions with respect to $\sigma$.  Assume moreover that 
\be
\label{III.5-reg}
\|DE^\sigma_{\vec{\gamma}}-DE^{\tau}_{\vec{\gamma}}\|_{\vec{\gamma}}\le C(\sigma)\ \delta(|\sigma-\tau|)\ f(  E^\sigma(\vec{\gamma}))\quad,
\ee
where $C(\sigma)\in L^\infty_{loc}((0,1))$, $\delta\in L^\infty_{loc}({\R}_+)$ and goes to zero at $0$ and $f\in L^\infty_{loc}({\R})$. Assume that for every $\sigma>0$ the functional $E^\sigma$ satisfies the Palais Smale condition. Let ${\mathcal A}$ be an admissible family of ${\mathcal M}$ and denote
\[
\beta(\sigma):=\inf_{A\in {\mathcal A}}\ \sup_{\vec{\gamma}\in A}\ E^\sigma(\vec{\gamma})\quad.
\]
Then there exists a sequence $\sigma_k\rightarrow 0$ and $\vec{\gamma}_k\in {\mathcal M}$ such that
\[
E^{\sigma_k}(\vec{\gamma}_k)=\beta(\sigma_k)\quad,\quad DE^{\sigma_k}(\vec{\gamma}_k)=0\quad.
\]
Moreover $\vec{\gamma}_k$ satisfies the so called ``entropy condition''
\[
\p_{\sigma_k}E^{\sigma_k}(\vec{\gamma}_k)=o\lf(\frac{1}{\sigma_k\ \log\lf(\frac{1}{\sigma_k}\rg)}\rg)\quad.
\]
\hfill $\Box$
\end{Th}
\noindent A proof of this theorem is given for instance in \cite{Ri6}. Applying theorem~\ref{th-III.3-reg} to our framework gives.
\begin{Th}
\label{th-min-max-struwe}
Let $p>1$ and ${\mathcal A}$ be an admissible family  in ${\mathcal E}_{\Sigma,p}(N^n)$  such that
\be
\label{0-01}
\inf_{A\in {\mathcal A}}\ \max_{\vec{\Phi}\in A}\mbox{Area}(\vec{\Phi})=\beta^0>0\quad.
\ee
Then there exists $\sigma_k\rightarrow 0$ and a family $\vec{\Phi}_k$ of critical points of $A^{\sigma_k}_p$ satisfying
 \[
\lim_{k\rightarrow +\infty} \mbox{Area}(\vec{\Phi}_k)=\beta^0\quad\mbox{ and }\quad \sigma_k^2\int_{\Sigma}\lf[1+|\vec{\mathbb I}_{\vec{\Phi}_{k}}|^2\rg]^{p}\ dvol_{g_{\vec{\Phi}_k}}=o\lf(\frac{1}{\log\sigma_k^{-1}}\rg)\quad.
 \]
\hfill $\Box$ 
 \end{Th}
 
\subsection{The first variation of the viscous energies $A^\sigma_p$.}
Let $\vec{\Phi}$ be a smooth immersion from a closed 2-dimensional manifold $\Sigma$ into the unit sphere $S^{3}\subset {\R}^4$,
let $\vec{w}$ be an infinitesimal immersion satisfying $\vec{w}\cdot\vec{\Phi}\equiv 0$ and denote $\vec{\Phi}_t:$ a sequence
of  immersions into $S^3$ such that $d\vec{\Phi}/dt(0) =\vec{w}$.
The Gauss map of the immersion is given in local coordinates by 
\be
\label{III.4}
\vec{n}_t= \star_{{\R}^4}\lf(\vec{\Phi}_t\wedge\frac{{\p_{x_1}\vec{\Phi}_t\wedge\p_{x_2}\vec{\Phi}_t}}{|\p_{x_1}\vec{\Phi}_t\wedge\p_{x_2}\vec{\Phi}_t|}\rg)\quad.
\ee
Assuming $\vec{\Phi}$ is expressed locally in conformal coordinates and denote $e^{\la}=|\p_{x_1}\vec{\Phi}|=|\p_{x_2}\vec{\Phi}|$. We have
\be
\label{III.4a}
\vec{n}_t=\vec{n}+t\ (a_1\,\vec{e}_1+a_2\,\vec{e}_2+b\,\vec{\Phi}) +o(t)\quad,
\ee
where $\vec{e}_i=e^{-\la}\,\p_{x_i}\vec{\Phi}$. Since $\vec{n}_t\cdot\p_{x_i}\vec{\Phi}_t\equiv 0$ and $\vec{n}_t\cdot\vec{\Phi}_t\equiv 0$ we have
\be
\label{III.5}
\begin{array}{l}
\ds\frac{d\vec{n}}{dt}(0)=-\vec{n}\cdot\vec{w}\ \vec{\Phi}-\sum_{i=1}^2\vec{n}\cdot\p_{x_i}\vec{w}\ e^{-\la}\ \vec{e}_i\\[5mm]
\ds\quad\quad\quad=-\vec{n}\cdot\vec{w}\ \vec{\Phi}-\lf<\vec{n}\cdot d\vec{w}\, ,\, d\vec{\Phi}\rg>_{\gP}\quad.
\end{array}
\ee
Since $g_{ij}:=\p_{x_i}\vec{\Phi}\cdot\p_{x_j}\vec{\Phi}$, we have
\be
\label{III.6}
\frac{dg_{ij}}{dt}(0)=\p_{x_i}\vec{w}\cdot\p_{x_j}\vec{\Phi}+\p_{x_j}\vec{w}\cdot\p_{x_i}\vec{\Phi}\quad.
\ee
Since $\sum_ig_{ki}\,g^{ij}=\delta_{kj}$ and $g_{ki}=e^{2\la}\, \delta_{ki}$, we have
\be
\label{III.7}
\frac{d g^{ij}}{dt}(0)=-e^{-4\la}\ \lf[\p_{x_i}\vec{\Phi}\cdot\p_{x_j}\vec{w}+\p_{x_j}\vec{\Phi}\cdot\p_{x_i}\vec{w}\rg]\quad.
\ee
We have also using (\ref{III.5}) and (\ref{III.7})
\be
\label{III.8}
\begin{array}{l}
\ds\frac{d|d\vec{n}|_{\gP}^2}{dt}=\frac{d}{dt}\lf(\sum_{i,j=1}^2g^{ij}\p_{x_i}\vec{n}\cdot\p_{x_j}\vec{n}\rg)\\[5mm]
\ds\quad=-2\ \lf<d\vec{\Phi}\,\dot{\otimes}\,d\vec{w}\, ,\,d\vec{n}\,\dot{\otimes}\,d\vec{n}\rg>_{\gP}+2\,\lf<d\frac{d\vec{n}}{dt}\,;\, d\vec{n}\rg>_{\gP}\\[5mm]
\ds\quad=-2\ \lf<d\vec{\Phi}\,\dot{\otimes}\,d\vec{w}\, ,\,d\vec{n}\,\dot{\otimes}\,d\vec{n}\rg>_{\gP}+4\, \vec{H}\cdot\vec{w}-2\, \lf<d\lf<\vec{n}\cdot d\vec{w}\, ,\, d\vec{\Phi}\rg>_{\gP}\, ;\, d\vec{n}\rg>_{\gP}\quad,
\end{array}
\ee
where $\vec{H}$ is the mean-curvature vector given by
\[
\vec{H}:=\frac{1}{2}\sum_{i,j=1}^2g^{ij}\vec{\mathbb I}_{ij}\quad.
\]
and $\vec{\mathbb I}_{\vP}$ denotes the second fundamental form
\[
\vec{\mathbb I}_{\vP}=\sum_{i,j=1}^2\vec{\mathbb I}_{ij}\ dx_i\otimes dx_j=-\sum_{i,j=1}^2\  \p_{x_i}\vec{\Phi}\cdot\p_{x_j}\vec{n}\ \vec{n}\ dx_i\otimes dx_j\quad.
\]
Finally, we have $dvol_{g_{\vec{\Phi}}}= \sqrt{g_{11}g_{22}-g_{12}^2}\ dx_1\wedge dx_2$, hence
\be
\label{III.8-1}
\frac{d}{dt}(dvol_{g_{\vec{\Phi}}})(0)=\lf[\sum_{i=1}^2\p_{x_i}\vec{\Phi}\cdot\p_{x_i}\vec{w}\rg]\ dx_1\wedge dx_2=\lf<d\vec{\Phi}\,;d\vec{w}\rg>_{g_{\vec{\Phi}}}\ dvol_{g_{\vec{\Phi}}}\quad.
\ee
using (\ref{III.8}) and (\ref{III.8-1}) we obtain
\be
\label{III.11}
\ds\lf.\frac{d}{dt}\mbox{Area}(\vec{\Phi}_t)\rg|_{t=0}=\int_\Sigma  \lf<d\vec{\Phi}\,;d\vec{w}\rg>_{g_{\vec{\Phi}}}\ dvol_{g_{\vec{\Phi}}}\quad.
\ee
For any $p>1$ we denote
\[
F_p(\vec{\Phi}):=\int_{\Sigma}\lf[1+|\vec{\mathbb I}_{\vP}|_{\gP}^2\rg]^{p}\ dvol_{g_{\vec{\Phi}}}\quad.
\]
Using (\ref{III.5}) and (\ref{III.8}) we have
\be
\label{III.12}
\begin{array}{l}
\ds\lf.\frac{d}{dt}F_p(\vec{\Phi}_t)\rg|_{t=0}=\int_\Sigma f^p\ \lf<d\vec{\Phi}\,;d\vec{w}\rg>_{g_{\vec{\Phi}}}\ dvol_{g_{\vec{\Phi}}}
-2\, p\, \int_\Sigma f^{p-1}\ \lf<d\vec{\Phi}\,\dot{\otimes}\,d\vec{w}\, ,\,d\vec{n}\,\dot{\otimes}\,d\vec{n}\rg>_{\gP}\ dvol_{g_{\vec{\Phi}}}\\[5mm]
\ds\quad-\,2\, p\, \int_\Sigma f^{p-1}\ \lf<d\lf<\vec{n}\cdot d\vec{w}\, ,\, d\vec{\Phi}\rg>_{\gP}\, ;\, d\vec{n}\rg>_{\gP}\ dvol_{g_{\vec{\Phi}}}
+\,4\, p\, \int_\Sigma f^{p-1}\ \vec{H}\cdot\vec{w} \ dvol_{g_{\vec{\Phi}}}\quad,
\end{array}
\ee
where $f:=\lf[1+|\vec{\mathbb I}_{\vP}|^2\rg]$.
\subsection{The almost conservation laws satisfied by the critical points of $A^\sigma_p(\vP)$.}

The fact that $A^\sigma_p$ is $C^1$ in ${\mathcal E}_{\Sigma,p}$ is quite standard for $p>1$. Indeed, in local coordinates the functional has the form
\[
\int_{\Sigma}e(\vec{\Phi},\nabla\vec{\Phi},\nabla^2\vec{\Phi})\ dx^2\quad,
\]
where $e$ is a $C^\infty$ function.
Let $\vP$ be a critical point in ${\mathcal E}_{\Sigma,p}$ of $A^\sigma_p$. We then have
\be
\label{III.14}
\begin{array}{l}
\ds\vec{\Phi}\wedge d^{\ast_{\gP}}\lf[\lf[1+\sigma^2\,f^p\rg] \ d\vec{\Phi}\rg]-\,2\, p\,\sigma^2\ \vec{\Phi}\wedge\,  d^{\ast_{\gP}}\lf[\ d^{\ast_{\gP}}\lf[ f^{p-1}\ d\vec{n}\rg]\cdot d\vec{\Phi}\ \vec{n}\rg]\\[5mm]
\ds-\,2\, p\,\sigma^2\ \vec{\Phi}\wedge\,  d^{\ast_{\gP}}\lf[\ f^{p-1}\ (d\vec{n}\,\dot{\otimes}\,d\vec{n})\res_{\gP} d\vec{\Phi}\rg] +\ 4\, p\ \sigma^2\ f^{p-1}\ \vec{\Phi}\wedge\vec{H}=0\quad
\mbox{ in }{\mathcal D}'(\Sigma)\quad,
\end{array}
\ee
where $f:=\lf[1+|\vec{\mathbb I}_{\vP}|^2\rg]$ as above, $(d\vec{n}\,\dot{\otimes}\,d\vec{n})\res_{\gP} d\vec{\Phi}$ is the contraction given in local conformal coordinates by 
\[
(d\vec{n}\,\dot{\otimes}\,d\vec{n})\res_{\gP} d\vec{\Phi}:=e^{-2\la}\ \sum_{i,j=1}^2\ \p_{x_i}\vec{n}\cdot\p_{x_j}\vec{n}\ \p_{x_j}\vP\ dx_i\quad,
\]
 and $d^{\ast_{\gP}}$ is the adjoint of $d$ for the $L^2$ norm on $\Sigma$ with respect to the metric $\gP$ induced by the immersion $\vec{\Phi}$. It coincides with $-\,e^{-2\la} \mbox{ div}\cdot$ in conformal coordinates. In conformal coordinates again  the equation becomes then
\be
\label{III.14-a}
\begin{array}{l}
\ds\vec{\Phi}\wedge\mbox{div}\lf[\lf[1+\sigma^2 f^p\rg]\nabla\vec{\Phi}-2\, p\,\sigma^2\  e^{-2\la}\, f^{p-1}\ \lf<(\nabla\vec{n}\dot{\otimes}\nabla\vec{n};\nabla\vec{\Phi}\rg>\rg.\\[5mm]
\ds\lf.+2\,p\,\sigma^2\ e^{-2\la}\ \mbox{div}\lf[f^{p-1}\,\nabla\vec{n}\rg]\cdot\nabla\vec{\Phi}\ \vec{n}\rg]- 4\, p\, \sigma^2\ f^{p-1}\ \vec{\Phi}\wedge\vec{H}=0
\end{array}
\ee
We rewrite the first term in the second line.
\be
\label{A-2}
\begin{array}{l}
\ds 2\,p\,\sigma^2\ e^{-2\la}\, \mbox{div}\lf[f^{p-1}\,\nabla\vec{n}\rg]\cdot\nabla\vec{\Phi}\ \vec{n} \\[5mm]
\ds \quad\quad\quad=2\,p\,\sigma^2\ e^{-2\la}\, \mbox{div}\lf[f^{p-1}\,\lf[\nabla\vec{n}+\,H\,\nabla\vec{\Phi}\rg]\rg]\cdot\nabla\vec{\Phi}\ \vec{n} -2\,p\,\sigma^2\ \nabla\lf[f^{p-1}\ H\ \rg]\ \vec{n}\quad.
\end{array}
\ee
The {\it trace free part } of the second fundamental form is denoted
\[
\vec{\mathbb I}^{\,0}:=\vec{\mathbb I}- \vec{H}\ g\quad.
\]
In  coordinates and in codimension 1 one has
\[
\vec{\mathbb I}^{\,0}={\mathbb I}^{\,0}\, \vec{n}=-\sum_{i,j=1}^2\lf[\p_{x_i}\vec{n}\cdot\p_{x_j}\vec{\Phi}+H\,\p_{x_i}\vec{\Phi}\cdot\p_{x_j}\vec{\Phi}\rg]\ dx_{i}\otimes dx_j\quad.
\]
For any $k=1,2$ after some computations we obtain
\[
\begin{array}{l}
\ds\sum_{i=1}^2\p_{x_i}\lf[f^{p-1}\ \lf[\p_{x_i}\vec{n}+H\,\p_{x_i}\vec{\Phi}\rg]\rg]\cdot\p_{x_k}\vec{\Phi}\ \vec{n}=-\p_{x_k}\lf[f^{p-1}\ {\mathbb I}^0_{k,k}   \rg]\ \vec{n}-\p_{x_{k+1}}\lf[f^{p-1}\ {\mathbb I}^0_{k+1,k}   \rg]\ \vec{n}\quad.
\end{array}
\]
Denoting $\ov{\nabla}\cdot:=(\p_{x_1}\cdot,-\p_{x_2}\cdot)$ and $(\ov{\nabla})^\perp\cdot:=(\p_{x_2}\cdot,\p_{x_1}\cdot)$, we have then
\be
\label{A-3}
2\,p\,\sigma^2\ e^{-2\la}\, \mbox{div}\lf[f^{p-1}\,\lf[\nabla\vec{n}+\,H\,\nabla\vec{\Phi}\rg]\rg]\cdot\nabla\vec{\Phi}\  =-2\, p\,\sigma^2\ e^{-2\la}\, 
\lf[\ov{\nabla}\lf[ f^{p-1}\, {\mathbb I}^0_{11}\rg] + (\ov{\nabla})^\perp\lf[ f^{p-1}\, {\mathbb I}^0_{12}\rg]\ \rg]\quad.
\ee
Combining (\ref{A-2}) and (\ref{A-3}) gives
\be
\label{A-3-a}
\begin{array}{l}
\ds 2\,p\,\sigma^2\ e^{-2\la}\, \mbox{div}\lf[f^{p-1}\,\nabla\vec{n}\rg]\cdot\nabla\vec{\Phi}\ \vec{n} =-2\,p\,\sigma^2\ \nabla\lf[f^{p-1}\ \vec{H}\ \rg]
+\,2\,p\,\sigma^2\ f^{p-1}\ {H}\ \ \nabla\vec{n}\\[5mm]
-\,2\, p\,\sigma^2\ e^{-2\la}\, 
\lf[\ov{\nabla}\lf[ f^{p-1}\, {\mathbb I}^0_{11}\rg] + (\ov{\nabla})^\perp\lf[ f^{p-1}\, {\mathbb I}^0_{12}\rg]\rg]\ \vec{n}\quad.
\end{array}
\ee
So the equation (\ref{III.14-a}) becomes
\be
\label{III.14-b}
\begin{array}{l}
\ds\vec{\Phi}\wedge\mbox{div}\lf[\lf[1+\sigma^2 f^p\rg]\nabla\vec{\Phi}-2\,p\,\sigma^2\ \nabla\lf[f^{p-1}\ \vec{H}\ \rg]-2\, p\,\sigma^2\  e^{-2\la}\, f^{p-1}\ \lf<\nabla\vec{n}\dot{\otimes}\nabla\vec{n};\nabla\vec{\Phi}\rg>\rg.\\[5mm]
\ds\lf.+\,2\,p\,\sigma^2\ f^{p-1}\ {H}\ \ \nabla\vec{n}\ -\,2\, p\,\sigma^2\ e^{-2\la}\, 
\lf[\ov{\nabla}\lf[ f^{p-1}\, {\mathbb I}^0_{11}\rg] + (\ov{\nabla})^\perp\lf[ f^{p-1}\, {\mathbb I}^0_{12}\rg]\rg]\ \vec{n}\rg]=4\, p\, \sigma^2\ f^{p-1}\ \vec{\Phi}\wedge\vec{H}\quad.
\end{array}
\ee
The equation (\ref{III.14-b}) can be rewritten in an exact divergence free equation of the form $\mbox{div}(\vec{\Phi}\wedge\cdots)=0$, that is in an exact conservation law which is coming from the $SO(4)$ invariance of the problem in the target.
However, since we are interested in general targets, we don't want to take advantage of the ``roundness'' of $S^3$ and we shall rewrite (\ref{III.14-b}) in an ``almost conservation law''
which is more generic and which holds in ${\mathcal D}'(\Sigma)$ . It  is due this time to the translation invariance of the integrand of $F_p$ in ${\R}^4$ in relation with the Noether theorem as observed in \cite{Be}. However the fact that we don't get exactly get
a conservation law is coming from the fact that the constraint to take values into the closed sub-manifold $S^3$ is not translation invariant. This pointwize constraint is ``generating'' additional
terms (i.e. the last term in the l.h.s. and the full r.h.s. of (\ref{III.15}) ) in comparison to the identity we would get if we would release this constraint. Nevertheless these additional terms  happen to be of much lower degree and are not going to perturb the arguments in the section below as if we would be dealing with an exact conservation law. This is why we are speaking about an ``almost conservation law''.
\be
\label{III.15}
\begin{array}{l}
-\mbox{div}\lf[\lf[1+\sigma^2 f^p\rg]\nabla\vec{\Phi}-2\,p\,\sigma^2\ \nabla\lf[f^{p-1}\ \vec{H}\ \rg]-2\, p\,\sigma^2\  e^{-2\la}\, f^{p-1}\ \lf<\nabla\vec{n}\dot{\otimes}\nabla\vec{n};\nabla\vec{\Phi}\rg>\rg.\\[5mm]
\ds\lf.+\,2\,p\,\sigma^2\ f^{p-1}\ {H}\ \ \nabla\vec{n}\ -\,2\, p\,\sigma^2\ e^{-2\la}\, 
\lf[\ov{\nabla}\lf[ f^{p-1}\, {\mathbb I}^0_{11}\rg] + (\ov{\nabla})^\perp\lf[ f^{p-1}\, {\mathbb I}^0_{12}\rg]\rg]\ \vec{n}\rg]+4\, p\, \sigma^2\ f^{p-1}\ \vec{H}
\\[5mm]
\quad=\lf[1+\sigma^2\ (1-p)\ f^p+\, p\,\sigma^2 f^{p-1}\rg]\ |\nabla\vec{\Phi}|^2\,\vec{\Phi}\quad.
\end{array}
\ee 
Finally we end up this section by quoting the following  theorem
\begin{Th}
\label{th-reg}
Let $p\ge1$ and $\vec{\Phi}$ be an element in the space ${\mathcal E}_{\Sigma,p}$ of $W^{2,2p}-$immersions of a closed surface $\Sigma$.
Assume $\vec{\Phi}$ is a critical point of $A^\sigma_p(\vP)$ then $\vec{\Phi}$ is $C^\infty$ in any conformal parametrization.\hfill $\Box$
\end{Th}
\begin{Rm}
A proof of theorem~\ref{th-reg} has been given in \cite{KLL}  and for $C^1$ into the euclidian space.
The method of proof in \cite{KLL} relies on the work of J.Langer with the decomposition of the immersion into
the union of graphs.  See also \cite{RiB} for a  proof making use of the underlying conservation laws.
\hfill $\Box$
\end{Rm}
\subsection{Proof of theorem~\ref{th-min-max}.} Combine theorem~\ref{th-min-max-struwe} and theorem~\ref{th-main}, this gives theorem~\ref{th-min-max}.\hfill $\Box$

\section{The passage to the limit $\sigma\rightarrow 0$ with controlled conformal class.}
\reset
The goal of the present section is to prove the following theorem
\begin{Th}
\label{th-limit}
Let $p>1$ and let $\vec{\Phi}_{k}$ be a sequence of critical points of $A_p^{\sigma_k}$ in the class ${\mathcal E}_{\Sigma,p}$ where $\sigma_k\rightarrow 0$ and satisfying
\be
\label{VII.3}
0<\limsup_{k\rightarrow +\infty} \mbox{Area}(\vec{\Phi}_{k})<+\infty\quad,
\ee
and
\be
\label{VII.2}
\sigma_k^2\ F_p(\vec{\Phi}_{k})=\sigma_k^2\ \int_{\Sigma}\lf[1+|{\mathbb I}_{\vec{\Phi}_{k}}|^{2}_{g_{\vec{\Phi}_{k}}}\rg]^p dvol_{g_{\vec{\Phi}_{k}}}=o\lf(\frac{1}{\log(1/\sigma_k)}\rg)\quad.
\ee
Assume moreover that the conformal class associated to $(\Sigma,g_{\vec{\Phi}_k})$ is precompact in the moduli space, then, modulo extraction
of a subsequence, there exists a closed riemann surface $(S,h_0)$ with genus$(S)\le$genus$(\Sigma)$, a weakly conformal 
map $\vec{\Phi}_\infty$ from $S$ into $N^n$ and an integer valued map $N\in L^\infty(S,{\N})$ such that
\[
\lim_{k\rightarrow +\infty} A^{\sigma_k}(\vec{\Phi}_k)=\frac{1}{2}\int_SN\,|d\vec{\Phi}_\infty|^2_{h_0}\ dvol_{h_0}\quad.
\]
Moreover the push forward of $S$ by $\vec{\Phi}_\infty$ together with the multiplicity $N$ defines  an oriented stationary integer varifold and the oriented varifold $|T_k|$ equal to the push-forward by $\vec{\Phi}_k$  of $\Sigma$ converges in the sense of Radon measures towards the oriented stationary integer varifold associated to $\vec{\Phi}_\infty$. The surface $S$ is moreover either equal to the union of $\Sigma$ with finitely many copies of $S^2$ or is equal to finitely many copies of $S^2$. \hfill $\Box$
\end{Th}
\begin{Rm}
\label{remark-genus}
Observe that in theorem~\ref{th-limit}, due to the assumption about the controlled conformal class, there can be a genus jump genus$(S)<$genus$(\Sigma)$ only if the area
vanishes on the main part of the Riemann surface and $\Phi_\infty(S)$ is going to be a bouquet of minimal sphere. This cannot be excluded a priori \hfill $\Box$
\end{Rm}

{\bf In this section we shall then assume that $\vec{\Phi}_k$ is conformal from a sequence of riemannian surfaces $(\Sigma,g_k)$ into $S^3$ for which the underlying Riemann structure  is pre-compact in the moduli space of $\Sigma$.}

In order to prove theorem~\ref{th-limit} we shall need several lemma. 

\medskip

\begin{Lm}
\label{lm-VII.1} {\bf [Monotonicity Formula]}
Under the assumptions of theorem~\ref{th-limit} the sequence of varifolds $|T_k|$  equal to the push forward of $\Sigma$ by 
$\vec{\Phi}_{k}$ converges, modulo extraction of a subsequence,  towards  a  stationary  varifold. In particular, introducing the Radon measure in $S^3$ given by 
\be
\label{VII.3-a}
<\mu_k,\varphi>:=\int_\Sigma\varphi(\vec{\Phi}_k)\ dvol_{g_{\vec{\Phi}_{k}}}\quad,
\ee
$\mu_k$ converges modulo extraction of a subsequence to a limiting Radon measure $\mu_\infty$ satisfying the following monotonicity formula
\be
\label{VII.4}
\forall \ \vq\in\mbox{supp}( \mu_\infty)\quad\forall\, r>0\quad\quad\frac{d}{dr}\lf[\frac{e^{C\, r}\mu_\infty(B_r(\vq))}{r^2}\rg]\ge 0\quad.
\ee
for some $C>0$ independent of $\vq$ and $r$. \hfill $\Box$
\end{Lm}
{\bf Proof of lemma~\ref{lm-VII.1}.}
The monotonicity formula for the limiting measure $\mu_\infty$ is a direct consequence of the fact that $|T_k|$ converges towards  a  stationary  varifold (see \cite{Al} and \cite{Si}). So it would suffices to prove this last fact in order to get
(\ref{VII.4}). However the proof of both statements (that can be proven independently of each other) are very similar. In the first case it suffices to prove that for any vector field $\vec{X}$ 
we have
\be
\label{M-0}
\lim_{k\rightarrow +\infty}\int_{M_k} \mbox{div}_{M_k}\vec{X}\ d{\mathcal H}^2=
\lim_{k\rightarrow +\infty}\int_{\Sigma}\lf[\sum_{i=1}^4\lf<\p_{y_i}\vec{X}(\vec{\Phi}_k)\ \nabla\Phi^i_k,\nabla\vec{\Phi}_k\rg>-\vec{X}(\vec{\Phi}_k)\cdot\vec{\Phi}_k\ |\nabla\vec{\Phi}_k|^2\rg]\ dx^2=0\quad,
\ee
where $M_k:=\vec{\Phi}_{k}(\Sigma)$ and $\vec{\Phi}_k=(\Phi^1_k,\cdots,\Phi^4_k)$. The computations for proving (\ref{M-0}) are more or less the same as the one for proving (\ref{VII.4}) and we shall only
present the later since we shall revisit them in the forthcoming lemma~\ref{lm-quanti}.

The explicit mention of the indices $\sigma_k$ and $k$ can be deleted when there is no possible confusion. For any $\vec{q}\in S^3$ and any radius $r$ small enough, {\it Simon's monotonicity formula} ( see \cite{Si}  chapter 4) applied
to  $\vec{\Phi}(\Sigma)$ (which is smooth immersion for any $k$) which is seen as a {\it varifold} from ${\R}^4$ gives
\be
\label{M-1}
\begin{array}{l}
\ds\frac{d}{dr}\lf[\frac{1}{r^2}\int_{\vec{\Phi}^{-1}(B^4_r(\vec{q}))}\ dvol_{g_{\vec{\Phi}}}\rg]=\frac{d}{dr}\lf[\int_{\vec{\Phi}^{-1}(B^4_r(\vec{q}))}\ \frac{|(\vec{n}\wedge\vec{\Phi})\res(\vec{\Phi}-\vec{q})|^2}{|\vec{\Phi}-\vec{q}|^4} \  dvol_{g_{\vec{\Phi}}} \rg]\\[5mm]
\ds\quad -\frac{1}{2\,r^3}\ \int_{\vec{\Phi}^{-1}(B^4_r(\vec{q}))}(\vec{\Phi}-\vec{q})\cdot d^{\ast_g}d\vec{\Phi}\ dvol_{g_{\vec{\Phi}}} \\[5mm]
\ds\quad\quad\ge  -\frac{1}{2\,r^3}\ \int_{\vec{\Phi}^{-1}(B^4_r(\vec{q}))}(\vec{\Phi}-\vec{q})\cdot d^{\ast_g}d\vec{\Phi}\ dvol_{g_{\vec{\Phi}}} \quad,
\end{array}
\ee
where we have used that the first term in the r.h.s. of (\ref{M-1}) is non negative \footnote{Indeed we are taking the derivative of an integral of a positive integrand over a bigger and bigger set.}. Thanks to equation (\ref{III.15})  we obtain
\be
\label{M-2}
\begin{array}{l}
\ds-\int_{\vec{\Phi}^{-1}(B^4_r(\vec{q}))}(\vec{\Phi}-\vec{q})\cdot d^{\ast_g}d\vec{\Phi}\ dvol_{g_{\vec{\Phi}}}=\int_{\vec{\Phi}^{-1}(B^4_r(\vec{q}))}(\vec{\Phi}-\vec{q})\cdot \Delta\vec{\Phi}\ dx^2\\[5mm]
\ds=-\int_{\vec{\Phi}^{-1}(B^4_r(\vec{q}))}(\vec{\Phi}-\vec{q})\cdot \mbox{div}\lf[    \sigma^2\ f^{p-1}\ [f\,\nabla\vec{\Phi}-2\,p\ \lf({H}\ \nabla\vec{n}-e^{-2\la}<\nabla\vec{n}\dot{\otimes}\nabla\vec{n};\nabla\vec{\Phi}>\rg)] \rg]\, dx^2\\[5mm]
\ds +\, 2\  p\ \sigma^2\  \int_{\vec{\Phi}^{-1}(B^4_r(\vec{q}))}(\vec{\Phi}-\vec{q})\cdot \mbox{div}\lf[ e^{-2\la}\ \lf[\ov{\nabla}\lf[ f^{p-1}\, {\mathbb I}^0_{11}\rg] + (\ov{\nabla})^\perp\lf[ f^{p-1}\, {\mathbb I}^0_{12}\rg]\ \rg]\vec{n}\rg]\ dx^2\\[5mm]
\ds +\, 2\ p\ \sigma^2\  \int_{\vec{\Phi}^{-1}(B^4_r(\vec{q}))}(\vec{\Phi}-\vec{q})\cdot \Delta\lf[f^{p-1}\ \vec{H}\rg] \ dx^2\\[5mm]
\ds -\int_{\vec{\Phi}^{-1}(B^4_r(\vec{q}))}\lf[1+ \, (1-\,p)\ \sigma^2\ f^{p}+ p\, \sigma^2\ f^{p-1}\rg]\ (\vec{\Phi}-\vec{q})\cdot\vec{\Phi}\,|\nabla\vec{\Phi}|^2\ dx^2\quad.
\end{array}
\ee
Regarding the  last line, observe in one hand that $(\vec{\Phi}-\vec{q})\cdot\vec{\Phi}=1-\cos{(\vec{\Phi},\vec{q})}\,=O(r^2)$ hence
\be
\label{M-3}
\lf|\frac{1}{r^3}\int_{\vec{\Phi}^{-1}(B^4_r(\vec{q}))}(\vec{\Phi}-\vec{q})\cdot\vec{\Phi}\,|\nabla\vec{\Phi}|^2\ dx^2\rg|\le \frac{C}{r}\int_{\vec{\Phi}^{-1}(B^4_r(\vec{q}))}\ dvol_{g_{\vec{\Phi}}}
\ee
and in the other hand, again for fixed $r$ and $\vec{q}$, as $k\rightarrow +\infty$
\be
\label{M-4}
\begin{array}{l}
\ds\lf|\int_{\vec{\Phi}^{-1}(B^4_r(\vec{q}))}\lf[ \, (1-p)\ \sigma^2\ f^{p}+\, p\, \sigma^2\ f^{p-1}\rg]\ (\vec{\Phi}-\vec{q})\cdot\vec{\Phi}\,|\nabla\vec{\Phi}|^2\ dx^2\rg|\\[5mm]
\ds\quad\quad\le C\ \sigma^2\ F_p(\vec{\Phi})+C\ \sigma^2 \ M(T)^{1/p}\ [F_p(\vec{\Phi})]^{1-1/p}\quad\rightarrow 0\quad.
\end{array}
\ee
Integrating by parts each of the two first lines in the r.h.s. of (\ref{M-2}) gives 
\be
\label{M-6}
\begin{array}{l}
\ds -\int_{\vec{\Phi}^{-1}(B^4_r(\vec{q}))}(\vec{\Phi}-\vec{q})\cdot \mbox{div}\lf[    \sigma^2\ f^{p-1}\ [f\,\nabla\vec{\Phi}-2\,p\ \lf({H}\ \nabla\vec{n}-e^{-2\la}<\nabla\vec{n}\dot{\otimes}\nabla\vec{n};\nabla\vec{\Phi}>\rg)] \rg]\, dx^2\\[5mm]
\ds +\, 2\  p\ \sigma^2\  \int_{\vec{\Phi}^{-1}(B^4_r(\vec{q}))}(\vec{\Phi}-\vec{q})\cdot \mbox{div}\lf[ e^{-2\la}\ \lf[\ov{\nabla}\lf[ f^{p-1}\, {\mathbb I}^0_{11}\rg] + (\ov{\nabla})^\perp\lf[ f^{p-1}\, {\mathbb I}^0_{12}\rg]\ \rg]\vec{n}\rg]\ dx^2\\[5mm]
\ds=\sigma^2\ \int_{\vec{\Phi}^{-1}(B^4_r(\vec{q}))}\ f^{p-1}\ \lf[f\,|\nabla\vec{\Phi}|^2-2\,p\ {H}\ \nabla\vec{n}\cdot\nabla\vec{\Phi}+2\, p\,(f-1)\, e^{2\la}] \rg]\ dx^2\\[5mm]
\ds-\ \sigma^2\ \int_{\vec{\Phi}^{-1}(\p B^4_r(\vec{q}))}  f^p\ \p_\nu\vec{\Phi}\cdot (\vec{\Phi}-\vec{q})-\, 2\, p\ f^{p-1}\ H\ \p_\nu\vec{n}\cdot (\vec{\Phi}-\vec{q})dl\\[5mm]
\ds+\ 2\ p\ \sigma^2\ \int_{\vec{\Phi}^{-1}(\p B^4_r(\vec{q}))}\ f^{p-1}\  \lf<\p_\nu\vec{n}\cdot\nabla\vec{n},\nabla\vec{\Phi}\cdot(\vec{\Phi}-\vec{q})\rg>\ dl \\[5mm]
\ds+\ 2\ p\ \sigma^2\  \int_{\vec{\Phi}^{-1}(\p B^4_r(\vec{q}))}\ e^{-2\la}\ (\vec{\Phi}-\vec{q})\cdot\vec{n} \lf[ \nu_1\ \p_{x_1}\lf[ f^{p-1}\, {\mathbb I}^0_{11}\rg] - \nu_2\ \p_{x_2}\lf[ f^{p-1}\, {\mathbb I}^0_{11}\rg] \rg]\ dl\\[5mm]
\ds+\ 2\ p\ \sigma^2\  \int_{\vec{\Phi}^{-1}(\p B^4_r(\vec{q}))}\ e^{-2\la}\ (\vec{\Phi}-\vec{q})\cdot\vec{n} \lf[  \nu_1\ \p_{x_2}\lf[ f^{p-1}\, {\mathbb I}^0_{12}\rg] + \nu_2\ \p_{x_1}\lf[ f^{p-1}\, {\mathbb I}^0_{12}\rg] \rg]\ dl\quad,
\end{array}
\ee
where $\nu$ is the outward unit (in the coordinates ) normal to ${\vec{\Phi}^{-1}(B^4_r(\vec{q}))}$ and is given explicitly by
\[
\nu=(\p_{x_1}|\vec{\Phi}-\vec{q}|\,,\p_{x_2}|\vec{\Phi}-\vec{q}|)/|\nabla |\vec{\Phi}-\vec{q}||\quad.
\]
This is nothing but the normalized  gradient of the function distance to $\vec{q}$. We clearly have
\be
\label{M-7}
\lim_{k\rightarrow+\infty}\sigma^2\ \int_{\vec{\Phi}^{-1}(B^4_r(\vec{q}))}\ f^{p-1}\ \lf[f\,|\nabla\vec{\Phi}|^2-2\,p\ {H}\ \nabla\vec{n}\cdot\nabla\vec{\Phi}+\,2\,p\,e^{2\la} \, (f-1)] \rg]\ dx^2=0\quad.
\ee
Multiplying (\ref{M-6}) by $\chi(r)/r^3$ where $\chi$ is an arbitrary compactly supported function  in ${\R}_+^\ast$ and integrating over ${\R}_+^\ast$ gives successively
\be
\label{M-8}
\begin{array}{l}
\ds\sigma^2\ \int_{{\R}_+}\ \chi(r)\ \frac{dr}{r^3}\int_{\vec{\Phi}^{-1}(\p B^4_r(\vec{q}))}  \lf[f^p\ \p_\nu\vec{\Phi}\cdot (\vec{\Phi}-\vec{q})-2\, p\ f^{p-1}\ H\ \p_\nu\vec{n}\cdot (\vec{\Phi}-\vec{q})\rg]\ dl\\[5mm]
\ds+\,2\,p\,\sigma^2\,\int_{{\R}_+}\ \chi(r)\ \frac{dr}{r^3}\int_{\vec{\Phi}^{-1}(\p B^4_r(\vec{q}))} f^{p-1}\  \lf<\p_\nu\vec{n}\cdot\nabla\vec{n},\nabla\vec{\Phi}\cdot(\vec{\Phi}-\vec{q})\rg>\ dl\\[5mm]
\ds= \sigma^2\ \int_\Sigma\ \chi(|\vec{\Phi}-\vec{q}|)\ \lf[f^p\ \frac{|\nabla|\vec{\Phi}-\vec{q}||^2}{|\vec{\Phi}-\vec{q}|^2}-2\, p\ f^{p-1}\ H \frac{\nabla|\vec{\Phi}-\vec{q}|}{|\vec{\Phi}-\vec{q}|^3}\cdot \lf<\nabla\vec{n}\cdot (\vec{\Phi}-\vec{q})\rg> \rg]\, dx^2\\[5mm]
\ds +\ 2\, p\, \sigma^2\ \int_\Sigma\ \chi(|\vec{\Phi}-\vec{q}|) \ f^{p-1}\  \lf<\frac{\nabla|\vec{\Phi}-\vec{q}|}{|\vec{\Phi}-\vec{q}|^3}\cdot\nabla\vec{n}, \nabla\vec{n}\,\nabla\vec{\Phi}\cdot(\vec{\Phi}-\vec{q}) \rg> \ dx^2 \\[5mm]
\quad\quad\longrightarrow 0\quad\mbox{ as }k\rightarrow +\infty\quad,
\end{array}
\ee
where we have bound the r.h.s. of (\ref{M-8}) by a constant depending on $\chi$ times $\sigma^2 F_p(\vec{\Phi})$ . We also obtain 
\be
\label{M-9}
\begin{array}{l}
\ds-\ 2\ p\ \sigma^2\   \int_{{\R}_+}\ \chi(r)\ \frac{dr}{r^3}\int_{\vec{\Phi}^{-1}(\p B^4_r(\vec{q}))}\ e^{-2\la}\ (\vec{\Phi}-\vec{q})\cdot \vec{n}\lf[ \nu_1\ \p_{x_1}\lf[ f^{p-1}\, {\mathbb I}^0_{11}\rg]\rg]\ dl\\[5mm]
\ds+\ 2\ p\ \sigma^2\   \int_{{\R}_+}\ \chi(r)\ \frac{dr}{r^3}\int_{\vec{\Phi}^{-1}(\p B^4_r(\vec{q}))}\ e^{-2\la}\ (\vec{\Phi}-\vec{q})\cdot \vec{n}\ \lf[ \nu_2\ \p_{x_2}\lf[ f^{p-1}\, {\mathbb I}^0_{11}\rg] \rg]\ dl\\[5mm]
\ds=-\ p\ \sigma^2\ \int_{\Sigma}\chi(|\vec{\Phi}-\vec{q}|)\ \frac{(\vec{\Phi}-\vec{q})}{|\vec{\Phi}-\vec{q}|^4}\cdot \vec{n}\ \lf[ e^{-2\la}\ \p_{x_1}|\vec{\Phi}-\vec{q}|^2\ \p_{x_1}\lf[ f^{p-1}\, {\mathbb I}^0_{11}\rg] \rg]\ dx^2\\[5mm]
\ds\quad+\ p\ \sigma^2\ \int_{\Sigma }\chi(|\vec{\Phi}-\vec{q}|)\ \frac{(\vec{\Phi}-\vec{q})}{|\vec{\Phi}-\vec{q}|^4}\cdot\vec{n}\  \lf[e^{-2\la}\ \p_{x_2}|\vec{\Phi}-\vec{q}|^2\ \p_{x_2}\lf[ f^{p-1}\, {\mathbb I}^0_{11}\rg] \rg]\ dx^2\quad.
\end{array}
\ee
Integrating by parts  the r.h.s of (\ref{M-9}), we have 
\be
\label{M-10}
\begin{array}{l}
\ds-\ p\ \sigma^2\ \int_{\Sigma}\chi(|\vec{\Phi}-\vec{q}|)\ \frac{(\vec{\Phi}-\vec{q})}{|\vec{\Phi}-\vec{q}|^4}\cdot\vec{n}\ \lf[ e^{-2\la}\ \p_{x_1}|\vec{\Phi}-\vec{q}|^2\ \p_{x_1}\lf[ f^{p-1}\, {\mathbb I}^0_{11}\rg]\rg]\, dx^2\\[5mm]
\ds\quad+\ p\ \sigma^2\ \int_{\Sigma}\chi(|\vec{\Phi}-\vec{q}|)\ \frac{(\vec{\Phi}-\vec{q})}{|\vec{\Phi}-\vec{q}|^4}\cdot \vec{n}\ \lf[e^{-2\la}\ \p_{x_2}|\vec{\Phi}-\vec{q}|^2\ \p_{x_2}\lf[ f^{p-1}\, {\mathbb I}^0_{11}\rg] \rg]\, dx^2\\[5mm]
\ds\quad=\ p\ \sigma^2\ \int_{\Sigma} f^{p-1}\, {\mathbb I}^0_{11}\ \ov{\nabla}\lf[\chi(|\vec{\Phi}-\vec{q}|)\ \frac{(\vec{\Phi}-\vec{q})\cdot\vec{n}}{|\vec{\Phi}-\vec{q}|^4}\rg] e^{-2\la}\ \nabla|\vec{\Phi}-\vec{q}|^2\ dx^2\\[5mm]
\ds\quad+\ p\ \sigma^2\ \int_{\Sigma} f^{p-1}\, {\mathbb I}^0_{11}\cdot\chi(|\vec{\Phi}-\vec{q}|)\ \frac{(\vec{\Phi}-\vec{q})\cdot\vec{n}}{|\vec{\Phi}-\vec{q}|^4}\ \ov{\nabla}\lf[e^{-2\la}\ \nabla|\vec{\Phi}-\vec{q}|^2\rg]\ dx^2\quad.
\end{array}
\ee
We recall that we have respectively
\be
\label{M-5}
\ov{\nabla}\cdot\lf(e^{-2\la}\ \nabla\vec{\Phi}\rg)=2\,e^{-2\la}\ \vec{\mathbb I}_{11}^0\quad\mbox{ and }\quad (\ov{\nabla})^\perp\cdot\lf(e^{-2\la}\ \nabla\vec{\Phi}\rg)=\,2\,e^{-2\la}\ \vec{\mathbb I}_{12}^0\quad.
\ee
Combining these identities with the fact that $\vec{\Phi}$ is conformal we deduce that
\be
\label{M-10-a}
\begin{array}{l}
\ds\ov{\nabla}\lf[e^{-2\la}\ \nabla|\vec{\Phi}-\vec{q}|^2\rg]=2\,\ov{\nabla}\lf[e^{-2\la}\ \nabla(\vec{\Phi}-\vec{q})\rg]\cdot(\vec{\Phi}-\vec{q})+2\, e^{-2\la}\ \nabla(\vec{\Phi}-\vec{q})\cdot \ov{\nabla}(\vec{\Phi}-\vec{q})\\[5mm]
\ds\quad\quad=4\ e^{-2\la}\ \vec{\mathbb I}_{11}^0\cdot(\vec{\Phi}-\vec{q})\quad.
\end{array}
\ee
Combining (\ref{M-10}) and (\ref{M-10-a}) and observing that we have the following  pointwise upper bound
\[
\lf|\ov{\nabla}\lf[\chi(|\vec{\Phi}-\vec{q}|)\ \frac{(\vec{\Phi}-\vec{q})\cdot\vec{n}}{|\vec{\Phi}-\vec{q}|^4}\rg] \rg|\le C\ \lf[\|\chi'\|_{\infty} \ d^{-3}_{\chi}+\|\chi\|_{\infty}d^{-4}_{\chi}\rg]\ |\nabla\vec{\Phi}|(x)+\|\chi\|_{\infty}\ d^{-3}_{\chi}\ |\nabla\vec{n}|(x)\quad,
\]
where $d_\chi$ is the distance of the support of $\chi$ to zero we deduce
\be
\label{M-10-b}
\begin{array}{l}
\ds\lf|-\ p\ \sigma^2\ \int_{\Sigma}\chi(|\vec{\Phi}-\vec{q}|)\ \frac{(\vec{\Phi}-\vec{q})}{|\vec{\Phi}-\vec{q}|^4}\cdot \lf[ e^{-2\la}\ \p_{x_1}|\vec{\Phi}-\vec{q}|^2\ \p_{x_1}\lf[ f^{p-1}\, \vec{\mathbb I}^0_{11}\rg]\rg]\, dx^2\rg.\\[5mm]
\ds\quad\lf.+\ p\ \sigma^2\ \int_{\Sigma}\chi(|\vec{\Phi}-\vec{q}|)\ \frac{(\vec{\Phi}-\vec{q})}{|\vec{\Phi}-\vec{q}|^4}\cdot  \lf[e^{-2\la}\ \p_{x_2}|\vec{\Phi}-\vec{q}|^2\ \p_{x_2}\lf[ f^{p-1}\, \vec{\mathbb I}^0_{11}\rg] \rg]\, dx^2\rg|\\[5mm]
\ds\quad\quad\le C_{\chi}\ \sigma^2\ F_p(\vec{\Phi})+C_{\chi}\ \sigma^2 \ M(T)^{1/p}\ [F_p(\vec{\Phi})]^{1-1/p}\quad\rightarrow 0\quad.
\end{array}
\ee
The control of the last term of the r.h.s. of (\ref{M-6})  is performed  similarly to the preceding one  following each step between (\ref{M-9}) and (\ref{M-10-b}). So finally deduce that for any $\chi$ compactly supported in ${\R}_+^\ast$ we have 
\be
\label{M-11}
\begin{array}{l}
\ds - \int_{{\R}_+}\ \chi(r)\ \frac{dr}{r^3}\int_{\vec{\Phi}^{-1}(B^4_r(\vec{q}))}(\vec{\Phi}-\vec{q})\cdot \mbox{div}\lf[    \sigma^2\ f^{p-1}\ [f\,\nabla\vec{\Phi}-\,2\,p\ \lf({H}\ \nabla\vec{n}-e^{-2\la}<\nabla\vec{n}\dot{\otimes}\nabla\vec{n};\nabla\vec{\Phi}>\rg)] \rg]\ dx^2\\[5mm]
\ds +\, 2\  p\ \sigma^2\   \int_{{\R}_+}\ \chi(r)\ \frac{dr}{r^3}\int_{\vec{\Phi}^{-1}(B^4_r(\vec{q}))}(\vec{\Phi}-\vec{q})\cdot \mbox{div}\lf[ e^{-2\la}\ \lf[\ov{\nabla}\lf[ f^{p-1}\, {\mathbb I}^0_{11}\rg] + (\ov{\nabla})^\perp\lf[ f^{p-1}\, {\mathbb I}^0_{12}\rg]\ \rg]\vec{n}\rg]\ dx^2\\[5mm]
\quad\rightarrow 0\quad.
\end{array}
\ee
It remains to bound
\be
\label{M-12}
\begin{array}{l}
\ds - \int_{{\R}_+}\ \chi(r)\ \frac{dr}{r^3}\sigma^2\  \int_{\vec{\Phi}^{-1}(B^4_r(\vec{q}))}(\vec{\Phi}-\vec{q})\cdot \Delta\lf[f^{p-1}\ \vec{H}\rg] \ dx^2\\[5mm]
\ds\quad =\int_{{\R}_+}\ \chi(r)\ \frac{dr}{r^3}\sigma^2\  \int_{\vec{\Phi}^{-1}(B^4_r(\vec{q}))}\nabla(\vec{\Phi}-\vec{q})\cdot \nabla\lf[f^{p-1}\ \vec{H}\rg] \ dx^2\\[5mm]
\ds\quad-\int_{{\R}_+}\ \chi(r)\ \frac{dr}{r^3}\sigma^2\  \int_{\vec{\Phi}^{-1}(\p B^4_r(\vec{q}))}(\vec{\Phi}-\vec{q})\cdot \p_\nu\lf[f^{p-1}\ \vec{H}\rg] \ dl\quad.
\end{array}
\ee
The last integral in the r.h.s. of (\ref{M-12}) is equal to
\be
\label{M-13}
\begin{array}{l}
\ds-\int_{{\R}_+}\ \chi(r)\ \frac{dr}{r^3}\sigma^2\  \int_{\vec{\Phi}^{-1}(\p B^4_r(\vec{q}))}(\vec{\Phi}-\vec{q})\cdot \p_\nu\lf[f^{p-1}\ \vec{H}\rg] \ dl\\[5mm]
\ds\quad= -\sigma^2\ \int_\Sigma\chi(|\vec{\Phi}-\vec{q}|)\ \nabla|\vec{\Phi}-\vec{q}|\cdot \lf<\nabla\lf[f^{p-1} \vec{H}\rg],\frac{\vec{\Phi}-\vec{q}}{|\vec{\Phi}-\vec{q}|^3}\rg>\ dx^2\quad.
\end{array}
\ee
We observe that since $\vec{H}\cdot\nabla\vec{\Phi}=0$
\[
 \lf<\nabla\lf[f^{p-1} \vec{H}\rg],\frac{\vec{\Phi}-\vec{q}}{|\vec{\Phi}-\vec{q}|^3}\rg>= \nabla\lf<f^{p-1} \vec{H},\frac{\vec{\Phi}-\vec{q}}{|\vec{\Phi}-\vec{q}|^3}\rg>+\, 3\, \lf<f^{p-1} \vec{H},\frac{\vec{\Phi}-\vec{q}}{|\vec{\Phi}-\vec{q}|^4}\rg>\ \nabla|\vec{\Phi}-\vec{q}|\quad.
\]
Hence, after integrating by parts we obtain from (\ref{M-13})
\be
\label{M-13-a}
\begin{array}{l}
\ds-\int_{{\R}_+}\ \chi(r)\ \frac{dr}{r^3}\sigma^2\  \int_{\vec{\Phi}^{-1}(\p B^4_r(\vec{q}))}(\vec{\Phi}-\vec{q})\cdot \p_\nu\lf[f^{p-1}\ \vec{H}\rg] \ dl\\[5mm]
\ds\quad= \sigma^2\ \int_\Sigma\chi(|\vec{\Phi}-\vec{q}|)\ \Delta|\vec{\Phi}-\vec{q}|\ \lf<f^{p-1} \vec{H},\frac{\vec{\Phi}-\vec{q}}{|\vec{\Phi}-\vec{q}|^3}\rg>\ dx^2\\[5mm]
\ds\quad+\,\sigma^2\ \int_\Sigma\lf[\chi'(|\vec{\Phi}-\vec{q}|)-3\, \frac{\chi(|\vec{\Phi}-\vec{q}|)}{|\vec{\Phi}-\vec{q}|}\rg]\ |\nabla|\vec{\Phi}-\vec{q}||^2\ \lf<f^{p-1} \vec{H},\frac{\vec{\Phi}-\vec{q}}{|\vec{\Phi}-\vec{q}|^3}\rg>\ dx^2\quad.
\end{array}
\ee
we observe that in the domain where $\chi(|\vec{\Phi}-\vec{q}|)\ne 0$ we have
\[
\Delta |\vec{\Phi}-\vec{q}|=\frac{(\vec{\Phi}-\vec{q})\cdot \Delta\vec{\Phi}}{|\vec{\Phi}-\vec{q}|}+\frac{|\nabla\vec{\Phi}|^2}{|\vec{\Phi}-\vec{q}|}-\frac{|\nabla|\vec{\Phi}-\vec{q}||^2}{|\vec{\Phi}-\vec{q}|}\quad,
\]
and using the fact that $\Delta\vec{\Phi}=-\vec{\Phi}\,|\nabla\vec{\Phi}|^2+\vec{H}\ |\nabla\vec{\Phi}|^2$ we finally obtain
\be
\label{M-14}
\Delta |\vec{\Phi}-\vec{q}|=-\frac{1-\vec{q}\cdot \vec{\Phi}}{|\vec{\Phi}-\vec{q}|}+\frac{(\vec{\Phi}-\vec{q})\cdot\vec{H}\, |\nabla\vec{\Phi}|^2}{|\vec{\Phi}-\vec{q}|}+\frac{|\nabla\vec{\Phi}|^2}{|\vec{\Phi}-\vec{q}|}-\frac{|\nabla|\vec{\Phi}-\vec{q}||^2}{|\vec{\Phi}-\vec{q}|}\quad.
\ee
Hence combining (\ref{M-13}), (\ref{M-13-a}) and (\ref{M-14}) we obtain
\be
\label{M-15}
\begin{array}{l}
\ds\lf|\int_{{\R}_+}\ \chi(r)\ \frac{dr}{r^3}\sigma^2\  \int_{\vec{\Phi}^{-1}(\p B^4_r(\vec{q}))}(\vec{\Phi}-\vec{q})\cdot \p_\nu\lf[f^{p-1}\ \vec{H}\rg] \ dl\rg|\\[5mm]
\ds\quad\quad\le C_\chi\ \sigma^2\ F_p(\vec{\Phi})+C_\chi\ \sigma^2 \ M(T)^{1/p}\ [F_p(\vec{\Phi})]^{1-1/p}\quad\rightarrow 0\quad.
\end{array}
\ee
Taking now the first integral in the r.h.s. of (\ref{M-12}) we have
\be
\label{M-16}
\begin{array}{l}
\ds\int_{{\R}_+}\ \chi(r)\ \frac{dr}{r^3}\sigma^2\  \int_{\vec{\Phi}^{-1}(B^4_r(\vec{q}))}\nabla(\vec{\Phi}-\vec{q})\cdot \nabla\lf[f^{p-1}\ \vec{H}\rg] \ dx^2\\[5mm]
\ds\quad\quad=-\int_{{\R}_+}\ \chi(r)\ \frac{dr}{r^3}\sigma^2\  \int_{\vec{\Phi}^{-1}(B^4_r(\vec{q}))}\Delta\vec{\Phi}\cdot f^{p-1}\ \vec{H} \ dx^2\\[5mm]
\ds\quad\quad+\sigma^2\ \int_\Sigma \chi(|\vec{\Phi}-\vec{q}|)\ \nabla|\vec{\Phi}-\vec{q}|\cdot \lf<\nabla(\vec{\Phi}-\vec{q}),f^{p-1}\ \vec{H}\rg>\ dx^2\quad.
\end{array}
\ee
So we have also
\be
\label{M-17}
\begin{array}{l}
\ds\lf|\int_{{\R}_+}\ \chi(r)\ \frac{dr}{r^3}\sigma^2\  \int_{\vec{\Phi}^{-1}(B^4_r(\vec{q}))}\nabla(\vec{\Phi}-\vec{q})\cdot \nabla\lf[f^{p-1}\ \vec{H}\rg] \ dx^2\rg|\\[5mm]
\ds\quad\quad\le C_\chi\ \sigma^2\ F_p(\vec{\Phi})+C_\chi\ \sigma^2 \ M(T)^{1/p}\ [F_p(\vec{\Phi})]^{1-1/p}\quad\rightarrow 0\quad.
\end{array}
\ee
Combining (\ref{M-12}) and (\ref{M-15}) and (\ref{M-17}) we have
\be
\label{M-18}
\lf|\int_{{\R}_+}\ \chi(r)\ \frac{dr}{r^3}\sigma^2\  \int_{\vec{\Phi}^{-1}(B^4_r(\vec{q}))}(\vec{\Phi}-\vec{q})\cdot \Delta\lf[f^{p-1}\ \vec{H}\rg] \ dx^2\rg|\rightarrow 0\quad.
\ee
Combining now (\ref{M-2}), (\ref{M-3}), (\ref{M-4}), (\ref{M-11}) and (\ref{M-18}) we have that for any fixed non negative $\chi(r)$ compactly supported in ${\R}^\ast_+$ and any $\vec{q}\in {\R}^4$
\be
\label{M-19}
\begin{array}{l}
\ds-\int_0^\infty\chi'(r)\ dr\ \frac{1}{r^2}\int_{\vec{\Phi}_k^{-1}(B^4_r(\vec{q}))}\ dvol_{g_{\vec{\Phi}_k}}\ge\\[5mm]
\ds\quad\quad -C\ \int_0^\infty\chi(r)\ dr\ \frac{1}{r}\int_{\vec{\Phi}_k^{-1}(B^4_r(\vec{q}))}\ dvol_{g_{\vec{\Phi}_k}}-C_\chi\ \sigma^2\ F_p(\vec{\Phi})-C_\chi\ \sigma^2 \ M(T)^{1/p}\ [F_p(\vec{\Phi})]^{1-1/p}\quad,
\end{array}
\ee
for some constant $C_\chi$ depending on $\chi$.
Taking $\mu_k$ the Radon measure on ${\R}^4$ given by (\ref{VII.3-a}) this can be rewritten as
\[
-\int_0^\infty\chi'(r)\ dr\ \frac{1}{r^2}\mu_k(B^4_r(\vec{q}))\ge -C\ \int_0^\infty\chi(r)\ dr\ \frac{1}{r}\mu_k(B^4_r(\vec{q}))+o_k(1)\quad.
\]
We extract a subsequence such that $\mu_k$ converges weakly in Radon measure and we finally obtain that for any fixed non negative $\chi(r)$ compactly supported in ${\R}^\ast_+$ and any $\vec{q}\in {\R}^4$
\[
-\int_0^\infty\chi'(r)\ dr\ \frac{1}{r^2}\mu_\infty(B^4_r(\vec{q}))\ge -C\ \int_0^\infty\chi(r)\ dr\ \frac{1}{r}\mu_\infty(B^4_r(\vec{q}))\quad,
\]
which classically implies (\ref{VII.4}) and lemma~\ref{lm-VII.1} is proved. \hfill $\Box$

\medskip
A rather direct consequence of the proof of the limiting monotonicity formula is given by the following non concentration result.
\begin{Lm}
\label{lm-concent}{\bf [Non Collapsing Lemma]}
Let $p>1$ and $0<\sigma<1$.  There exists $\delta>0$ and $\ep>0$ such that for any critical point $\vec{\Phi}$ of $A^\sigma$ satisfying
\be
\label{M-19-a}
 \sigma^2 F_p(\vec{\Phi})\le \ep\ \mbox{Area}(\vec{\Phi})\quad,
\ee
then
\be
\label{M-19-b}
\mbox{diam}(\vec{\Phi}(\Sigma))>\delta\quad.
\ee
\hfill $\Box$
\end{Lm}
\noindent{\bf Proof of lemma~\ref{lm-concent}.} Assume (\ref{M-19-a}) for some $\ep$ fixed later. Let $1>\delta>0$ and choose $\chi_\delta= (r-\delta)^+$ on $[0,1+\delta]$ identically equal to $1$ on $[1+\delta, 2+\delta]$ and equal to $(3+\delta-r)^+$ for $r>2+\delta$. Assuming that the whole immersed surface is included in a ball $B^4_\delta(\vec{q})$, the inequality (\ref{M-19}) gives then
\be
\label{M-19-c}
-\, \mbox{Area}(\vec{\Phi})\, \lf[\int_\delta^{1+\delta} \frac{dr}{r^2}-\int_{2+\delta}^{3+\delta} \frac{dr}{r^2}\rg]\ge -\, C\ \mbox{Area}(\vec{\Phi}) \int_\delta^{3+\delta} \frac{dr}{r}- C_\delta\,\ep^{1-1/p} \ \mbox{Area}(\vec{\Phi}) \quad.
\ee
Dividing by $\mbox{Area}(\vec{\Phi}) $ we obtain
\[
C\ \log\frac{1}{\delta}\ge \frac{1}{\delta}-\frac{1}{4}- C_\delta\,\ep^{1-1/p}\quad.
\]
Assume that $\delta$ is small enough in such a way that $C\ \log\frac{1}{\delta}< \frac{1}{\delta}-{1}$, choosing $\ep>0$ such that  $C_\delta\,\ep^{1-1/p}<3/4$ we obtain a contradiction.
This proves lemma~\ref{lm-concent}.\hfill $\Box$

\medskip

 The next result establishes  a uniform lower bound of the limiting area for any sequence of immersions satisfying the assumptions of theorem~\ref{th-limit}. This result is the ``work-horse" in our proof of the main theorem and shall be used crucially at several steps. 
Precisely we have the following result
\begin{Lm}
\label{lm-quanti} {\bf[Global Energy Quantization]}
Let $p>1$. For every $\La>0$ there exists $Q_0(\La)>0$ and $\sigma(\La)>0$ such that the following holds. Let $\Sigma$ be a closed surface and let $\vec{\Phi}$ be a  critical points of $A_p^{\sigma}$ for $\sigma<\sigma(\La)$ and satisfying
\be
\label{VII.2-aa}
\sigma^2\ F_p(\vec{\Phi})=\sigma^2\ \int_{\Sigma}[1+|{\mathbb I}_{\vec{\Phi}}|^{2}_{g_{\vec{\Phi}}}]^p\ dvol_{g_{\vec{\Phi}}}\le \frac{\La}{\log(1/\sigma)} \mbox{Area}(\vec{\Phi})\quad.
\ee
then, 
\be
\label{VII.2-ab}
 \mbox{Area}(\vec{\Phi})\ge Q_0(\La)\quad.
\ee
\end{Lm}
\noindent{\bf Proof of lemma~\ref{lm-quanti}.} We denote as usual
\[
f(\sigma)=\frac{\sigma^2\ F_p(\vec{\Phi})}{\mbox{Area}(\vec{\Phi})}\quad.
\]
Let $\eta>0$ to be fixed later. For any $\vec{q}\in\vec{\Phi}(\Sigma)$ we consider the 4-dimensional ball in ${\R}^4$, $B^4_\sigma(\vec{q})$ centered at $\vec{q}$ with radius $\sigma$. We consider the subset $E_\eta$ of $\vec{\Phi}(\Sigma)$ given by
\[
E_\eta:=\lf\{\vec{q}\in\vec{\Phi}(\Sigma)\subset S^3\quad;\quad\sigma^{-2}\int_{B^4_\sigma(\vec{q})\cap\vec{\Phi}(\Sigma)}\ dvol_{g_{\vec{\Phi}}}<\eta\rg\}\quad.
\]
From the covering $(B^4_\sigma(\vec{q}))_{\vec{q}\in E_\eta}$ we extract a Besicovitch sub-covering $(B^4_\sigma(\vec{q}_i))_{i\in I}$ such that each point in ${\R}^4$ is covered by at most $N$ balls where $N$ is a universal number. A corollary of Simon's monotonicity formula (see corollary 5.12 \cite{Ri1} and take $T=\sigma$) gives for each $i\in I$
\be
\label{VII.201}
\sigma^{-2}\int_{B^4_\sigma(\vec{q})\cap\vec{\Phi}(\Sigma)}\ dvol_{g_{\vec{\Phi}}}\ge \frac{2\pi}{3}-\frac{1}{2}\int_{B^4_\sigma(\vec{q}_i)}|\vec{H}_{\vec{\Phi}}^{{\R}^4}|^2\ dvol_{g_{\vec{\Phi}}}\quad.
\ee
Considering $\eta=\pi/3$ this imposes
\be
\label{VII.202}
\int_{B^4_\sigma(\vec{q}_i)}|\vec{H}_{\vec{\Phi}}^{{\R}^4}|^2\ dvol_{g_{\vec{\Phi}}}>\frac{2\,\pi}{3}\quad.
\ee
Hence
\be
\label{VII.203}
\ds\int_{\cup_{i\in I}B^4_\sigma(\vec{q}_i)}|\vec{H}_{\vec{\Phi}}^{{\R}^4}|^2\ dvol_{g_{\vec{\Phi}}}\ge \frac{1}{N}\sum_{i\in I}\int_{B^4_\sigma(\vec{q}_i)}|\vec{H}_{\vec{\Phi}}^{{\R}^4}|^2\ dvol_{g_{\vec{\Phi}}}\ge \frac{2\,\pi}{3\,N}\ \mbox{card}{I}\quad.
\ee
Combining (\ref{VII.2-aa}) and (\ref{VII.203}) we obtain
\be
\label{VII.204}
\sigma^2\ \frac{2\,\pi}{3\,N}\ \mbox{card}{I}\le f(\sigma)\ \mbox{Area}(\vec{\Phi})\quad.
\ee
So we have
\be
\label{VII.205}
\int_{E_{\pi/3}}\ dvol_{g_{\vec{\Phi}}}\le\int_{\cup_{i\in I}B^4_\sigma(\vec{q}_i)}\ dvol_{g_{\vec{\Phi}}}\le \frac{\pi}{3}\ \sigma^2 \ \mbox{card}{I}\le f(\sigma)\ \mbox{Area}(\vec{\Phi})\quad.
\ee
Let $1>\delta>0$ to be fixed later . Consider now for $j\in \{1,2\cdots \log_2 \sigma^{-1}\}$. We use the notation
\[
A(j,\vec{q}):=\int_{B^4_{2^{j} \sigma}(\vec{q})\cap\vec{\Phi}(\Sigma)}\ dvol_{g_{\vec{\Phi}}}\quad\mbox{and}\quad F(j,\vec{q}):=\sigma^2\,\int_{B^4_{2^j\,\sigma}(\vec{q})\cap\vec{\Phi}(\Sigma)}\ \lf[1+|{\mathbb I}_{\vec{\Phi}}|^{2}_{g_{\vec{\Phi}}}\rg]^p\ dvol_{g_{\vec{\Phi}}}\quad,
\]
\[
G^j_\delta:=\lf\{
\begin{array}{l}
\ds\vec{q}\in\vec{\Phi}(\Sigma)\setminus E_{\pi/3}\ ;\  \frac{(2^{-2\,j}\ A(j+1,\vec{q}))^{1/p}\ F(j+1,\vec{q})^{1-1/p}+F(j,\vec{q})}{A(j,\vec{q})}\ge \frac{f(\sigma)}{\delta}\\[3mm]
\ds \mbox{and }\quad A(j+1,\vec{q})\le 3\,\pi\, 2^{2j+2}\ \sigma^2\quad.
\end{array}
\rg\}
\]
For each $j\in \{1,2,\cdots ,\log_2 \sigma^{-1}\}$ and for any $\vec{q}\in G^j_\delta$ we consider the closed balls $B^4_{2^j \sigma}(\vec{q})$. The following  covering of $G_\delta:=\cup_{j\in \{1,2\cdots \log_2 \sigma^{-1}-1\}} G_\delta^j$
\[
\lf((B^4_{2^j \sigma}(\vec{q}))_{\vec{q}\in G^j_\delta})\rg)_{j=1,2,\cdots ,\log_2 \sigma^{-1}}
\] 
realizes a Besicovitch covering of $G_\delta$. By the mean of Besicovitch theorem, we extract  a Besicovitch sub-covering
 \[
\lf((B^4_{2^{j} \sigma}(\vec{q}_i))_{i\in I_j}\rg)_{j=1\cdots \log_2 \sigma^{-1}}
\]
 of $G_\delta$ such that each point in ${\R}^4$ is covered by at most ${\mathcal N}$ balls   where ${\mathcal N}$ is a universal number\footnote{Observe that it is not clear whether for each $j$ the sub-familly 
$(B^4_{2^{j} \sigma}(\vec{q}_i))_{i\in I_j}$ covers the whole $G_\delta^j$ but at least the union of these families cover $G_\delta$.} . In other words we have
 \be
 \label{rep-1}
 \lf\|\sum_{j=1}^{\log_2\sigma^{-1}-1}\sum_{i\in I_j} {\mathbf 1}_{B^4_{2^{j} \sigma}(\vec{q}_i)}\ \rg\|_{L^\infty({\R}^4)}\le\ {\mathcal N}\quad.
 \ee
 For any ${j=1\cdots \log_2 \sigma^{-1}}$
 the balls $B^4_{2^{j} \sigma}(\vec{q}_i)$ for $i\in I_j$ have all the same radius, moreover each point of ${\R}^4$ is covered by at most ${\mathcal N}$ of such balls. Hence by doubling each of these
 balls and considering $B^4_{2^{j+1} \sigma}(\vec{q}_i)$, since they all have the \underbar{same radius} there exists a universal number\footnote{Observe that a-priori each point of ${\R}^4$ can be covered by at most
 ${\mathfrak N}\, \log_2\sigma^{-1}$ of the double balls $\lf((B^4_{2^{j+1} \sigma}(\vec{q}_i))_{i\in I_j}\rg)_{j=1\cdots \log_2 \sigma^{-1}}$.}   ${\mathfrak N}$ such that
\[
\sup_{j=1\cdots \log_2 \sigma^{-1}}\lf\|\sum_{i\in I_j} {\mathbf 1}_{B^4_{2^{j+1} \sigma}(\vec{q}_i)}\rg\|_{L^\infty({\R}^4)}\le\ {\mathfrak N}\quad,
\]
where $ {\mathbf 1}_{B^4_{2^{j+1} \sigma}(\vec{q}_i)} $ is the characteristic function of the ball $B^4_{2^{j+1} \sigma}(\vec{q}_i)$. We have for any $\al>0$ that
\be
\label{VII.205-ab}
\lf\|\sum_{j=1}^{\log_2\sigma^{-1}-1}\sum_{i\in I_j} {\mathbf 1}_{B^4_{2^{j+1} \sigma}(\vec{q}_i)}\ 2^{\al j}\rg\|_{L^\infty({\R}^4)}\le C\ {\mathfrak N}\ \sum_{j=0}^{\log_2\sigma^{-1}}2^{\al\, j}\le C\ {\mathfrak N}\ \sigma^{-\al}\quad,
\ee
 For any  $j\in\{1,2\cdots \log_2 \sigma^{-1}\}$ and $\vec{q}\in G_\delta^j$, the whole support of $\vec{\Phi}(\Sigma)$ cannot be included in $B^4_{2^j\,\sigma}(\vec{q})$ otherwise we would contradict
 the non collapsing lemma~\ref{lm-concent} for $\sigma$ small enough. Hence, since $\vec{q}\in\vec{\Phi}({\Sigma})$ for any radius $r\in (2^j\sigma,2^{j+1}\sigma)$ 
 we have $\vec{\Phi}(\Sigma)\cap \p B_r(\vec{q})\ne\emptyset$ and we can apply lemma~\ref{lm-A.2}. Hence we deduce
 \be
\label{VII.205-b}
0<\ep_0(4)<\int_{B^4_{2^{j+1} \sigma}(\vec{q})} |\vec{\mathbb I}^{\,{\R}^4}_{\vec{\Phi}}|^{2}\ dvol_{g_{\vec{\Phi}}}\quad.
\ee
Since $A(j+1,\vec{q})\le 3\,\pi\, 2^{2j+2}\ \sigma^2$ inequality (\ref{VII.205-b}) implies
\be
\label{VII.205-c}
\begin{array}{l}
\ds\frac{A(j+1,\vec{q})}{2^{2\,j+2}}\le \frac{3\,\pi\ \sigma^2}{\ep_0(4)}\ \int_{B^4_{2^{j+1} \sigma}(\vec{q})} |\vec{\mathbb I}^{\,{\R}^4}_{\vec{\Phi}}|^{2}\ dvol_{g_{\vec{\Phi}}}\\[3mm]
\ds\quad\quad\le \frac{3\,\pi\ \sigma^2}{\ep_0(4)}\ A(j+1,\vec{q})^{1-1/p}\ \lf(\int_{B^4_{2^{j+1} \sigma}(\vec{q})}[1+ |{\mathbb I}_{\vec{\Phi}}|^{2}]^p\ dvol_{g_{\vec{\Phi}}}\rg)^{1/p}
\end{array}
\ee
 and 
we deduce  that
\be
\label{VII.205-d}
\frac{A(j+1,\vec{q})}{2^{2\,j}}\le C\ (2^{j+1}\,\sigma)^{2p-2} \ F(j+1,\vec{q})\quad.
\ee
So for $\vec{q}\in G_\delta^j$ we have combining the definition of $G_\delta^j$ with (\ref{VII.205-d})
\be
\label{VII.205-e}
\frac{f(\sigma)}{\delta}\ A(j,\vec{q})\le F(j,\vec{q})+\, C\ (2^{j+1}\,\sigma)^{2-2/p}\ F(j+1,\vec{q})
\ee
summing this identity with respect to $j\in J$ we obtain
\be
\label{VII.206-bb}
\begin{array}{l}
\ds\frac{f(\sigma)}{\delta}\int_{G_\delta}\ dvol_{g_{\vec{\Phi}}}\le \frac{f(\sigma)}{\delta}\sum_{j=1}^{\log_2\sigma^{-1}}\int_{G^j_\delta}\ dvol_{g_{\vec{\Phi}}}    \le\frac{f(\sigma)}{\delta}\ \sum_{j=1}^{\log_2\sigma^{-1}}\sum_{i\in I_j}\int_{B^4_{2^{j} \sigma}(\vec{q}_i)}\ dvol_{g_{\vec{\Phi}}}\\[5mm]
\ds\quad\quad\le \sum_{j=1}^{\log_2\sigma^{-1}}\sum_{i\in I_j}\sigma^2\ \int_{B^4_{2^{j_i} \sigma}(\vec{q}_i)}\ \lf[1+|{\mathbb I}_{\vec{\Phi}}|^{2}_{g_{\vec{\Phi}}}\rg]^p\  dvol_{g_{\vec{\Phi}}}\\[5mm]
\ds\quad\quad+\sigma^2\ \int_{\Sigma} \sum_{j=1}^{\log_2\sigma^{-1}}\sum_{i\in I_j} {\mathbf 1}_{B^4_{2^{j+1} \sigma}(\vec{q}_i)}\ 2^{\al j}\ \sigma^\al \lf[1+|{\mathbb I}_{\vec{\Phi}}|^{2}_{g_{\vec{\Phi}}}\rg]^p\  dvol_{g_{\vec{\Phi}}}\quad,
\end{array}
\ee
where $\al:=2-2p$. Using now (\ref{rep-1}) and (\ref{VII.205-ab}), we then deduce
\be
\label{VII.206}
\begin{array}{l}
\ds\frac{f(\sigma)}{\delta}\int_{G_\delta}\ dvol_{g_{\vec{\Phi}}}\le C\,\sigma^2\ \int_\Sigma \lf[1+|{\mathbb I}_{\vec{\Phi}}|^{2}_{g_{\vec{\Phi}}}\rg]^p\  dvol_{g_{\vec{\Phi}}}=\, C\, f(\sigma)\ \int_\Sigma dvol_{g_{\vec{\Phi}}}\quad.
\end{array}
\ee
We deduce from (\ref{VII.205}) and (\ref{VII.206})
\be
\label{VII.207}
\int_{E_{\pi/3}\cup G_\delta}\ dvol_{g_{\vec{\Phi}}}\le (C\ \delta +f(\sigma))\ \int_{\Sigma}\ dvol_{g_{\vec{\Phi}}}\quad.
\ee
Since $f(\sigma)\rightarrow 0$ as $\sigma\rightarrow 0$, by taking any $0<\delta<1/C$ we have that for $\sigma$ small enough $\vec{\Phi}(\Sigma)\setminus(E_{\pi/3}\cup G_\delta)\ne\emptyset$. Let now $\vec{q}\in \vec{\Phi}(\Sigma)\setminus(E_{\pi/3}\cup G_\delta)$. Take  $j_0=j({\vec{q}})$ the \underbar{largest} index such that
\[
\ds\int_{B^4_{2^{j_0}\sigma}(\vec{q})}\ dvol_{g_{\vec{\Phi}}}\ge 2^{2\,j_0}\sigma^2\ \pi/3\ \quad.
\]
Since $\vec{q}\in \vec{\Phi}(\Sigma)\setminus(E_{\pi/3}\cup G_\delta)$ we must have
\be
\label{VII.208}
\lf\{
\begin{array}{l}
\ds\int_{B^4_{2^{j_0}\sigma}(\vec{q})}\ dvol_{g_{\vec{\Phi}}}\ge 2^{2\,j_0}\sigma^2\ \pi/3\quad\quad\mbox{ and }\\[5mm]
\ds\forall j\ge j_0\quad\quad \frac{f(\sigma)}{\delta}\ \int_{B^4_{2^j \sigma}(\vec{q})\cap\vec{\Phi}(\Sigma)}\ dvol_{g_{\vec{\Phi}}}\ge \sigma^2\,\int_{B^4_{2^j\,\sigma}(\vec{q})\cap\vec{\Phi}(\Sigma)}\ \lf[ 1+|{\mathbb I}_{\vec{\Phi}}|^{2}_{g_{\vec{\Phi}}}\rg]^p\ dvol_{g_{\vec{\Phi}}}\\[5mm]
\ds+\, \lf[\sigma^2\,\int_{B^4_{2^{j+1}\,\sigma}(\vec{q})\cap\vec{\Phi}(\Sigma)}\ \lf[1+|{\mathbb I}_{\vec{\Phi}}|^{2}_{g_{\vec{\Phi}}}\rg]^p\ dvol_{g_{\vec{\Phi}}}\rg]^{1-1/p}\ \lf[ 2^{-2\,j}\int_{B^4_{2^{j+1} \sigma}(\vec{q})\cap\vec{\Phi}(\Sigma)}\ dvol_{g_{\vec{\Phi}}}\rg]^{1/p}\quad.
\end{array}
\rg.
\ee
Let $j\in \{j_0,\cdots,\log_2\sigma^{-1}-1\}$ and let $\chi$ be an arbitrary smooth function, bounded by 1, supported in $[2^{j-2}\sigma,2^{j+1}\sigma]$ and such that $|\chi'|\le\ C\ 2^{-j}\sigma^{-1}$. We can estimate each error terms between (\ref{M-1}) and (\ref{M-19})  in the computations of the monotonicity formula at fixed $k$ between (\ref{M-1}) and (\ref{M-19}) by the mean of the area  we obtain
\be
\label{VII.209}
\begin{array}{l}
\ds-\int_0^{+\infty}\chi' (r)\frac{dr}{r^2}\int_{\vec{\Phi}^{-1}(B^4_r(\vec{q}))}\ dvol_{g_{\vec{\Phi}}}\ge
\ds\quad -\  C\ \int_0^{+\infty}\chi(r)\ dr\,\lf[\frac{1}{r}+\frac{o_\sigma(1)}{r^2}\rg]\int_{\vec{\Phi}^{-1}(B^4_r(\vec{q}))}\ dvol_{g_{\vec{\Phi}}}\\[5mm]
\ds-\ C\ \int_0^{+\infty}\chi(r)\frac{dr}{r^3}\int_{\vec{\Phi}^{-1}(B^4_r(\vec{q}))}\ \sigma^2\ \lf[1+|{\mathbb I}_{\vec{\Phi}}|^{2}_{g_{\vec{\Phi}}}\rg]^p\ dvol_{g_{\vec{\Phi}}}\\[5mm]
\ds-\ C\, 2^{-3 j}\,\sigma^{-3} \lf[\sigma^2\,\int_{B^4_{2^{j+1}\,\sigma}(\vec{q})\cap\vec{\Phi}(\Sigma)}\ \lf[1+|{\mathbb I}_{\vec{\Phi}}|^{2}_{g_{\vec{\Phi}}}\rg]^p\ dvol_{g_{\vec{\Phi}}}\rg]^{1-1/p}\ \lf[ 2^{-2\,j}\int_{B^4_{2^{j+1} \sigma}(\vec{q})\cap\vec{\Phi}(\Sigma)}\ dvol_{g_{\vec{\Phi}}}\rg]^{1/p}\quad.
\end{array}
\ee
 Using (\ref{VII.208}) we deduce that for any $r\in [2^{j_0}\,\sigma, 1/2]$
\be
\label{VII.209-a}
\begin{array}{l}
\ds\frac{d}{dr}\lf[ \frac{1}{r^2}\int_{\vec{\Phi}^{-1}(B^4_r(\vec{q}))}\ dvol_{g_{\vec{\Phi}}}\rg]\ge
\ds\quad -\  \lf[\frac{C}{r}+\frac{o_\sigma(1)}{r^2}\rg]\int_{\vec{\Phi}^{-1}(B^4_r(\vec{q}))}\ dvol_{g_{\vec{\Phi}}}\\[5mm]
\ds\quad \quad\quad\quad-\,C\ \frac{f(\sigma)}{ \delta}\  \frac{1}{r^3}\int_{\vec{\Phi}^{-1}(B^4_r(\vec{q}))}\ dvol_{g_{\vec{\Phi}}}\quad.
\end{array}
\ee
Let $Y(r):=\frac{1}{r^2}\int_{\vec{\Phi}^{-1}(B^4_r(\vec{q}))}\ dvol_{g_{\vec{\Phi}}}$, this ordinary differential inequality gives, for $\sigma$ small enough, the existence of $C>0$ independent of $r$ and $\sigma$ and $\delta$, such that for $r\in  [2^{j_0}\sigma, 1/2]$
\be
\label{VII.210}
\frac{d}{dr}\lf[e^{C\, r} r^{\frac{C\,f(\sigma)}{ \delta}} Y\rg]\ge 0\quad.
\ee
Integrating between $2^{j_0}\,\sigma$ and $1/2$ gives
\[
e^{C/2}\ Y(1/2)\ 2^{-\frac{C\,f(\sigma)}{ \delta}}\ge e^{C\,2^{j_0}\,\sigma} (2^{j_0}\,\sigma)^{\frac{C\,f(\sigma)}{ \delta}} Y(2^{j_0}\,\sigma)\quad.
\]
Using the fact that $\vec{q}\in\vec{\Phi}(\Sigma)\setminus E_{\pi/3}$ we have then using the first line in (\ref{VII.208})
\be
\label{VII.211}
Y(1/2)\ge e^{-C/2}\ 2^{\frac{3\,C\,f(\sigma)}{ \pi}}\ e^{C\,2^{j_0}\,\sigma} (2^{j_0}\,\sigma)^{\frac{3\,C\,f(\sigma)}{ \pi}}\ \frac{\pi}{3}\quad.
\ee
Since $f(\sigma)\,\log_2\sigma^{-1})\le \La$ we have $(2^{j_0}\,\sigma)^{\frac{C\,f(\sigma)}{ \delta}}= 2^{ C\, f(\sigma)\, \delta^{-1}\ \log_2(2^{j_0}\,\sigma)} \ge 2^{ C\, f(\sigma)\, \delta^{-1}\ \log_2\sigma}\ge 2^{-\,C\,\delta^{-1}\,\La} $ . So $Q_0:= 2^{-\,C\,\delta^{-1}\,\La}\,e^{-C/2}\ \pi/3$ satisfies (\ref{VII.2-ab}) and the lemma~\ref{lm-quanti} is proved.\hfill $\Box$

\medskip

We now introduce two definitions. First we define the {\it Oscillation set}.
\begin{Dfi}
\label{df-concentration}
Let $\vec{\Phi}_k$ be a sequence of conformal smooth immersions from\footnote{Recall that in this section we are assuming that the underlying conformal class to $(\Sigma,g_k)$ is 
precompact in the moduli space.} $(\Sigma,g_k)$, critical points of $$A^{\sigma_k}_p(\vP):=\mbox{Area}(\vec{\Phi})+\sigma_k^2\, F_p(\vec{\Phi})=\int_\Sigma \lf[1+\sigma_k^2\ [1+\gI]^p\rg]\ dvol_{g_{\vec{\Phi}}}$$ in the space of
weak immersions into $S^3$ and for $\sigma_k\rightarrow 0$. Assume 
\[
\vec{\Phi}_k\rightharpoonup\vec{\Phi}_\infty\quad\quad\mbox{weakly in }W^{1,2}(\Sigma,S^3)\quad,
\]
where $\Sigma$ is equipped with a reference metric $g_0$. Assume the sequence of Riemann surfaces $(\Sigma,g_{\vec{\Phi}_k})$ is pre-compact in the moduli space of conformal structures on $\Sigma$ and assume
\[
\nu_{k}:=|d\vec{\Phi}_{k}|_{h_{k}}^2\ dvol_{h_{k}}=|\nabla\vec{\Phi}_{k}|^2\ dx^2\rightharpoonup \nu_\infty\quad\mbox{ in Radon measures}
\]
The oscillation set ${\mathcal O}\subset\Sigma$ is the set of points $x\in\Sigma$ such that
\be
\label{M-29-aa}
{\mathcal O}:=\lf\{
\begin{array}{c}
\ds x\in\Sigma\ ;\ \quad\nu_\infty(B_\rho(x))\ne 0\quad\forall\ \rho>0\\[5mm]
\ds\mbox{and }\quad\liminf_{\rho\rightarrow 0}\frac{\int_{B_{2\rho}(x)}|d\vec{\Phi}_\infty|_{g_0}^2\ dvol_{g_0} }{\nu_\infty(\ov{B_\rho(x)})}=0
\end{array}
\rg\}\quad.
\ee

\hfill $\Box$
\end{Dfi}
Now we define the {\it vanishing set} ${\mathcal V}$.
\begin{Dfi}
\label{df-vanish}
Let $\vec{\Phi}_k$ be a sequence of conformal smooth immersions from $(\Sigma,g_k)$, critical points of $$A^{\sigma_k}_p(\vP):=\mbox{Area}(\vec{\Phi})+\sigma_k^2\, F_p(\vec{\Phi})=\int_\Sigma \lf[1+\sigma_k^2\ [1+\gI]^p\rg]\ dvol_{g_{\vec{\Phi}}}$$ in the space of
weak immersions into $S^3$ and for $\sigma_k\rightarrow 0$.  We assume $(\Sigma,g_k)$ to be pre-compact in the moduli space of conformal structures on $\Sigma$. Denote
\be
\label{M-45-bbb}
f(\sigma_k):=\frac{\ds\sigma_k^2\,\int_\Sigma\  \lf[1+|{\mathbb I}_{\vec{\Phi}_k}|^{2}_{g_{\vec{\Phi}_k}}\rg]^p\  dvol_{g_{\vec{\Phi}_k}}}{\ds\int_\Sigma\ dvol_{g_{\vec{\Phi}_k}}}\quad.
\ee
We call the ''vanishing set'' the subset $\Sigma_0$ of $\Sigma$ given by
\be
\label{M-45-aa}
\Sigma_0:=\lf\{x\in \Sigma\quad;\quad\liminf_{r\rightarrow 0}\ \limsup_{k\rightarrow +\infty}\frac{\ds f(\sigma_k)\int_{B_r(x)}\ dvol_{g_{\vec{\Phi}_k}}}{\ds\sigma_k^2\,\int_{B_r(x)}\  \lf[1+|{\mathbb I}_{\vec{\Phi}_k}|^{2}_{g_{\vec{\Phi}_k}}\rg]^p\  dvol_{g_{\vec{\Phi}_k}}}=0\rg\}\quad.
\ee
\hfill $\Box$
\end{Dfi}
We will need later on the following lemma which justifies the denomination {\it vanishing set}.
\begin{Lm}
\label{lm-vanish}{\bf[No Limiting Measure on the Vanishing Set]}
Let $\vec{\Phi}_k$ be a sequence of conformal smooth immersions from $(\Sigma,g_k)$ into $S^3$, critical points of $$A^{\sigma_k}_p(\vP):=\mbox{Area}(\vec{\Phi})+\sigma_k^2\, F_p(\vec{\Phi})=\int_\Sigma \lf[1+\sigma_k^2\ [1+\gI]^p\rg]\ dvol_{g_{\vec{\Phi}}}$$ in the space of
weak immersions into $S^3$  for $\sigma_k\rightarrow 0$. We assume $(\Sigma,g_k)$ is strongly pre-compact in the Moduli space of $\Sigma$. Assume 
\[
\vec{\Phi}_k\rightharpoonup\vec{\Phi}_\infty\quad\quad\mbox{weakly in }W^{1,2}(\Sigma,S^3)\quad,
\]
and assume the following sequence of Radon measure weakly converges
\[
\nu_{k}:=|d\vec{\Phi}_{k}|_{g_{k}}^2\ dvol_{g_{k}}\ \rightharpoonup \nu_\infty\quad,
\]
then we have
\be
\label{M-45-555}
\nu_\infty(\Sigma_0)=0\quad.
\ee
\hfill $\Box$
\end{Lm}
\noindent{\bf Proof of lemma~\ref{lm-vanish}.}
We have
\be
\label{M-46-aa}
\begin{array}{l}
\ds\forall\ x\in \Sigma_0\quad\forall\, \delta>0\quad \exists\, k_{x,\delta}\in{\N}\quad\exists\,r_{x,\delta}>0\\[5mm]
\ds\quad\quad\mbox{s. t.}\quad \quad\forall\, k\ge k_{x,\delta}\quad\quad\quad\frac{\ds f(\sigma_k)\int_{B_{r_x}(x)}\ dvol_{g_{\vec{\Phi}_k}}}{\ds\sigma_k^2\,\int_{B_{r_x}(x)}\  \lf[1+|{\mathbb I}_{\vec{\Phi}_k}|^{2}_{g_{\vec{\Phi}_k}}\rg]\  dvol_{g_{\vec{\Phi}_k}}}<\delta\quad.
\end{array}
\ee
For any $\delta>0$ and $j\in {\N}$ we denote
\[
\Sigma_0^j(\delta):=\lf\{x\in\Sigma_0\quad;\quad k_{x,\delta}\le j\rg\}
\]
We have clearly $\Sigma_0=\cup_{j\in {\N}}\Sigma_0^j(\delta)$. From the covering $(\ov{B_{r_{x,\delta}}(x)})_{x\in \Sigma_0^j(\delta)}$ we extract a Besicovitch sub-covering of $\Sigma_0^j(\delta)$ that we denote $(\ov{B_{r_{x_i,\delta}}(x_i)})_{i\in I}$ in such a way that
any point of $\Sigma$ is covered by at most $N$ balls from this sub-covering. We have for all $k\ge j$
\[
\int_{B_{r_{x_i}(x_i)}}\ dvol_{g_{\vec{\Phi}_k}}\le \frac{\delta}{f(\sigma_k)}\ \sigma_k^2\,\int_{B_{r_{x_i}}(x_i)}\  \lf[1+|{\mathbb I}_{\vec{\Phi}_k}|^{2}_{g_{\vec{\Phi}_k}}\rg]^p\  dvol_{g_{\vec{\Phi}_k}}\quad.
\]
Summing over $i\in I$ gives
\be
\label{M-47-aa}
\begin{array}{l}
\ds\nu_k\lf(\bigcup_{i\in I}\ov{B_{r_{x_i}(x_i)}}\rg)\le \sum_{i\in I}\int_{B_{r_{x_i}(x_i)}}\ dvol_{g_{\vec{\Phi}_k}}\le \frac{\delta}{f(\sigma_k)}\ \sigma_k^2\, \sum_{i\in I}\int_{B_{r_{x_i}}(x_i)}\  \lf[1+|{\mathbb I}_{\vec{\Phi}_k}|^{2}_{g_{\vec{\Phi}_k}}\rg]^p\  dvol_{g_{\vec{\Phi}_k}}\\[5mm]
\ds\quad\le N\ \frac{\delta}{f(\sigma_k)}\ \sigma_k^2\int_{\cup_{i\in I}B_{r_{x_i}(x_i)}}\lf[1+|{\mathbb I}_{\vec{\Phi}_k}|^{2}_{g_{\vec{\Phi}_k}}\rg]^p\  dvol_{g_{\vec{\Phi}_k}}\le N\ \delta \int_\Sigma\ dvol_{g_{\vec{\Phi}_k}}\quad.
\end{array}
\ee
This implies that
\be
\label{M-48-aa}
\nu_\infty(\Sigma_0^j(\delta))\le \limsup_{k\rightarrow +\infty}\nu_k\lf(\bigcup_{i\in I}\ov{B_{r_{x_i}(x_i)}}\rg)\le N\ \delta\ \nu_\infty(\Sigma)\quad.
\ee
This inequality is independent of $j$ and since $\Sigma_0^j(\delta)\subset \Sigma_0^{j+1}(\delta)$ we deduce that 
\be
\label{M-49-aa}
\nu_\infty(\Sigma_0)\le N\ \delta\ \nu_\infty(\Sigma)\quad.
\ee
Since this holds for any $\delta>0$ we have proven 
\be
\label{M-50-aa}
\nu_\infty(\Sigma_0)=0\quad.
\ee
This completes the proof of lemma~\ref{lm-vanish}.\hfill $\Box$

\medskip

The next goal is to prove the following orthogonal decomposition of the limiting measure $\nu_\infty$.
\begin{Lm}
\label{lm-quantization} {\bf[Structure of the Limiting Measure]}
Under the assumptions of theorem~\ref{th-limit},  we have the existence of finitely many points $a_1\cdots a_n$ in $\Sigma$ such that the measure $\nu_\infty$ decomposes orthogonally
as follows
\be
\label{M-30}
\nu_\infty= m(x) \ {\mathcal L}^2+\sum_{i=1}^n\ \al_i\ \delta_{a_i}\quad,
\ee
where ${\mathcal L}^2$ is the Lebesgue measure on $\Sigma$ equipped with the reference metric $g_0$, $m$ is an $L^1$ function with respect to the Lebesgue measure and $\al_i$ are positive numbers bounded from below
by the universal positive number $Q_0=\lim_{\Lambda\rightarrow 0}Q_0(\La)$ given by lemma~\ref{lm-quanti}.\hfill $\Box$
\end{Lm}
\noindent{\bf Proof of lemma~\ref{lm-quantization}.}
\noindent{\bf Step 1} We prove that
\be
\label{M-39-a}
\int_{{\mathcal O}}|d\vec{\Phi}_\infty|_{g_0}^2\ dvol_{g_0}=0\quad,
\ee
Indeed, for any $\ep>0$ to any $x\in {\mathcal O}$ we assign $r_x$ such that
\be
\label{M-40}
\int_{B_{r_x}(x)}|d\vec{\Phi}_\infty|_{g_0}^2\ dvol_{g_0}\le\int_{B_{2\, r_x}(x)}|d\vec{\Phi}_\infty|_{g_0}^2\ dvol_{g_0}\le\ep\ \nu_\infty(\ov{B_{r_x}(x)})\quad.
\ee
Extracting a Besicovitch covering $(\ov{B_{r_i}(x_i)})_{i\in I})$ such that each point of $\Sigma$ is covered by at most $N$ balls from the covering. We obtain that
\be
\label{M-40-a}
\int_{\cup_{i\in I}B_{r_i}(x_i)}|d\vec{\Phi}_\infty|_{g_0}^2\ dvol_{g_0}\le \ep\sum_{i\in I}\nu_\infty(\ov{B_{r_i}(x_i)})\le\ep\ N\ \nu_\infty(\Sigma)\quad.
\ee
and since this holds for any $\ep>0$ we obtain (\ref{M-39-a}).

\medskip

\noindent{\bf Step 2} : Proof of the {\bf absolute continuity of $\nu_\infty$} with respect to the Lebesgue measure away from the oscillation set ${\mathcal O}$. Precisely we prove in this step
\be
\label{M-31}
\nu_\infty\res (\Sigma\setminus{\mathcal O})=m \ d{\mathcal L}^2\quad,
\ee
where $m\in L^1(\Sigma)$.

\medskip

Let $\ep>0$. Following (\ref{M-40-a}), we first include ${\mathcal O}$ in an open subset ${\mathcal O}^\ep$ such that
\be
\label{M-40-c}
\int_{{\mathcal O}^\ep}|d\vec{\Phi}_\infty|_{g_0}^2\ dvol_{g_0}\le \ep\quad.
\ee
Let $x\in \Sigma^\ep:=\Sigma\setminus{\mathcal O}^\ep$ then there exists $\delta_x>0$ such that
\[
\inf_{\rho>0}\frac{\ds\int_{B_{2\rho}(x)}|d\vec{\Phi}_\infty|_{g_0}^2\ dvol_{g_0} }{\nu_\infty(\ov{B_{\rho}(x)})}\ge\delta_x\quad.
\]
We denote $F_j:=\{x\in \Sigma\setminus{\mathcal O}\quad; \delta_x>2^{-j}\}$. We then have
\[
\Sigma\setminus{\mathcal O}=\bigcup_{j\in{\N}}F_j\quad.
\]
Let $G$ be a closed subset of $\Sigma^\ep:=\Sigma\setminus{\mathcal O}^\ep$ such that ${\mathcal H}^2(G)=0$. We claim that
\be
\label{M-33}
\nu_\infty(G)=0\quad.
\ee
Since $\Sigma^\ep:=\Sigma\setminus{\mathcal O}^\ep$ is closed $G$ is compact. Let  $\al>0$ to be fixed later on. Since ${\mathcal H}^2(G)=0$ and since $G$ is compact
\[
\exists\ \beta>0\quad\quad\mbox{s.t. }\quad {\mathcal H}^2(G_\beta)\le \al\quad\mbox{where }G_\beta:=\{x\in \Sigma\ ; \ \mbox{dist}(x,G)<\beta\}\quad.
\]
Indeed the closeness of $G$ implies $G:=\cap_{n\in {\N}} G_{1/n}$, $G_{1/n}$ is decreasing for the inclusion and fundamental properties of Hausdorff measures give then ${\mathcal H}^2(G)=\lim_{n\rightarrow +\infty}{\mathcal H}^2(\cap_{n\in {\N}} G_{1/n})$. Let $j\in {\N}$.
From the covering $(\ov{B_{\beta/2}(x)})_{x\in G\cap F_j}$ we extract a Vitalli covering $(\ov{B_{\beta/2}(x_i)})_{i\in I}$ in such a way that the balls $\ov{B_{\beta/6}(x_i)}$ are disjoint. Since all the balls have the same radius $\beta/2$ with centers at distances at least $\beta/3$ each point of $\Sigma$ is covered by at most $N$ balls $B_{\beta}(x_i)$ where $N$ is a universal number. since each $x_j\in F_j$
\be
\label{M-33b}
\nu_\infty(\ov{B_{\beta/2}(x)})\le 2^{j+1}\ \int_{B_{\beta}(x)}|d\vec{\Phi}_\infty|_{g_0}^2\ dvol_{g_0} \quad.
\ee
Since all the balls $B_{\beta}(x_i)$ are included in $G_{\beta}$ we have
\be
\label{M-34}
{\mathcal H}^2\lf(\bigcup_{i\in I}B_{\beta}(x_i)\rg)\le\ \al\quad.
\ee
We have moreover
\be
\label{M-35}
\begin{array}{l}
\ds\nu_\infty(G\cap F_j)\le\sum_{i\in I}\nu_\infty\lf(\bigcup\ov{B_{\beta/2}(x_i)}\rg)\le 2^{j+1}\ \sum_{i\in I}\int_{B_{\beta}(x_i)}|d\vec{\Phi}_\infty|_{g_0}^2\ dvol_{g_0}\\[5mm]
\ds\quad\quad\quad\le 2^{j+1}\ N\ \int_{\cup_{i\in I}B_{\beta}(x_i)}|d\vec{\Phi}_\infty|_{g_0}^2\ dvol_{g_0}\le 2^{j+1}\ N\ \int_{G_\beta}|d\vec{\Phi}_\infty|_{g_0}^2\ dvol_{g_0}\quad.
\end{array}
\ee
Since $|d\vec{\Phi}_\infty|_{g_0}^2\ dvol_{g_0}$ is absolutely continuous with respect to the Lebesgue measure, for any $\eta>0$ there exists $\al>0$ such that
\be
\label{M-36}
\forall\ E\ \mbox{ measurable }\quad\quad{\mathcal H}^2(E)\le\al\quad\Longrightarrow\quad\int_{E}|d\vec{\Phi}_\infty|_{g_0}^2\ dvol_{g_0}\le\eta\quad.
\ee
Hence we finally get combining (\ref{M-34}), (\ref{M-35}) and (\ref{M-36})
\be
\label{M-37}
\ds\nu_\infty(G\cap F_j)\le 2^{j+1}\ N\ \eta\quad.
\ee
For any $j\in {\N}$ the inequality (\ref{M-37}) holds for any $\eta>0$ thus $\nu_\infty(G\cap F_j)=0$ and we deduce (\ref{M-33}). Since (\ref{M-33}) holds true for any closed measurable subset of $\Sigma^\ep:=\Sigma\setminus {\mathcal O}^\ep$, then using the fundamental property of Radon measures saying that
\[
\forall \ G\ \mbox{measurable }\quad \nu_\infty(G)=\sup\{\nu_\infty(K)\ ;\ K\subset G\ ; K\ \mbox{compact}\}\quad,
\]
 we obtain that $\nu_\infty$ for any measurable subset $G$ of $\Sigma\setminus {\mathcal O}^\ep$ satisfying on $\Sigma\setminus {\mathcal O}^\ep$
is absolutely continuous with respect to the Lebesgue measure. By making $\ep$ go to zero this implies (\ref{M-31}).



\medskip

\noindent{\bf Step 3} : Detecting the {\bf ''bubbles''}. In this step we are just splitting the {\bf oscillation set} ${\mathcal O}$ into it's {\bf vanishing part} ${\mathcal O}_0:=\Sigma_0\cap {\mathcal O}$ and the {\bf bubble part} ${\mathcal B}$ where we recall that the $\Sigma_0$ is the so called {\it vanishing set} defined in definition~\ref{df-vanish} :
\[
{\mathcal B}:={\mathcal O}\setminus\lf({\mathcal O}\bigcap\Sigma_0\rg)\quad.
\] 
Recall that we have proved in lemma~\ref{lm-vanish}  $\nu_\infty(\Sigma_0)=0$ hence
\be
\label{M-50-aa-a}
\nu_\infty({\mathcal O}_0)=0\quad.
\ee


\noindent{\bf Step 4}: {\bf Finiteness of the bubble set ${\mathcal B}$}. Precisely in this step we are proving that for the constant $Q_0>0$ given by lemma~\ref{lm-quanti}
then
\be
\label{M-54}
\forall\,x\in {\mathcal B}\quad\quad \forall\, r>0\quad \nu_\infty(B_r(x))\ge Q_0\quad.
\ee
Once (\ref{M-54}) will be established we can then deduce that ${\mathcal B}$ is made of finitely many points. Let then $x\in{\mathcal B}$, then  there exists $\delta_x>0$  and $r_x>0$ that can be taken as small as one wants
such that
\be
\label{M-54-aa}
\forall\, r<r_x\quad\quad \limsup_{k\rightarrow +\infty}\frac{\ds f(\sigma_k)\int_{B_r(x)}\ dvol_{g_{\vec{\Phi}_k}}}{\ds\sigma_k^2\,\int_{B_r(x)}\  \lf[1+|{\mathbb I}_{\vec{\Phi}_k}|^{2}_{g_{\vec{\Phi}_k}}\rg]^p\  dvol_{g_{\vec{\Phi}_k}}}\ge\delta_x>0\quad.
\ee
Let $0<r_c<r_{x}$ to be fixed later, let $\vec{\Phi}_{k'}$ a sequence for which
\be
\label{M-55}
\forall\, k'\in {\N}\quad \frac{\ds f(\sigma_{k'})\int_{B_{r_c}(x)}\ dvol_{g_{\vec{\Phi}_{k'}}}}{\ds\sigma_{k'}^2\,\int_{B_{r_c}(x)}\  \lf[1+|{\mathbb I}_{\vec{\Phi}_{k'}}|^{2}_{g_{\vec{\Phi}_{k'}}}\rg]^p\  dvol_{g_{\vec{\Phi}_{k'}}}}\ge\frac{\delta_x}{2}\quad.
\ee
By assumption (\ref{VII.2}) from theorem~\ref{th-limit} we have that $f(\sigma)=o(1/\log\sigma^{-1})$ we are ''almost'' fulfilling the assumptions of lemma~\ref{lm-quanti} except that we have a surface with boundary $B_r(x)$ and not a closed surface.
So we have to choose a ''nice'' cut $r_c$ in such a way to be able to apply the arguments of lemma~\ref{lm-quanti}.

Since $x\in {\mathcal O}$, by definition, for any $\eta>0$ there exists $\rho<r_x$ such that 
\be
\label{M-56}
\eta\ \nu_\infty(B_\rho(x))\ge \int_{B_{2\rho}}|\nabla\vec{\Phi}_\infty|^2\ dx^2\quad.
\ee
Using Fubini and the mean-value theorem we can find $r\in [\rho,2\rho]$ such that
\be
\label{M-57}
\begin{array}{l}
\ds\lim_{k\rightarrow +\infty}\|\vec{\Phi}_k(x)-\vec{\Phi}_k(y)\|^2_{(L^\infty(\p B_r(x_1)))^2}=\|\vec{\Phi}_\infty(x)-\vec{\Phi}_\infty(y)\|^2_{(L^\infty(\p B_r(x_1)))^2}\\[5mm]
\ds\quad\quad\le\lf[\int_{\p B_r(x_1)}|\nabla\vec{\Phi}_\infty|\ dl\le\rg]^2\le 8\pi \ \int_{B_{2\rho}(x_1)} |\nabla\vec{\Phi}_\infty|^2\ dx^2\quad.
\end{array}
\ee
 We take this $r=r_c$ to be our ''nice cut''. We can assume 
\[
s:=\sqrt{8\pi \ \int_{B_{2\rho}(x_1)} |\nabla\vec{\Phi}_\infty|^2\ dx^2}>0\quad,
\]
the case $s=0$ could be treated in a similar way but we would have to introduce a new small parameter... Let $\vec{q}_0:=\vec{\Phi}_\infty(x_2)$ for some 
fixed arbitrary $x_2\in \p B_r(x_1)$. For $k$ large enough we have that
\be
\label{M-58}
\vec{\Phi}_k(\p B_{r_c}(x_1))\subset B^4_{2s}(\vec{q}_0)\quad.
\ee
Let $R>4$ to be fixed later. The monotonicity formula (\ref{VII.4}) and (\ref{M-56}) imply that
\be
\label{M-59}
\mu_\infty(B^4_{R\,s}(\vec{q}_0))\le C\ R^2\, s^2\le C\, R^2\, \eta \, \nu_\infty(B_\rho(x))\quad.
\ee
Hence for $\eta$ chosen in such a way that $ C\, R^2\,\eta<1/2$ we have that for $k'$ large enough (recall that $k'$ is the sequence satisfying (\ref{M-55}) for our ''nice cut'' $r_c$ which is fixed now)
\[
\int_{B_{r_c}(x)\setminus (\vec{\Phi}_{k'})^{-1}(B^4_{R\,s}(\vec{q}_0))}\ dvol_{g_{\vec{\Phi}_{k'}}}\ge 4^{-1}\int_{B_{r_c}(x)}\ dvol_{g_{\vec{\Phi}_{k'}}}\quad.
\]
Taking the same notations of the proof of lemma~\ref{lm-quanti} where $\Sigma$ is replaced by $B_{r_c}(x)$ we can then find $\vec{q}_1\in \vec{\Phi}_{k'}(B_{r_c}(x))\setminus (E_{\pi/3}\cup F_\delta\cup B^4_{R\,s}(\vec{q}_0))$.
As in the proof of lemma~\ref{lm-quanti} we shall apply the monotonicity formula centered at this point $\vec{q}_1$ but we will remove from $\vec{\Phi}_{k'}(B_{r_c}(x))$ the balls $B^4_{t\,s}(\vec{q}_0)$ for $t\in [2,4]$. 
The monotonicity formula with boundary (see for instance \cite{Ri4}) gives for all $r>0$
\be
\label{M-60}
\begin{array}{l}
\ds\frac{d}{dr}\lf[\frac{1}{r^2}\int_{B_{r_c}(x)\cap\vec{\Phi}^{-1}(B^4_r(\vec{q}_1)\setminus B^4_{t\,s}(\vec{q}_0))}\ dvol_{g_{\vec{\Phi}}}\rg]\\[5mm]
\ds\quad=\frac{d}{dr}\lf[\int_{B_{r_c}(x)\cap\vec{\Phi}^{-1}(B^4_r(\vec{q}_1)\setminus B^4_{t\,s}(\vec{q}_0))}\ \frac{|(\vec{n}\wedge\vec{\Phi})\res(\vec{\Phi}-\vec{q}_1)|^2}{|\vec{\Phi}-\vec{q}_1|^4} \  dvol_{g_{\vec{\Phi}}} \rg]\\[5mm]
\ds\quad -\frac{1}{2\,r^3}\ \int_{B_{r_c}(x)\cap\vec{\Phi}^{-1}(B^4_r(\vec{q}_1)\setminus B^4_{t\,s}(\vec{q}_0)}(\vec{\Phi}-\vec{q}_1)\cdot d^{\ast_g}d\vec{\Phi}\ dvol_{g_{\vec{\Phi}}} \\[5mm]
\ds\quad-\frac{1}{r^3}\ \int_{{\R}^4} <\vec{q}-\vec{q}_1,\vec{\nu}> \ d{\mathcal H}^1\res \lf[\vec{\Phi}(B_{r_c}(x))\cap B^4_r(\vec{q}_1)\cap\partial B^4_{t\,s}(\vec{q}_0)\rg]\\[5mm]

\ds\quad\quad\ge  -\frac{1}{2\,r^3}\ \int_{B_{r_c}(x)\cap\vec{\Phi}^{-1}(B^4_r(\vec{q}_1)\setminus B^4_{t\,s}(\vec{q}_0))}(\vec{\Phi}-\vec{q}_1)\cdot d^{\ast_g}d\vec{\Phi}\ dvol_{g_{\vec{\Phi}}} \\[5mm]
\ds\quad-\frac{1}{r^3}\ \int_{{\R}^4} <\vec{q}-\vec{q}_1,\vec{\nu}> \ d{\mathcal H}^1\res \lf[\vec{\Phi}(B_{r_c}(x))\cap B^4_r(\vec{q}_1)\cap\partial B^4_{t\,s}(\vec{q}_0)\rg]\quad,
\end{array}
\ee
where $\vec{\nu}$ is the outward unit tangent to the surface $\vec{\Phi}_k(B_{r_c}(x))\setminus B^4_{t\,s}(\vec{q}_0)$ along the boundary $$\p(\vec{\Phi}_k(B_{r_c}(x))\setminus B^4_{t\,s}(\vec{q}_0))=\vec{\Phi}_k(B_{r_c}(x))\cap \partial B^4_{t\,s}(\vec{q}_0)$$ and perpendicular to this boundary\footnote{Observe that $\vec{\Phi}_k(\p B_{r_c}(x))\subset B^4_{t\,s}(\vec{q}_0)$ so there is no contribution from $\vec{\Phi}_k(\p B_{r_c}(x))$ outside $B^4_{t\,s}(\vec{q}_0)$.}.
We consider $\chi(t)$ a smooth non negative function supported in $[1,2]$ satisfying $\int_2^4\chi(t)\ dt=1$, $\chi\le 1$ and $|\chi'|\le 1$. We multiply the inequality (\ref{M-60}) by $\chi(t)$ and we integrate between $2$ and $4$
this gives, after observing that the first term in the r.h.s. of (\ref{M-60}) is non negative \footnote{Indeed we are taking the derivative of an integral of a positive integrand over a bigger and bigger set.},
\be
\label{M-61}
\begin{array}{l}
\ds\frac{d}{dr}\lf[\frac{1}{r^2}\int_2^4\chi(t)\ dt\int_{B_{r_c}(x)\cap\vec{\Phi}^{-1}(B^4_r(\vec{q}_1)\setminus B^4_{t\,s}(\vec{q}_0))}\ dvol_{g_{\vec{\Phi}}}\rg]\\[5mm]
\ds\quad\quad\ge  -\frac{1}{2\,r^3}\ \int_2^4\chi(t)\ dt\int_{B_{r_c}(x)\cap\vec{\Phi}^{-1}(B^4_r(\vec{q}_1)\setminus B^4_{t\,s}(\vec{q}_0))}(\vec{\Phi}-\vec{q}_1)\cdot d^{\ast_g}d\vec{\Phi}\ dvol_{g_{\vec{\Phi}}} \\[5mm]
\ds\quad-\frac{1}{r^3}\ \int_2^4\chi(t)\ dt\int_{{\R}^4} <\vec{q}-\vec{q}_1,\vec{\nu}> \ d{\mathcal H}^1\res \lf[\vec{\Phi}(B_{r_c}(x))\cap B^4_r(\vec{q}_1)\cap\partial B^4_{t\,s}(\vec{q}_0)\rg]\quad.
\end{array}
\ee
By substituting $d^{\ast_g}d\vec{\Phi}$ with it's expression deduced from (\ref{III.15}), exactly as in the proof of the monotonicity formula (\ref{VII.4}) and as in the proof of lemma~\ref{lm-quanti} the new terms involving $\sigma$
coming from the boundaries $\partial B^4_{t\,s}(\vec{q}_0)$ in the first integral of the r-h-s of (\ref{M-61}) tend to zero as $k$ tends to infinity since the distance between the center $\vec{q}_1$ and this boundary
is bounded from below by $R\, s>0$ independently of $\sigma$. So it remains then to estimate the last term in (\ref{M-61}). This is done as follows
\be
\label{M-62}
\begin{array}{l}
\ds\lf|\frac{1}{r^3}\ \int_2^4\chi(t)\ dt\int_{{\R}^4} <\vec{q}-\vec{q}_1,\vec{\nu}> \ d{\mathcal H}^1\res \lf[B_{r_c}(x)\cap \vec{\Phi}^{-1}(B^4_r(\vec{q}_1))\cap \vec{\Phi}^{-1}(\partial B^4_{t\,s}(\vec{q}_0))\rg]\rg|\\[5mm]
\ds\quad\quad\le \frac{2\,|\vec{q}_1-\vec{q}_0|}{r^3}\int_2^4 dt\ {\mathcal H}^1(\partial B^4_{t\,s}(\vec{q}_0))\\[5mm]
\ds\quad\quad\le \frac{2\,|\vec{q}_1-\vec{q}_0|}{r^3}\int_{\vec{\Phi}^{-1}(B^4_{4\,s}\setminus B^4_{2\,s})}\ 
\frac{|d|\vec{\Phi}-\vec{q}_0||_{g_{\vec{\Phi}}}}{s}\ dvol_{g_{\vec{\Phi}}}\le C\ \frac{|\vec{q}_1-\vec{q}_0|}{r^3}\ s\quad,
\end{array}
\ee
where we used successively the coarea formula for the function $|\vec{\Phi}-\vec{q}_0|/s$ and the monotonicity formula (\ref{VII.4}) in the last inequality. Observe that this term appears only for
 $r>\mbox{dist}(B^4_{4\,s}(\vec{q}_0),\vec{q}_1)>|\vec{q}_0-\vec{q}_1|/2$. Hence the integral with respect to $r$ between $\sigma$ and $1/2$ gives
 \be
\label{M-63}
\begin{array}{l}
\ds\lf|\int_\sigma^{1/2}\frac{dr}{r^3}\ \int_2^4\chi(t)\ dt\int_{{\R}^4} <\vec{q}-\vec{q}_1,\vec{\nu}> \ d{\mathcal H}^1\res \lf[\vec{\Phi}(B_{r_c}(x))\cap B^4_r(\vec{q}_1)\cap\partial B^4_{t\,s}(\vec{q}_0)\rg]\rg|\\[5mm]
 \ds\quad\quad\le \frac{C}{|\vec{q}_1-\vec{q}_0|^2}\ |\vec{q}_1-\vec{q}_0|\ s\le \frac{C}{R}\quad.
\end{array}
\ee
The rest of the argument of the proof of lemma~\ref{lm-quanti} carries through and we get that
\[
\nu_\infty(B_{r_c}(x))\ge Q_0- C/R\quad.
\]
Since we can take $R$ as large as we want, we obtain (\ref{M-54}). Hence $\nu_\infty$ restricted to ${\mathcal O}$ is equal to a finite sum of Dirac masses and this last step concludes the proof of lemma~\ref{lm-quantization}.\hfill $\Box$

\medskip

We shall now prove the following lemma

\begin{Lm}
\label{lm-no-neck}{\bf[Absence of Energy in the Necks]} 
Let $\vec{\Phi}_k$ satisfying the assumptions of theorem~\ref{th-limit}. Let $1>\eta_k>0$, $1>\delta_k>0$ and $x_k\in \Sigma$ satisfying
\be
\label{M-64}
\lim_{k\rightarrow +\infty}\log\frac{\eta_k}{\delta_k}=+\infty\quad,
\ee
and such that
\be
\label{M-65}
\lim_{k\rightarrow 0}\sup_{ j\in\{ 1\cdots\log_2 (\eta_k/\delta_k)\}}\nu_k(B_{2^{j+1}\delta_k}(x_k)\setminus B_{2^{j}\delta_k}(x_k) )=0\quad.
\ee
Then
\be
\label{M-66}
\lim_{k\rightarrow 0}\nu_k(B_{\eta_k}(x_k)\setminus B_{\delta_k}(x_k) )=0\quad.
\ee
\hfill $\Box$
\end{Lm}
\noindent{\bf Proof of lemma~\ref{lm-no-neck}.}  We argue by contradiction. If (\ref{M-66}) does not hold we can then find a subsequence that we denote still $\vec{\Phi}_k$ such that
\be
\label{M-67}
\lim_{k\rightarrow 0}\nu_k(B_{\eta_k}(x_k)\setminus B_{\delta_k}(x_k) )=A>0\quad.
\ee
Let $Q_0$ be the universal constant in the lemma~\ref{lm-quanti}. We can assume without loss of generality that 
\be
\label{M-67-bb}
A< Q_0\quad.
\ee
Indeed, if this would not be the case
 we would replace $\delta_k$ by a larger number that we keep denoting $\delta_k$ and  since (\ref{M-65}) holds we necessarily have (\ref{M-64}) for this new $\delta_k$. We have for $k$ large enough
\be
\label{M-68}
\frac{\ds\sigma_k^2\,\int_\Sigma\  \lf[1+|{\mathbb I}_{\vec{\Phi}_k}|^{2}_{g_{\vec{\Phi}_k}}\rg]^p\  dvol_{g_{\vec{\Phi}_k}}}{\ds\int_{B_{\eta_k}(x_k)\setminus B_{\delta_k}(x_k)}\ dvol_{g_{\vec{\Phi}_k}}}\le \frac{2\, \nu_\infty(\Sigma)}{A}\, f(\sigma_k)\quad.
\ee
Following the approach of step 5 of the proof of lemma~\ref{lm-quantization}, we first select  2 ''good cuts'' at the two ends of the annulus. So we choose respectively $\delta_{k,c}\in [\delta_k,2\delta_k]$ and
$\eta_{k,c}\in [\eta_k/2,\eta_k]$ such that we have respectively 
\[
s_k^2:=\lf[\int_{\p B_{\delta_{k,c}}(x_k)} |\nabla \vec{\Phi}_k|\ dl\rg]^2\le\ \pi\ \nu_k(B_{2\,\delta_k}(x_k) \setminus B_{\delta_k}(x_k) )\longrightarrow 0\quad,
\]
and 
\[
t_k^2:=\lf[\int_{\p B_{\eta_{k,c}}(x_k)} |\nabla \vec{\Phi}_k|\ dl\rg]^2\le\ \pi\ \nu_k(B_{2\,\eta_k}(x_k) \setminus B_{\eta_k}(x_k) )\longrightarrow 0\quad.
\]
Let $x_{1,k}\in \p B_{\delta_{k,c}}(x_k)$ and $x_{2,k}\in\p B_{\delta_{k,c}}(x_k)$ arbitrary. We have respectively
\be
\label{M-69}
\vec{\Phi}_k(\p B_{\delta_{k,c}}(x_k))\subset B^4_{s_k}(\vec{\Phi}_k(x_{1,k}))\quad\mbox{ and }\quad\vec{\Phi}_k(\p B_{\eta_{k,c}}(x_k))\subset B^4_{t_k}(\vec{\Phi}_k(x_{2,k}))\quad.
\ee
Arguing as in the proof of the non collapsing lemma~\ref{lm-concent}, which is a corollary of the monotonicity formula, there exists $s>0$ fixed such that
\[
\max_{\vec{q}\in {\R}^4}\,\mu_\infty(B^4_s(\vec{q}))<A/4\quad.
\]
We then have for $k$ large enough 
\[
\mu_k\lf(\vec{\Phi}_k(B_{\eta_{k,c}}(x_k)\setminus B_{\delta_{k,c}}(x_k))\setminus\lf(B^4_s(\vec{\Phi}_k(x_{1,k}))\cup B^4_s(\vec{\Phi}_k(x_{2,k}))\rg)\rg)\ge A/2\quad.
\]
As in the step 5 of the proof of lemma~\ref{lm-quantization}, we adopt the notations from the proof of lemma~\ref{lm-quanti} and replacing $\Sigma$ by the annulus $B_{\eta_{k,c}}(x_k)\setminus B_{\delta_{k,c}}(x_k)$,
 we can find $\vec{q}_k$ such that
\[
\vec{q}_k\in \vec{\Phi}_k(B_{\eta_{k,c}}(x_k)\setminus B_{\delta_{k,c}}(x_k))\setminus \lf(E_{\pi/3}\cup G_\delta\cup B^4_s(\vec{\Phi}_k(x_{1,k}))\cup B^4_s(\vec{\Phi}_k(x_{2,k})\rg)\quad.
\]
We can carry over one by one the computation of the monotonicity formula centered at $\vec{q}$, controlling the boundary terms induced by the two cuts $\vec{\Phi}_k(\partial B_{\eta_{k,c}}(x_k))$ and $\vec{\Phi}_k(\partial B_{\eta_{k,c}}(x_k))$ which stay at a distance bounded from bellow with respect to $\vec{q}_k$, following the approach of the end of the  step 5 of the proof of lemma~\ref{lm-quantization}. It is here even simpler since 
the lengths of the cuts $s_k$ and $t_k$ shrink to zero in the present case. Hence we obtain
\[
A=\lim_{k\rightarrow 0}\nu_k(B_{\eta_k}(x_k)\setminus B_{\delta_k}(x_k) )\ge Q_0\quad,
\]
which contradicts (\ref{M-67-bb}). This concludes the proof of lemma~\ref{lm-no-neck}.\hfill $\Box$

\medskip

\noindent{\bf Defining the Bubble Tree} Because of the previous quantization property, together with the no-neck energy property, following a classical combinatorics argument (In the style of proposition III.1 in  \cite{BeRi} - see also \cite{MoRi}), after extracting an ad-hoc subsequence, one can construct a family of sequences of smooth conformal injections $(\psi_k^i)_{i=1\cdots L}$
from $S^i\setminus \bigcup_{j=1}^{n_i} B_\ep(a_j^i)$ (for any $\ep$ for $k$ large enough) into $(\Sigma,{g}_{\vec{\Phi}_k})$, equipped with a strongly converging constant curvature metric $h^i_k$, in such a way that
\[
\nu_k^j:=dvol_{g_{\vec{\Phi}_k\circ\psi_k^j}}\rightharpoonup\nu^j_\infty=m^j\ dvol_{h^j_\infty}\quad\quad\mbox{ as Radon measures on } S^i\setminus \bigcup_{j=1}^{n_i} B_\ep(a_j^i)\quad.
\]
for any $\ep$ and 
\[
\sum_{i=1}^L\nu_\infty^i(S^i)=\mu_\infty(S^3)\quad.
\]
In the case for instance when the conformal class of $(\Sigma,{g}_{\vec{\Phi}_k})$ is controlled, the first bubble is given by $\Sigma$ itself and the others are $S^2$.
Except for the next lemma where we are working in the junction regions between several bubbles, the so called {\it neck regions}, we shall be working on a single bubble
that we shall generically denote $\Sigma$.

\medskip

\begin{Lm}
\label{lm-L^1-conv} {\bf [Construction of an Approximating Sequence]}
 Assume the hypothesis of theorem~\ref{th-limit} are fulfilled and that we have extracted subsequences such that $\vec{\Phi}_k$ converges weakly towards $\vec{\Phi}_\infty$
 in $W^{1,2}(\Sigma)$ and $\nu_k$ converges towards $\nu_\infty$ satisfying (\ref{M-30}) where ${\mathcal B}:=\{a_1\cdots a_l\}$ the blow-up set. Let $\phi$ be a function in $C^\infty_0(B^2_1(0))$ satisfying $\int_{B^2_1(0)}\phi(x)\ dx^2=1$ and denote $\phi_t(x):=t^{-2}\phi(x/t)$. Then, modulo  extraction of a subsequence, the family  of smooth maps $\phi_r\star\vec{\Phi}_\infty$, converging \underbar{strongly} in $W^{1,2}_{loc}(\Sigma\setminus {\mathcal B})$ to $\vec{\Phi}_\infty$ as $r$ goes to zero, satisfies
\be
\label{M-84-bb-b}
\lim_{r\rightarrow 0}\limsup_{k\rightarrow+\infty}\int_{\Sigma\setminus \cup_{l=1}^nB_\ep(a_l)}|\vec{\Phi}_k-\phi_r\star\vec{\Phi}_\infty |\ dvol_{g_{{\Phi}_k}}=0\quad.
\ee
\hfill $\Box$
\end{Lm}
\noindent{\bf Proof of lemma~\ref{lm-L^1-conv}.} 
Let $\ep>0$. Let $x\in \Sigma\setminus \cup_{l=1}^nB_\ep(a_l)$ arbitrary and $r>0$ such that there exists $k_{x,r}$ such that
\be
\label{s-000}
\forall \, k\ge k_{x,r}\quad\int_{B_{4\,r}(x)}|\nabla\vec{\Phi}_k|^2\ dx^2<\ep\quad.
\ee
As before, we use Fubini and the mean value theorem to extract a slice  $r_k\in(r,2r)$ such that
\[
\begin{array}{l}
\ds\|\vec{\Phi}_\infty(x)-\vec{\Phi}_\infty(y)\|^2_{L^\infty(\p B_{r_k}(x))^2}\le C\ \int_{B_{2\rho}}|\nabla\vec{\Phi}_\infty|^2\ dx^2\le\ep\\[5mm]
\ds\quad\mbox{ and }\quad\int_{\p B_{r_k}(x)}|\nabla\vec{\Phi}_k|^2\le \frac{C}{r}\int_{B_{2r}(x)}|\nabla\vec{\Phi}_k|^2\ dx^2< \frac{C\,\ep}{r}\\[5mm]
\ds\lf\|\vec{\Phi}_\infty(x)-\vec{\Phi}_\infty^\rho(x)\rg\|_{L^\infty(\p B_{r_k}(x))}\le\ep\quad\mbox{where }\quad\vec{\Phi}_\infty^\rho(x):=\frac{1}{|B_{2r}|}\int_{B_{2r}(x)}\vec{\Phi}_\infty\ dx^2\quad.
\end{array}
\]
Because of the weak $W^{1,2}$ convergence of $\vec{\Phi}_k$ towards $\vec{\Phi}_\infty$, and because of the uniform $W^{1,2}-$bound on $\p B_{r_k}(x)$ of $\vec{\Phi}_k(r_k,\theta)$,
by {\it Rellich Kondrachov } compact embedding theorem, $\vec{\Phi}_k(r_k,\theta)-\vec{\Phi}_\infty(r_k,\theta)$ converges to zero in $L^\infty$ norm. We then choose $k_{x,r}$ such that
\[
\forall \ k\ge k_{x,r}\quad\quad\|\vec{\Phi}_k-\vec{\Phi}_\infty\|_{L^\infty(\p B_{r_k}(x))}\le \sqrt{\ep}\quad.
\]
Denote $\Sigma_k^r(x):= B_{r_k}(x)\setminus\vec{\Phi}_k^{-1}(B^4_{R\sqrt{\ep}}(\vec{\Phi}_\infty^r(x))$   and assume that
\[
\frac{\ds\sigma_k^2\int_{\Sigma^r_k(x)}(1+|{\mathbb I}_{\vec{\Phi}_k}|^2)^p\ dvol_{g_{\vec{\Phi}_k}}}{\ds\int_{\Sigma_k^r(x)}\ dvol_{g_{\vec{\Phi}_k}}}\le \frac{1}{\log\sigma_k^{-1}}\quad.
\]
Again we can then argue word by word as in the proof of lemma~\ref{lm-quanti} for the surface $\Sigma^r_k(x)$ until (\ref{VII.208}) in order to find a point $\vec{q}$ in $\vec{\Phi}_k(\Sigma^r_k(x))\setminus(E_{\pi/3}\cup G_\delta)$. Once we have this point we perform the rest of the argument of lemma~\ref{lm-quanti} but for the  surface with boundary $\vec{\Phi}_k(B_{r_k}(x_k)\setminus\vec{\Phi}_k^{-1}(B^4_{R\ep}(\vec{\Phi}_\infty^r(x)))$. The boundary is going to generate a new term in the monotonicity formula 
\[
-\frac{1}{r^3}\ \int_{{\R}^4} <\vec{q}-\vec{q}_1,\vec{\nu}> \ d{\mathcal H}^1\res \lf[\vec{\Phi}_k(B_{r_k}(x))\cap\partial B^4_{t\,\ep}(\vec{q})\rg]\quad,
\]
for $t\in [2,4]$ that we treat exactly as in (\ref{M-62}) in order to get that for $k$ large enough $\int_{B_{r_k}(x_k)}|\nabla\vec{\Phi}_k|^2\ dx^2\ge Q_0-C/R$ which is a contradiction for $R$ large enough. Hence we have
\be
\label{sau-1}
\frac{\ds\sigma_k^2\int_{\Sigma^\rho_k(x)}(1+|{\mathbb I}_{\vec{\Phi}_k}|^2)^p\ dvol_{g_{\vec{\Phi}_k}}}{\ds\int_{\Sigma_k^r(x)}\ dvol_{g_{\vec{\Phi}_k}}}> \frac{1}{\log\sigma_k^{-1}}\quad,
\ee
and then
\be
\label{sau-2}
\begin{array}{l}
\ds\int_{B_{r_k}(x)}|\vec{\Phi}_k(y)-\vec{\Phi}^\rho_\infty(x)|\ |\nabla\vec{\Phi}_k|^2(y)\ dy^2\le R\, \sqrt{\ep}\,\int_{B_{r_k}(x)}\ |\nabla\vec{\Phi}_k|^2(y)\ dy^2\\[5mm]
\ds\quad\quad\quad+C\ \log\sigma_k^{-1}\ \sigma_k^2\int_{B_{r_k}(x)}(1+|{\mathbb I}_{\vec{\Phi}_k}|^2)^p\ dvol_{g_{\vec{\Phi}_k}}\quad.
\end{array}
\ee
Let $\phi$ be a function in $C^\infty_0(B^2_1(0))$ satisfying $\int_{B^2_1(0)}\phi(x)\ dx^2=1$ and denote $\phi_t(x):=t^{-2}\phi(x/t)$ we have for all $y\in B_r(x)$
\[
\phi_{r}\star\vec{\Phi}_\infty(y)-\vec{\Phi}_\infty^r(x)=\int_{z\in B_{2 r}(y)}\phi_r(y-z)\ \vec{\Phi}_\infty(z)\ dz^2-\int_{z\in B_{2r}(y)}\phi_r(y-z)\ \vec{\Phi}_\infty^r(x)\ dz^2\quad.
\]
Hence
\[
\begin{array}{l}
\ds|\phi_{r}\star\vec{\Phi}_\infty(y)-\vec{\Phi}_\infty^r(x)|\le\frac{C}{r^4}\int_{z\in B_{2r}(y)}\int_{v\in B_{2r}(x)}\lf|\phi\lf(\frac{y-z}{r}\rg)\rg|\ |\vec{\Phi}_\infty(z)-\vec{\Phi}_\infty(v)|\ dz^2 \ dv^2\\[5mm]
\ds\quad\quad\le \frac{C}{r^4}\int_{z\in B_{4r}(x)}\int_{v\in B_{4r}(x)}|\vec{\Phi}_\infty(z)-\vec{\Phi}_\infty(v)|\ dz^2 \ dv^2\quad.
\end{array}
\]
Thus, using Poincar\'e inequality on $B_{4r}(x)$
\be
\label{sau-3}
\forall \, x\in \Sigma\setminus \cup_{l=1}^nB_\ep(a_l)\quad\quad\|\phi_{r}\star\vec{\Phi}_\infty(y)-\vec{\Phi}_\infty^r(x)\|_{L^\infty(B_{r}(x))}^2\le C\,\int_{B_{4r}(x)}|\nabla\vec{\Phi}_\infty|^2(y)\ dy^2\quad.
\ee
Let $r$ such that
\[
\sup_{x\in \Sigma\setminus \cup_{l=1}^nB_\ep(a_l)}\nu_\infty(B_{4r}(x))\le \ep/2\quad.
\]
One takes a finite covering $(B_r(x_i))_{i\in I}$ of $\Sigma\setminus \cup_{l=1}^nB_\ep(a_l)$ by balls of fixed radius $r$ such that each point is covered by at most a universal number ${\mathfrak N}$ of balls of size $2 r$. Summing (\ref{sau-2}) gives for $k$ large enough
\[
\begin{array}{l}
\ds\sum_{i\in I}\int_{B_{r}(x_i)}|\vec{\Phi}_k(y)-\vec{\Phi}^r_\infty(x_i)|\ |\nabla\vec{\Phi}_k|^2(y)\ dy^2\\[5mm]
\ds\quad\quad\le R\, {\mathfrak N}\,\sqrt{\ep}\,\int_{\Sigma\setminus \cup_{l=1}^nB_{\ep_0}(a_l)}\ |\nabla\vec{\Phi}_k|^2(y)\ dy^2+C\  {\mathfrak N}\  \log\sigma_k^{-1}\ \sigma_k^2\int_{\Sigma}(1+|{\mathbb I}_{\vec{\Phi}_k}|^2)^p\ dvol_{g_{\vec{\Phi}_k}}\quad.
\end{array}
\]
Combining this inequality with (\ref{sau-3}) gives then
\be
\label{sau-4}
\begin{array}{l}
\ds\int_{\Sigma\setminus \cup_{l=1}^nB_{2\ep_0}(a_l)}|\vec{\Phi}_k-\phi_{r}\star\vec{\Phi}_\infty|(y)\ |\nabla\vec{\Phi}_k|^2(y)\ dy^2\le C\ \sqrt{\ep}\,\int_{\Sigma\setminus \cup_{l=1}^nB_{\ep_0}(a_l)}\ |\nabla\vec{\Phi}_k|^2(y)\ dy^2\\[5mm]
\ds\quad\quad+C\  {\mathfrak N}\ \log\sigma_k^{-1}\  \sigma_k^2\int_{\Sigma\setminus \cup_{l=1}^nB_\ep(a_l)}(1+|{\mathbb I}_{\vec{\Phi}_k}|^2)^p\ dvol_{g_{\vec{\Phi}_k}}\quad.
\end{array}
\ee
This concludes the proof of the lemma.\hfill $\Box$

\begin{Lm}
\label{lm-rect} {\bf[Rectifiability of the Limit]}
Let $\vec{\Phi}_k$ satisfying the assumptions of theorem~\ref{th-limit}. Then the limiting measure $\mu_\infty$ is supported by a rectifiable 2-dimensional subset $K$ of $S^3$ given by the image of the different bubbles by the $W^{1,2}$ map $\vec{\Phi}_\infty$. Precisely there exists a uniformly bounded ${\mathcal H}^2$ measurable function $\theta$ on $K$ such that
\be
\label{sau-4-0}
\mu_\infty=\theta\ d{\mathcal H}^2\res K\quad.
\ee
Moreover if we decompose $\mu_\infty=\sum_{i=1}^L\mu_\infty^i$ where each $\mu_\infty^i$ is the limiting measure produced by one bubble we have for each bubble
\be
\label{M-70}
\mu^i_\infty(\phi)=\int_\Sigma \phi(\vec{\Phi}_\infty)\ d\nu^i_\infty=\int_\Sigma \phi(\vec{\Phi}_\infty)\ m^i(x)\ dx^2\quad,
\ee
where $\nu_\infty^i=m^i\ d{\mathcal L}^2$. \hfill $\Box$ 
\end{Lm}
\noindent{\bf Proof of lemma~\ref{lm-rect}.} We first prove (\ref{M-70}). Let $\ep>0$. Using (\ref{sau-4}) we have the existence of $r$ such that, for $k$ large enough
\be
\label{sau-5}
\begin{array}{l}
\ds\int_{\Sigma\setminus \cup_{l=1}^nB_{2\ep_0}(a_l)}|\vec{\Phi}_k-\phi_{r}\star\vec{\Phi}_\infty|(y)\ |\nabla\vec{\Phi}_k|^2(y)\ dy^2\le C\ \sqrt{\ep}\,\int_{\Sigma\setminus \cup_{l=1}^nB_{\ep_0}(a_l)}\ |\nabla\vec{\Phi}_k|^2(y)\ dy^2\\[5mm]
\ds\quad\quad+C\  {\mathfrak N}\ \log\sigma_k^{-1}\  \sigma_k^2\int_{\Sigma\setminus \cup_{l=1}^nB_\ep(a_l)}(1+|{\mathbb I}_{\vec{\Phi}_k}|^2)^p\ dvol_{g_{\vec{\Phi}_k}}\quad.
\end{array}
\ee
Let $\varphi\in C^1({\R}^4)$ we have
\be
\label{sau-6}
\begin{array}{l}
\ds\mu^1_k(\varphi)=\int_{\Sigma\setminus\cup_{l=1}^nB_{\ep_0}(a_l)} \varphi(\vec{\Phi}_k)\ dvol_{g_{\vec{\Phi}_k}}=\int_{\Sigma\setminus\cup_{l=1}^nB_{\ep_0}(a_l)} \varphi(\phi_{r}\star\vec{\Phi}_\infty)\ dvol_{g_{\vec{\Phi}_k}}\\[5mm]
\ds\quad\quad+\int_{\Sigma\setminus\cup_{l=1}^nB_{\ep_0}(a_l)} \varphi(\vec{\Phi}_k)-\varphi(\phi_{r}\star\vec{\Phi}_\infty)\ dvol_{g_{\vec{\Phi}_k}}\quad,
\end{array}
\ee
where $\mu^1_k$ is the measure issued from $\vec{\Phi}_k$ restricted to $\Sigma\setminus {\mathcal B}$. We have in one hand by the convergence of Radon measures 
\be
\label{sau-7}
\begin{array}{l}
\ds\lim_{k\rightarrow +\infty}\int_{\Sigma\setminus\cup_{l=1}^nB_{\ep_0}(a_l)} \varphi(\phi_{r}\star\vec{\Phi}_\infty)\ dvol_{g_{\vec{\Phi}_k}}=\nu_\infty(\varphi(\phi_{r}\star\vec{\Phi}_\infty))\\[5mm]
\ds\quad\quad=\int_{\Sigma\setminus\cup_{l=1}^nB_{\ep_0}(a_l)}\varphi(\phi_{r}\star\vec{\Phi}_\infty)\ m^1(x)\ dx^2\quad,
\end{array}
\ee
and in the other hand we have
\be
\label{sau-8}
\begin{array}{l}
\ds\lf|\int_{\Sigma\setminus\cup_{l=1}^nB_{\ep_0}(a_l)} \varphi(\vec{\Phi}_k)-\varphi(\phi_{r}\star\vec{\Phi}_\infty)\ dvol_{g_{\vec{\Phi}_k}}\rg|\\[5mm]
\ds\quad\le \|\nabla\varphi\|_\infty\ 
\int_{\Sigma\setminus \cup_{l=1}^nB_{2\ep_0}(a_l)}|\vec{\Phi}_k-\phi_{r}\star\vec{\Phi}_\infty|(y)\ |\nabla\vec{\Phi}_k|^2(y)\ dy^2\quad.
\end{array}
\ee
Combining (\ref{sau-5})...(\ref{sau-8}) we obtain
\be
\label{sau-9}
\limsup_{k\rightarrow +\infty}\lf|\mu_k^1(\varphi)-\nu_\infty(\varphi(\phi_{r}\star\vec{\Phi}_\infty))\rg|\le C_\varphi\ \ep\quad.
\ee
By taking $\ep$ smaller and smaller as well as $\rho$ gets smaller and smaller we obtain (\ref{M-70}). It remains to prove  (\ref{sau-4-0}). Because of the monotonicity formula $\mu^\infty$ vanishes
on any measurable set of ${\mathcal H}^2$ measure zero in $S^3$. 
Using  the quantitative Lusin type property for Sobolev maps of F.C. Liu (see \cite{Liu}) we deduce that for any $\al>0$ there exists a $C^1$ map $\vec{\Xi}^\al$ from $\Sigma$ into\footnote{The fact that we can apply Liu's result for maps into $W^{1,2}(\Sigma,S^3)$ comes
from the fact that smooth maps in $C^1(\Sigma,S^3)$ are dense in $W^{1,2}(\Sigma,S^3)$ for the $W^{1,2}-$topology.} $S^3$ and an open subset $B^\al$ of $\Sigma$ such that
\be
\label{M-42-i}
\lf\{
\begin{array}{l}
\ds{\mathcal H}^2(B^\al)\le \al\quad,\\[5mm]
\ds\vec{\Phi}_\infty=\vec{\Xi}^\al\quad\mbox{ on }\Sigma\setminus B^\al\quad\mbox{and}\quad d\vec{\Phi}_\infty=d\vec{\Xi}^\al\quad\mbox{ on }\Sigma\setminus B^\al\quad,\\[5mm]
\ds\|\vec{\Phi}_\infty-\vec{\Xi}^\al\|^2_{W^{1,2}(\Sigma)}\le\al\quad.
\end{array}
\rg.
\ee
The identity (\ref{M-70}) implies then
\[
\mu_\infty^i(\varphi)=\int_{\Sigma\setminus B^\al}\varphi(\vec{\Xi}^\al)\ d\nu_\infty^i+\int\varphi(\vec{\Phi}_\infty)\ d\nu^i_\infty\res B^\al
\]
Since $\vec{\Xi}^\al$ is $C^1$ on $\Sigma$ the measurable set $K^\al:=\vec{\Xi}^\al(\Sigma)$ is 2 rectifiable and there exists a measure $\tau^\al$ supported on $K^\al$ such that
\[
\mu_\infty^i(\varphi)=\int_{K^\al} \varphi(\vec{q})\ d\tau^\al(\vec{q})+\int\varphi(\vec{\Phi}_\infty)\ d\nu^i_\infty\res B^\al\quad.
\]
Observe that since $\nu_\infty^i$ is absolutely continuous with respect to the Lebesgue measure on $\Sigma$ we have
\[
\lim_{\al\rightarrow 0}\sup_{|E|<\al}\nu_\infty^i(E)=0\quad.
\]
Hence, by taking $K:=\cup_{n\in {\N}^\ast} K^{1/n}$, there exists a measure $\tau$ on $K$ such that $\mu_\infty:= \tau\res K$. Because of the monotonicity formula $\mu^\infty$ vanishes
on any measurable set of ${\mathcal H}^2$ measure zero in $S^3$ and hence $\tau$ is absolutely continuous with respect to $d{\mathcal H}^2\res K$ and there exist an ${\mathcal H}^2$
measurable function $\theta$ on $K$ such that (\ref{sau-4-0}) holds and this concludes the proof of lemma~\ref{lm-rect}.\hfill $\Box$

\medskip

\begin{Lm}
\label{vanish-rank<2}{\bf [Vanishing of the Limiting Measure on the Degenerating Set]}
Let ${\mathfrak L}_{\nabla\vec{\Phi}_\infty}$ be the  subset of $\Sigma\setminus{\mathcal B}$ of Lebesgue points  for $\nabla\vec{\Phi}_\infty$. We denote by ${\mathfrak L}^0_{\nabla\vec{\Phi}_\infty}$
the measurable subset of ${\mathfrak L}_{\nabla\vec{\Phi}_\infty}$ of points where the Lebesgue representative of $\nabla\vec{\Phi}_\infty$ has rank strictly less than 2. Then we have
\be
\label{M-70-aa}
\nu_\infty({\mathfrak L}^0_{\nabla\vec{\Phi}_\infty})=0\quad.
\ee
\hfill $\Box$
\end{Lm}
\noindent{\bf Proof of lemma~\ref{vanish-rank<2}.} Let $K$ be a compact subset of ${\mathfrak L}^0_{\nabla\vec{\Phi}_\infty}$ such that
\[
\nu_\infty(K)\ge 2^{-1}\nu_\infty({\mathfrak L}^0_{\nabla\vec{\Phi}_\infty})\quad.
\]
Let $\al>0$ and consider $B^\al$ and $\vec{\Xi}^\al$ satisfying (\ref{M-42-i}). We choose $\al$ small enough in such a way that
\[
\nu_\infty(K\setminus B^\al)\ge 2^{-1}\nu_\infty(K)\quad.
\]
Since $\int_{{\mathfrak L}^0_{\nabla\vec{\Phi}_\infty}}|\p_{x_1}\vec{\Phi}_\infty\times\p_{x_2}\vec{\Phi}_\infty|\ dx^2=0$ and since $\nabla\vec{\Phi}_\infty=\nabla\vec{\Xi}^\al$ on $\Sigma\setminus B^\al$ we have that ${\mathcal H}^2(\vec{\Xi}^\al(K\setminus B^\al))=0$. Observe that $\Om^\al:=\vec{\Xi}^\al(K\setminus B^\al)$ is compact in ${\R}^4$. Let $B_{\rho_i}(\vec{q}_i)$
be a finite covering of $\Om^\al$ such that $\sum_{i\in I} \rho_i^2\le  \al$. Let $\varphi^\al$ be a $C^1$ non negative function in ${\R}^4$, identically equal to one on $\Om^\al$, less than one and supported
in $\cup_{i\in I} B_{\rho_i}(\vec{q}_i)$. Because of the monotonicity formula we have
\be
\label{sau-10}
\int_{{\R}^4}\ \varphi^\al(\vec{q})\ d\mu^\infty(\vec{q})\le C\,\al\quad.
\ee
The formula (\ref{M-70}) and the fact that $\varphi^\al(\vec{\Xi}^\al)$ is identically equal to one on $K\setminus B^\al$ gives
\[
\nu_\infty(K\setminus B^\al)\le \int_{{\R}^4}\ \varphi^\al(\vec{q})\ d\mu^\infty(\vec{q})\quad,
\]
hence we obtain that $\nu_\infty({\mathfrak L}^0_{\nabla\vec{\Phi}_\infty})\le 4\, C\,\al$ for any $\al$ and
  this concludes the proof of lemma~\ref{vanish-rank<2}.\hfill $\Box$

\medskip

\medskip

\begin{Lm}
\label{lm-strong} {\bf[Convergence to an Integer Rectifiable Varifold]}
Under the assumptions of theorem~\ref{th-limit}, we have that  one we can extract a subsequence such that the integer varifold ${\mathbf v}_k$ associated to the current
$(\vec{\Phi}_k)_\ast[\Sigma]$ converges to an integer rectifiable varifold supported by a finite union of the images by  $W^{1,2}-$maps of  surfaces. More precisely
we have that on each bubble there exists a function $N^i\in L^\infty(S^i,{\N})$ such that 
\be
\label{targ-ai}
\nu_\infty^i=N^i\ |\p_{x_1}\vec{\Phi}_\infty\times\p_{x_2}\vec{\Phi}_\infty|\ dx^2\quad.
\ee  \hfill $\Box$
\end{Lm}
\noindent{\bf Proof of lemma~\ref{lm-strong}.}
Since we have proved that the necks contain no energy at the limit, it suffices to prove the convergence for $\vec{\Phi}_k$ restricted to $\Sigma\setminus\cup_{l=1}^n B_\ep(a_l)$. 
We denote by ${\mathbf v}_{\ep,k}$ the integer varifold  associated to the current $(\vec{\Phi}_k)_\ast[\Sigma\setminus\cup_{l=1}^n B_\ep(a_l)]$

\medskip

The proof of lemma~\ref{lm-strong} is a bit long and is therefore decomposed into two main parts. In the first part we establish the varifold convergence of ${\mathbf v}_{\ep,k}$ towards a limiting varifold ${\mathbf v}_{\ep,\infty}$ which is - as a Radon measure on the Grassman bundle of $TN^n$ -
absolutely continuous with respect to $(\vec{\Phi}_\infty)_\ast \delta_{T\Sigma\setminus\cup_{l=1}^n B_\ep(a_l)}$. The second step consists in proving the integrality of ${\mathbf v}_{\ep,\infty}$ 

\medskip

\noindent{\bf Step 1} : {\bf The convergence of ${\mathbf v}_{\ep,k}$ towards ${\mathbf v}_{\ep,\infty}<<(\vec{\Phi}_\infty)_\ast \delta_{T\Sigma\setminus\cup_{l=1}^n B_\ep(a_l)}$}.

\medskip

We fix $\al>0$ and we consider the map $\vec{\Xi}^\al$ and the open set $B^\al$ given by (\ref{M-42-i}). We choose a Lebesgue point $x$ for $\nabla\vec{\Phi}_\infty$ in $\Sigma\setminus B^\al$ such that
\be
\label{sau-10-a}
\lim_{r\rightarrow 0}\frac{|B_r(x)\setminus B^\al|}{|B_r(x)|}=1\quad.
\ee
We  assume that $x$ is not in the vanishing set $\Sigma_0$.We also assume that $x$ is not in the degenerating set ${\mathfrak L}^0_{\nabla\vec{\Phi}_\infty}$. These restrictions have no consequences since we have respectively $\nu_\infty<<{\mathcal L}^2$, $\nu_\infty(\Sigma_0)=0$ and $\nu_\infty({\mathfrak L}^0_{\nabla\vec{\Phi}_\infty})=0$. Such a point is a Lebesgue point for $x$ and one has
\be
\label{sau-11}
\lim_{r\rightarrow 0}\phi_r\star\vec{\Phi}_\infty(x)=\vec{\Xi}^\al(x)=\vec{\Phi}_\infty(x)\quad.
\ee
Without loss of generality, modulo the action of rotations, we assume that $\vec{\Xi}^\al(x)=\vec{\Phi}_\infty(x)=(0,0,1,0)$, that $\p_{x_1}\vec{\Xi}^\al(x)=\p_{x_1}\vec{\Phi}_\infty(x)=(a,0,0,0)$ and $\p_{x_2}\vec{\Xi}^\al(x)=\p_{x_2}\vec{\Phi}_\infty(x)=(b,c,0,0)$. We have $a\, c\ne 0$
since $\nabla\vec{\Phi}_\infty$ has rank 2. Moreover the approximate tangent plane at $\vec{\Phi}_\infty(x)$ coincides with Span$\{(1,0,0,0),(0,1,0,0)\}$. Observe that the existence of this approximate tangent plane and the
fact that $\vec{\Xi}^\al(x)$ is a regular point for $\vec{\Xi}^\al$ forces
Span$\{\p_{x_1}\vec{\Xi}^\al,\p_{x_2}\vec{\Xi}^\al\}=\{(1,0,0,0),(0,1,0,0)\}$ at any point in $(\vec{\Xi}^\al)^{-1}(\vec{\Xi}^\al(x))$.

\medskip

We recall that we adopt the notation $\vec{\Phi}=(\Phi^1,\Phi^2,\Phi^3,\Phi^4)$. We first have for the third coordinate 
\be
\label{M-92}
\begin{array}{l}
\ds\int_{B_r(x)}|\nabla\Phi_k^3|^2\ dy^2=\int_{B_r(x)}|\Phi^3_k\,\nabla\Phi_k^3|^2\ dy^2+\int_{B_r(x)}(1-|\Phi^3_k|^2)\,|\nabla\Phi_k^3|^2\ dy^2\\[5mm]
\ds \quad=\int_{B_r(x)}|\Phi^3_k\,\nabla\Phi_k^3|^2\ dy^2+\int_{B_r(x)}(|\Phi^3_\infty(x)|^2-|\Phi^3_k|^2)\,|\nabla\Phi_k^3|^2\ dy^2\quad.
\end{array}
\ee
We have $\Phi^3_k\, \nabla\Phi_k^3=-\Phi^1_k\, \nabla\Phi_k^1-\Phi^2_k\, \nabla\Phi_k^2-\Phi^4_k\, \nabla\Phi_k^4$ and since also for any $i=1\cdots 4$ we have
\be
\label{M-92-a}
|\nabla\Phi^i_k|^2\ dy_1\wedge dy_2\le\, 2\ dvol_{g_{\vec{\Phi}_k}}\quad,
\ee
and keeping in mind also $|\Phi_k^i|\le 1$, we deduce that (\ref{M-92}) gives
\[
\begin{array}{l}
\ds\int_{B_r(x)}|\nabla\Phi_k^3|^2\ dy^2\le 2\,\int_{B_r(x)} \lf[|\Phi_k^1(y)|^2+|\Phi_k^2(y)|^2+|\Phi_k^4(y)|^2\rg]\, dvol_{g_{\vec{\Phi}_k}}+4\, \int_{B_r(x)}|\Phi_k^3(y)-\Phi_\infty^3(x)|\ dvol_{g_{\vec{\Phi}_k}}\ .
\end{array}
\]
Since $\Phi_\infty^i(x)=0$ for $i\ne 3$ we have then
\be
\label{M-93}
\begin{array}{l}
\ds\int_{B_r(x)}|\nabla\Phi_k^3|^2\ dx^2\le \,10\  \int_{B_r(x)} |\vec{\Phi}_k-\vec{\Phi}_\infty(x)|\ dvol_{g_{\vec{\Phi}_k}}\quad.
\end{array}
\ee
We have
\be
\label{sau-12}
\int_{B_r(x)} |\vec{\Phi}_k-\vec{\Phi}_\infty(x)|\ dvol_{g_{\vec{\Phi}_k}}\le\int_{B_r(x)} |\vec{\Phi}_k-\phi_r\star\vec{\Phi}_\infty(x)|\ dvol_{g_{\vec{\Phi}_k}}+ |\vec{\Phi}_\infty(x)-\phi_r\star\vec{\Phi}_\infty(x)|\ \nu_k(B_r(x))\quad.
\ee
For any $\ep>0$, for $r$ small enough, using (\ref{sau-2}) and (\ref{sau-3}) we have the existence of a radius $r_k\in (\rho/2,\rho)$ such that

\be
\label{sau-13}
\begin{array}{l}
\ds\int_{B_{r_k}(x)}|\vec{\Phi}_k(y)-\phi_r\star\vec{\Phi}_\infty(y)|\ |\nabla\vec{\Phi}_k|^2(y)\ dy^2\le C\, \sqrt{\ep}\,\int_{B_{r_k}(x)}\ |\nabla\vec{\Phi}_k|^2(y)\ dy^2\\[5mm]
\ds\quad\quad\quad+C\ \log\sigma_k^{-1}\ \sigma_k^2\int_{B_{r_k}(x)}(1+|{\mathbb I}_{\vec{\Phi}_k}|^2)^p\ dvol_{g_{\vec{\Phi}_k}}\quad.
\end{array}
\ee
Since we are at a point which does not belong to the vanishing set we obtain, modulo extraction of a subsequence
\be
\label{M-84}
\lim_{r\rightarrow 0}\limsup_{k\rightarrow +\infty}\frac{\ds\int_{B_r(x)}|\vec{\Phi}_k(y)-\vec{\Phi}_\infty(x)|\ |\nabla\vec{\Phi}_k|^2(y)\ dy^2}{\ds\int_{B_r(x)}\ |\nabla\vec{\Phi}_k|^2(y)\ dy^2}=0\quad.
\ee
Combining (\ref{M-93})  with (\ref{M-84})  we obtain
\be
\label{M-94}
\lim_{r\rightarrow 0}\limsup_{k\rightarrow +\infty}\frac{\ds\int_{B_r(x)}|\nabla\Phi_k^3|^2\ dx^2}{\ds\int_{B_r(x)}\ dvol_{g_{\vec{\Phi}_k}}}=0\quad.
\ee

Since $x$ is a Lebesgue point for $\nabla\vec{\Phi}_\infty$ one has
\[
\int_{B_r(x)}|\nabla\vec{\Phi}(y)-\nabla\vec{\Xi}^\al(x)|^2=o(r^2)\quad.
\]
Then, using Fubini theorem together with the mean value theorem, for any $r>0$ and for each $k$ one can find a ``good slice'' $r_k(r)\in[2r,4r]$ such that
\be
\label{sau-11-a}
\lf\{
\begin{array}{l}
\ds\int_{0}^{2\pi}\lf|\p_\theta\lf(\vec{\Phi}_\infty(r_k(r),\theta)-r_k(r)\,\cos\theta\, \p_{x_1}\vec{\Xi}^\al(x)-r_k(r)\,\sin\theta\, \p_{x_2}\vec{\Xi}^\al(x)\rg)\rg|\ d\theta=o(r)\quad,\\[5mm]
\ds{\mathcal H}^1( \p B_{r_k(r)}(x)\cap B^\al)=o(r)\quad,\\[5mm]
\ds \int_{\p B_{r_k(r)}(x)}|\nabla\vec{\Phi}_k|^2\ dl_{\p B_{r_k(r)}}\le\frac{2}{r}\int_{B_{4\,r}(x)}|\nabla\vec{\Phi}_k|^2\ dx^2\quad.
\end{array}
\rg.
\ee
Since  
\be
\label{sau-11-b}
\|\vec{\Xi}^\al(r_k(r),\theta)-\vec{\Xi}^\al(x)-r_k(r)\,\cos\theta\ \p_{x_1}\vec{\Xi}^\al(x)-r_k(r)\,\sin\theta\ \p_{x_2}\vec{\Xi}^\al(x)\|_{L^\infty([0,2\pi])}= o(r)\quad,
\ee
 from (\ref{sau-11-a}) and (\ref{sau-11-b}) we deduce
\be
\label{sau-11-c}
\|\vec{\Phi}_\infty(r_k(r),\theta)-\vec{\Xi}^\al(x)-r_k(r)\,\cos\theta\ \p_{x_1}\vec{\Xi}^\al(x)-r_k(r)\,\sin\theta\ \p_{x_2}\vec{\Xi}^\al(x)\|_{L^\infty([0,2\pi])}= o(r)\quad.
\ee
Moreover since $\vec{\Phi}_k(r_k,\theta)-\vec{\Phi}_\infty(r_k,\theta)$ weakly in $H^{1/2}([0,2\pi])$ because of the last condition of (\ref{sau-11-a}), there exists $k_{x,r}\in{\N}$ such that
\be
\label{sau-11-d}
\forall \ k\ge k_{x,r}\quad\quad\|\vec{\Phi}_k(r_k(r),\theta)-\vec{\Xi}^\al(x)-r_k(r)\,\cos\theta\ \p_{x_1}\vec{\Xi}^\al(x)-r_k(r)\,\sin\theta\ \p_{x_2}\vec{\Xi}^\al(x)\|_{L^\infty([0,2\pi])}= o(r)\quad.
\ee

\medskip

Because of (\ref{sau-11-d}), there exists $k_{x,r}$ such that
\[
\forall\ k\ge k_{x,r}\quad \vec{\Phi}_k(\p B_{r_k(r)}(x))\subset B^4_{3\,|\nabla\vec{\Xi}^\al(x)|\, r}(\vec{\Phi}_\infty(x))\setminus B^4_{\gamma r}(\vec{\Phi}_\infty(x))\quad,
\]
where $\gamma:=\inf\{|\p_{x_1}\vec{\Xi}^\al(x),\p_{x_2}\vec{\Xi}^\al(x)\}$. For any $\tau>2\, |\nabla\vec{\Xi}^\al(x)|\ r$ we denote by $\om_k(\tau)$ the component of $\vec{\Phi}_k^{-1}(B^4_\tau(\vec{\Phi}_\infty(x)))$ containing $\p B_{r_k(r)}(x)$. Let
\[
\Om_k(\tau):=\om_k(\tau)\cup B_{r_k(r)}(x)\quad.
\]
Replacing $r$ by $\gamma^{-1} r/4\,|\nabla\vec{\Xi}^\al(x)|$ the corresponding ``good cut'' at $r_k(\gamma^{-1} r/4\,|\nabla\vec{\Xi}^\al(x)|)$ is sent by $\vec{\Phi}_k$ outside $B^4_{4\,|\nabla\vec{\Xi}^\al(x)|\, r}(\vec{\Phi}_\infty(x))$ hence, since $\p\Om_k(\tau)\subset \vec{\Phi}_k^{-1}(\p B^4_{\tau}(\vec{\Phi}_\infty))$
\be
\label{sau-11-e}
\forall \tau\in \lf[2\, |\nabla\vec{\Xi}^\al(x)|\ r,4\, |\nabla\vec{\Xi}^\al(x)|\ r\rg]\quad\quad \Om_k(\tau)\subset B_{\gamma^{-1} r/2\,|\nabla\vec{\Xi}^\al(x)|}(x)\quad.
\ee
We denote
\[
\Sigma_{k,r}:=\Om_k\lf(4\,|\nabla\vec{\Xi}^\al(x)|\ r\rg)\quad.
\]
Let $\chi^\al_{r,x}$ be a smooth non negative function on ${\R}^4$ supported in the ball $B^4_{4\, |\nabla\vec{\Xi}^\al(x)|\,r}(\vec{\Phi}_\infty(x))$, identically equal to one on $B^4_{3\, |\nabla\vec{\Xi}^\al(x)|\,r}(\vec{\Phi}_\infty(x))$
and such that $\|d^l\chi^\al_{r,x}\|_{L^{\infty}({\R}^4)}\le\ r^{-l}\ |\nabla\vec{\Xi}^\al(x)|^{-l}_\infty$ for $l=0,1,2$. We have in particular for $j=1\cdots 4$
\be
\label{sau-12}
\int_{B_{\gamma^{-1} r/2\,|\nabla\vec{\Xi}^\al(x)|}(x)}\ |\nabla\Phi_k^j|^2\ dx^2\ge\int_{\Sigma_{k,r}}\chi^\al_{r,x}(\vec{\Phi}_k)\ |\nabla\Phi_k^j|^2\ dx^2\ge \int_{B_r(x)}\ |\nabla\Phi_k^j|^2\ dx^2\quad.
\ee
Multiplying the 4th coordinate of equation (\ref{III.15}) by $\chi^\al_{r,x}(\vec{\Phi}_k)\ \Phi_k^4$ and integrating over $\Sigma$ gives
, arguing exactly as in the proof of lemma~\ref{lm-VII.1},
\be
\label{M-95}
\begin{array}{l}
\ds\int_{\Sigma_{k,r}} \chi^\al_{r,x}(\vec{\Phi}_k) |\nabla\Phi_k^4|^2\ dx^2=\int_{\Sigma_{k,r}} \chi^\al_{r,x}(\vec{\Phi}_k) |\Phi_k^4|^2\ |\nabla\Phi_k^4|^2\ dx^2\\[5mm]
\ds\quad\quad-\int_{\Sigma_{k,r}}\Phi_k^4\ \nabla (\chi^\al_{r,x}(\vec{\Phi}_k))\cdot \nabla\Phi_k^4+o_k(1)\quad.
\end{array}
\ee
We shall now define a radius $s_r=\delta(r)\, r$ where $\delta(r)=o_r(1)$ in the following way. Using Poincar\'e inequality as for proving (\ref{sau-3}) we have
\be
\label{sau-13}
\|\phi_{s_r}\star\vec{\Phi}_\infty-\vec{\Phi}_\infty^r\|_{L^\infty(\Sigma_{r,k})}^2\le \frac{C_x}{\delta_r^2}\int_{B_{\gamma^{-1} r/2\,|\nabla\vec{\Xi}^\al(x)|}(\vec{\Phi}_\infty)}|\nabla\vec{\Phi}_\infty|^2\ dx^2\quad,
\ee
where $C_x$ does not depend on $r$ but on $x$ only. Using the fact that, since $\vec{\Xi}^\al$ is $C^1$,
\be
\label{sau-13-a}
\begin{array}{l}
\ds r^{-2}\int_{B_{\gamma^{-1} r/2\,|\nabla\vec{\Xi}^\al(x)|}(x}|\nabla\vec{\Phi}_\infty-\nabla\vec{\Xi}^\al|^2\ dy^2=\ep(r)\\[5mm]
\ds\quad\mbox{ and }\quad|\nabla{\Xi}^{\al,4}|(x)=0\Rightarrow\|\nabla\Xi^{\al,4}\|_{L^\infty(B_r(x))}=o_r(1)\quad,
\end{array}
\ee
where $\ep(r)=o_r(1)$ by choosing $\delta^2(r):=\max\{{\|\nabla\Xi^{\al,4}\|_{L^\infty(B_r(x))}},\ep(r)^{1/2}\}$ we deduce from (\ref{sau-13})
\be
\label{sau-14}
\lf\|\phi_{s_r}\star{\Phi}^4_\infty-\frac{1}{B_r(x)}\int_{B_r(x)}\Phi^4_\infty\rg\|_{L^\infty(\Sigma_{r,k})}^2\le\, C_x\ \lf[\sqrt{\ep(r)}+ \|\nabla\Xi^{\al,4}\|_{L^\infty(B_r(x))}   \rg]\ r^2=o(r^2)\quad.
\ee
On $B_{r_k}(x)$ we decompose $\vec{\Phi}_\infty-\vec{\Xi}^\al=v+\psi$ such that $\Delta v=0$ in $B_{r_k}(x)$ and $\psi=0$ on $\p B_{r_k}(x)$. Because of (\ref{sau-11-c}) one has, using respectively the {\it maximum principle} and the {\it Dirichlet Principle},
\be
\label{sau-15}
\|v\|_{L^\infty(B_{r_k}(x))}=o(r)\quad\mbox{ and }\quad\int_{B_{r_k}(x)}|\nabla \psi|^2\le \int_{B_{r_k}(x)}|\nabla\vec{\Phi}_\infty-\nabla\vec{\Xi}^\al|^2\ dy^2=\ep(r)\ r^2\quad.
\ee
Sobolev-Poincar\'e inequality gives
\[
\frac{1}{|B_{r_k}(x)|}\int_{B_{r_k}(x)}|\psi|^2\le C\, \int_{B_{r_k}(x)}|\nabla\vec{\Phi}_\infty-\nabla\vec{\Xi}^\al|^2\ dx^2\quad.
\]
Combining this last fact with (\ref{sau-15}) gives
\[
\lf|\frac{1}{|B_{r_k}(x)|}\int_{B_{r_k}(x)} [\vec{\Phi}_\infty(y)-\vec{\Xi}^\al(y)]\ dy^2\rg|^2=o(r^2)\quad.
\]
This implies $\frac{1}{|B_{r_k}(x)|}\int_{B_{r_k}(x)} \vec{\Phi}^4_\infty(y)=o(r)$. Observe that similarly to the proof of (\ref{sau-3}) by the mean again of Poincar\'e inequality one has
\[
\lf|\frac{1}{|B_{r_k}(x)|}\int_{B_{r_k}(x)} {\Phi}^4_\infty(y)\ dy^2-  \frac{1}{|B_{r}(x)|}\int_{B_{r}(x)} {\Phi}^4_\infty(y)\ dy^2 \rg|^2\le C\ \int_{B_{2r}(x)}|\nabla\Phi_\infty^4|^2\ dy^2=o(r^2)\quad.
\]
Combining these two last estimates  with (\ref{sau-14}) we finally obtain
\be
\label{sau-16}
\lf\|\phi_{s_r}\star{\Phi}^4_\infty\rg\|_{L^\infty(\Sigma_{r,k})}=o(r)\quad.
\ee
We shall denote simply $\vec{\Phi}_{s_r}=\phi_{s_r}\star\vec{\Phi}_\infty$. Arguing now exactly as in the proof of lemma~\ref{lm-VII.1}, we have
\be
\label{M-95}
\begin{array}{l}
\ds\int_{\Sigma_{k,r}} \chi^\al_{r,x}(\vec{\Phi}_k) |\nabla\Phi_k^4|^2\ dx^2=\int_{\hat{\Sigma}_\ep} \chi^\al_{r,x}(\vec{\Phi}_k) |\Phi_k^4|^2\ |\nabla\Phi_k^4|^2\ dx^2\\[5mm]
\ds\quad\quad-\int_{\Sigma_{k,r}}\Phi_k^4\ \nabla (\chi^\al_{r,x}(\vec{\Phi}_k))\cdot \nabla\Phi_k^4+o_k(1)\quad.
\end{array}
\ee
Observe that from (\ref{sau-16}) one has $|\Phi^4_{s_r}|=o(r)$ hence $|\Phi_{s_r}^4\ \p_{z_j}\chi^\al_{r,x}(\vec{\Phi}_{s_r}))|\le o(r)\ r^{-1}=o(1)$. Thus we have
\be
\label{sau-17}
\begin{array}{l}
\ds\int_{\Sigma_{k,r}} [\chi^\al_{r,x}(\vec{\Phi}_k)-o_r(1)] |\nabla\Phi_k^4|^2\ dx^2=\int_{\Sigma_{k,r}} \chi^\al_{r,x}(\vec{\Phi}_k) [|\Phi_k^4|^2-|\Phi_{s_r}|^2]\ |\nabla\Phi_k^4|^2\ dx^2\\[5mm]
\ds\quad\quad-\sum_{j=1}^4\int_{\Sigma_{k,r}}[\Phi_k^4\  (\p_{z_j}\chi^\al_{r,x}(\vec{\Phi}_k))-\Phi_{s_r}^4\ \p_{z_j}\chi^\al_{r,x}(\vec{\Phi}_{s_r}))]\ \nabla\Phi_k^j\cdot \nabla\Phi_k^4+o_k(1)\quad.
\end{array}
\ee
Because of the first line in (\ref{sau-13-a}) one has
\[
\sup_{y\in \Sigma_{k,r}}\int_{B_{s_r}(y)}|\nabla\vec{\Phi}_\infty|^2(z)\ dz^2\le \ep(r)\ r^2+ C_x\ s_r^2\le C_x\ s_r^2\quad.
\]
Replacing $r$ by $s_r$ and $\ep$ by $s_r^2$ and $\Sigma\setminus\cup_{l=1}^n B_\ep(a_l)$ by $\Sigma_{k,r}$, one can transpose word by word the arguments from equation (\ref{s-000})
until equation (\ref{sau-4}) in order to obtain
\be
\label{sau-18}
\begin{array}{l}
\ds\int_{\Sigma_{k,r}}|\vec{\Phi}_k-\phi_{s_r}\star\vec{\Phi}_\infty|(y)\ |\nabla\vec{\Phi}_k|^2(y)\ dy^2\le C\ s_r\,\int_{\Sigma_{k,r}}\ |\nabla\vec{\Phi}_k|^2(y)\ dy^2\\[5mm]
\ds\quad\quad+C\  {\mathfrak N}\ \log\sigma_k^{-1}\  \sigma_k^2\int_{\Sigma_{k,r}}(1+|{\mathbb I}_{\vec{\Phi}_k}|^2)^p\ dvol_{g_{\vec{\Phi}_k}}\quad.
\end{array}
\ee
Combining (\ref{sau-17}) with (\ref{sau-18}) gives then
\be
\label{sau-17}
\begin{array}{l}
\ds\int_{\Sigma_{k,r}} [\chi^\al_{r,x}(\vec{\Phi}_k)-o_r(1)] |\nabla\Phi_k^4|^2\ dx^2\le C\ s_r\,\int_{\Sigma_{k,r}}\ |\nabla\vec{\Phi}_k|^2(y)\ dy^2\\[5mm]
\ds\quad\quad+C\  {\mathfrak N}\ \log\sigma_k^{-1}\  \sigma_k^2\int_{\Sigma_{k,r}}(1+|{\mathbb I}_{\vec{\Phi}_k}|^2)^p\ dvol_{g_{\vec{\Phi}_k}}\\[5mm]
\ds\quad\quad+C\ \int_{\Sigma_{k,r}}|\vec{\Phi}^4_k-\phi_{s_r}\star\vec{\Phi}^4_\infty|(y)\ |\p_z\chi^\al_{r,x}(\vec{\Phi}_k)\|\nabla\vec{\Phi}_k|^2(y)\ dy^2\\[5mm]
\ds\quad\quad+C\ \int_{\Sigma_{k,r}}|\vec{\Phi}^4_{s_r}|(y)\ |\p_z\chi^\al_{r,x}(\vec{\Phi}_k)- \p_z\chi^\al_{r,x}(\vec{\Phi}_{s_r})|\      \|\nabla\vec{\Phi}_k|^2(y)\ dy^2\quad.
\end{array}
\ee
Using the fact that $|\p_z\chi^\al_{r,x}|\le C\ r^{-1}$, that $|\p_z\chi^\al_{r,x}|\le C\ r^{-2}$ together with (\ref{sau-16}) and (\ref{sau-18}) again we finally obtain
\be
\label{sau-18}
\limsup_{k\rightarrow 0}\frac{\ds\int_{\Sigma_{k,r}} [\chi^\al_{r,x}(\vec{\Phi}_k)-o_r(1)] |\nabla\Phi_k^4|^2\ dx^2}{\ds\int_{\Sigma_{k,r}}\ |\nabla\vec{\Phi}_k|^2(y)\ dy^2}\le C\ [r^{-1}\, s_r+r^{-2}\ s_r^2]\quad.
\ee
Combining this fact with (\ref{sau-12}) and the fact that $s_r\,r^{-1}=o(1)$ we finally obtain
\be
\label{sau-19}
\limsup_{k\rightarrow 0}\frac{\ds\int_{B_r(x)} |\nabla\Phi_k^4|^2\ dx^2}{\ds\int_{B_r(x)}\ |\nabla\vec{\Phi}_k|^2(y)\ dy^2}=o_r(1)\quad.
\ee
Combining (\ref{M-94}) and (\ref{sau-19}) we have then
\be
\label{M-108}
\lim_{r\rightarrow 0}\lim_{k\rightarrow +\infty}\frac{\ds\int_{B_r(x)}[|\nabla\vec{\Phi}_k^1|^2+|\nabla\vec{\Phi}_k^2|^2]\ dx^2}{\ds\int_{B_r(x)}|\nabla\vec{\Phi}_k|^2\ dx^2}=1
\ee
as well as
\be
\label{M-108-b}
\lim_{\rho\rightarrow 0}\limsup_{k\rightarrow +\infty}\frac{\ds\int_{\vec{\Phi}_k^{-1}(B^4_\rho(\vec{\Phi}_\infty(x)))}[|\nabla\vec{\Phi}_k^1|^2+|\nabla\vec{\Phi}_k^2|^2]\ dx^2}{\ds\int_{\vec{\Phi}_k^{-1}(B^4_\rho(\vec{\Phi}_\infty(x)))}|\nabla\vec{\Phi}_k|^2\ dx^2}=1\quad.
\ee
Since $\vec{\Phi}_k$ is conformal we have then
\be
\label{M-109}
\lim_{r\rightarrow 0}\lim_{k\rightarrow +\infty}\frac{\ds\int_{B_r(x)} 2\ |\p_{x_1}\vec{\zeta}_k\wedge\p_{x_2}\vec{\zeta}_k|\ dx^2}{\ds\int_{B_r(x)}|\nabla\vec{\Phi}_k|^2\ dx^2}=1\quad,
\ee
where $\vec{\zeta}_k:=(\Phi_k^1,\Phi_k^2)$ and, combining (\ref{M-108}) with (\ref{M-109})
\be
\label{M-110}
\lim_{r\rightarrow 0}\lim_{k\rightarrow +\infty}\frac{\ds\int_{B_r(x)} 2\ |\p_{x_1}\vec{\zeta}_k\wedge\p_{x_2}\vec{\zeta}_k|\ dx^2}{\ds\int_{B_r(x)}|\nabla\vec{\zeta}_k|^2\ dx^2}=1\quad.
\ee
One difficulty at this stage is that we can not remove  the absolute values inside the upper integral of (\ref{M-110}).  If we would be able to do so, we would be proving the strong convergence
for $\nabla\vec{\Phi}_k$ towards $\nabla\vec{\Phi}_\infty$ and the lemma would be proven\footnote{Unfortunately we still don't know whether we can exchange the integration and the absolute values in (\ref{M-110}) at this stage of our study of the viscosity method.}. The rest of the argument consists in proving that  the limiting {\it un-oriented } varifold associated to the current $(\vec{\Phi}_k)_\ast[B_r(x)]$ is going to be equal, asymptotically as $r$ goes to zero, to an integer times $\vec{\Xi}^\al_\ast T_x\Sigma$. We formulate that differently. Denote by $\ti{G}_2(S^3)$ to be the Grassmanian of \underbar{oriented}
2 dimensional planes of the tangent bundle to $S^3$, $TS^3$. The image by $\vec{\Phi}_k$ of  $\Sigma^\al_\ep$,  induces
an \underbar{oriented integer rectifiable varifold} (see \cite{Hu})
$
\ti{{\mathbf v}}^\al_{\ep,k}\quad,
$
 where the choice of orientation of the tangent plane is taken to be the one induced by the push forward by the immersion $\vec{\Phi}_k$ of the one fixed on $\Sigma$. The sequence of oriented varifolds $\ti{{\mathbf v}}_k$  converges to a limiting 
 oriented varifold $\ti{\mathbf v}_\infty$ which is a limiting measure on the oriented 2-Grassmanian $\ti{G}_2(S^3)$. Denote by $T^+\Sigma$ the tangent bundle to $\Sigma$ with the positive orientation
 and $T^-\Sigma$ the same tangent bundle but with the opposite orientation. We see $\vec{\Xi}^\al_\ast(T^+\hat{\Sigma}^\al_\ep\cup T^-\hat{\Sigma}^\al_\ep)$ as a measurable subset  of $\ti{G}_2(S^3)$. With these notations, the identity (\ref{M-108}) is in fact equivalent to
 \be
 \label{M-111}
 \ti{\mathbf v}^\al_{\ep,\infty}\lf(\ti{G}_2(S^3)\setminus \vec{\Xi}^\al_\ast(T^+\hat{\Sigma}^\al_\ep\cup T^-\hat{\Sigma}^\al_\ep)\rg)=0\quad.
 \ee
The goal is now to prove
 \be
 \label{M-112}
 {\mathbf v}^\al_{\ep,\infty}= N_x\, \delta_{\vec{\Xi}^\al_\ast(T_x\hat{\Sigma}^\al_\ep)}\quad\mbox{ where } N_x\in {\N}^\ast\quad,
 \ee
where ${\mathbf v}^\al_{\ep,\infty}$ is the {\it un-oriented} varifold associated to $ \ti{\mathbf v}^\al_{\ep,\infty}$ and $\delta_{\vec{\Xi}^\al_\ast(T\hat{\Sigma}^\al_\ep)}$ is the Dirac
mass  at the {\it un-oriented} tangent plane $\vec{\Xi}^\al_\ast(T_x\hat{\Sigma}^\al_\ep)$.

\medskip

\noindent{\bf Step 2 : The integrality of ${\mathbf v}_{\ep,\infty}$ : The proof of (\ref{M-112}).} 

\medskip

To simplify the presentation, in order not to have to localize in the domain that would make the notations heavier, we shall assume that
\be
\label{m-112-b}
(\vec{\Xi}^\al)^{-1}\lf( \vec{\Xi}^\al(x)\rg)=\{x\}\quad.
\ee
For $i=1\cdots 4$ we denote by $\nabla^{\Sigma_k}y^i$ the vector-field tangent to $\Phi_k(\Sigma)$ given by the projection of the $i-$th canonical vector of ${\R}^4$ onto $(\vec{\Phi}_k)_\ast T\Sigma$. We also denote $\ast_k\nabla^{\Sigma_k}y^i$ the rotation by $\pi/2$ of this vector in the tangent plane to $\Phi_k(\Sigma)$, taking into account the orientation given by the push-forward by $\vec{\Phi}_k$
of the one we fixed on $\Sigma$. Denote by $(\vec{\ep}_i)_{i=1\cdots 4}$ the canonical basis of ${\R}^4$. The identity (\ref{M-108-b}) implies that
\be
\label{M-113}
\begin{array}{l}
\ds\limsup_{k\rightarrow +\infty}\int_{\vec{\Phi}_k^{-1}(B^4_\rho(\vec{\Phi}_\infty(x)))}
\mbox{dist}\lf(\frac{\p_{x_1}\vec{\Phi}_k\wedge\p_{x_2}\vec{\Phi}_k}{|\p_{x_1}\vec{\Phi}_k\wedge\p_{x_2}\vec{\Phi}_k|},\pm\ \vec{\ep}_1\wedge\vec{\ep}_2\rg)\ |\nabla\vec{\Phi}_k|^2\ dx^2=o(\rho^2)\quad.
\end{array}
\ee
recall $\mu_\infty(B^4_\rho(\vec{\Phi}_\infty(x)))\simeq \rho^2$.
This also implies
\be
\label{M-114}
\forall \ i=1,2\quad\limsup_{k\rightarrow+\infty}\int_{B^4_\rho(\vec{\Phi}_\infty(x))}|\nabla^{\Sigma_k}y_i-\vec{\ep}_i| \ d{\mathcal H}^2\res \vec{\Phi}_k(\Sigma)=o(\rho^2)
\quad.
\ee
For $(\p_{x_1}\vec{\Phi}_k\wedge\p_{x_2}\vec{\Phi}_k)\cdot(\vec{\ep}_1\wedge\vec{\ep}_2)\ne 0$ we denote $J_k=\mbox{sign}\lf((\p_{x_1}\vec{\Phi}_k\wedge\p_{x_2}\vec{\Phi}_k)\cdot(\vec{\ep}_1\wedge\vec{\ep}_2)\rg)$ otherwize we simply take $J_k=0$. Identity (\ref{M-113}) and (\ref{M-114}) imply
\be
\label{M-115}
\limsup_{k\rightarrow+\infty}\int_{B^4_\rho(\vec{\Phi}_\infty(x))}[|\ast_k\nabla^{\Sigma_k}y_1- J_k\,{\mathbf \ep}_2|+ |\ast_k\nabla^{\Sigma_k}y_2+J_k\,{\mathbf \ep}_1|]\ d{\mathcal H}^2\res \vec{\Phi}_k(\Sigma)
=o(\rho^2)\quad.
\ee
Let $\vec{T}^\rho_k$  be the following vector-valued one dimensional currents
\[
\forall\ \al\in {\Omega}^1({\R}^4)\quad\quad\lf<\vec{T}^\rho_k,\al\rg>:=\int_{B^4_\rho(\vec{\Phi}_\infty(x))\cap \vec{\Phi}_k(\Sigma)}\al\wedge \ast_k d\vec{y}=\int_{\vec{\Phi}_k^{-1}(B^4_\rho(\vec{\Phi}_\infty(x)))}\vec{\Phi}_k^\ast\al\wedge\ast \,d\vec{\Phi}_k\quad.
\]
Let $\varphi$ be a smooth function in $C^\infty_0(B^4_1(0))$ such that $\int_{{\R}^4}\varphi(y)\ dy^4=1$. Denote $\varphi_{\sigma_k}:=\sigma_k^{-4/p}\varphi(\cdot/\sigma^{1/p}_k)$.
 We recall the definition of the {\it $\sigma_k-$smoothing} $\varphi_{\sigma_k}\star\vec{T}^\rho_k$ of the current $\vec{T}^\rho_k$ (see \cite{Fe} 4.1.2)
\[
\forall\ \al\in {\Omega}^1({\R}^4)\quad\quad\lf<\varphi_{\sigma_k}\star\vec{T}^\rho_k,\al\rg>:=\int_{B^4_\rho(\vec{\Phi}_\infty(x))\cap \vec{\Phi}_k(\Sigma)}(\varphi_{\sigma_k}\star\al)\wedge \ast_k d\vec{y}\quad,
\] 
where $\al_{\sigma_k}:=\varphi_{\sigma_k}\star\al$  denotes the following convolution operation
\[
\al_{\sigma_k}=\varphi_{\sigma_k}\star\al:=\int_{{\R}^4}\varphi_{\sigma_k}(-z)\ \tau_{z}^\ast\al\ dz^4
\]
where $\tau_z(y)=y+z$. We shall use the following lemma
\begin{Lm}
\label{lm-cv-app-sigma}
{\bf[Convergence of the $\sigma_k-$Approximation of $\vec{T}^\rho_k$.]}
Under the previous notations we have
\be
\label{M-116}
\limsup_{k\rightarrow +\infty}{\sup_{\mbox{supp}(\phi)\subset B^4_\rho(\vec{\Phi}_\infty(x))\ ;\ \|d\phi\|_\infty\le 1}\lf<\vec{T}^\rho_k-\varphi_{\sigma_k}\star\vec{T}^\rho_k,d\phi\rg>}=0\quad.
\ee
\hfill $\Box$
\end{Lm}
\noindent{\bf Proof of lemma~\ref{lm-cv-app-sigma}.} Let $\phi$ be a lipschitz function supported in $B^4_\rho(\vec{\Phi}_\infty(x))$ with $\|d\phi\|_\infty\le 1$. We have
\[
\begin{array}{l}
\ds\lf<\vec{T}^\rho_k-\varphi_{\sigma_k}\star\vec{T}^\rho_k,d\phi\rg>=\int_{{\R}^4} dz\ \varphi_{\sigma_k}(-z)\int_{B^4_\rho(\vec{\Phi}_\infty(x))\cap \vec{\Phi}_k(\Sigma)} (d\phi-\tau_z^\ast d\phi)\wedge \ast_kd\vec{y}\\[5mm]
\ds\quad\quad\quad\quad=-\int_{{\R}^4} dz\ \varphi_{\sigma_k}(-z)\int_{B^4_\rho(\vec{\Phi}_\infty(x))\cap \vec{\Phi}_k(\Sigma)} (\phi(y)-\phi(y+z))\wedge d\ast_kd\vec{y}\quad.
\end{array}
\]
Using the fact that $\|d\phi\|_\infty\le 1$ and that $\varphi_{\sigma_k}$ is supported in $B^4_{\sigma^{1/p}_k}(0)$, we have
\[
\begin{array}{l}
\ds\lf|\lf<\vec{T}^\rho_k-\varphi_{\sigma_k}\star\vec{T}^\rho_k,d\phi\rg>\rg|\le \sigma^{1/p}_k\ \int_{\Sigma}[|\vec{H}_k|+1]\ dvol_{g_{\vec{\Phi}_k}}\\[5mm]
\ds\quad\quad\quad\quad\le \lf[\sigma^2_k\ \int_{\Sigma}[|\vec{H}_k|^{2p}+1]\ dvol_{g_{\vec{\Phi}_k}}\rg]^{1/2p}\ \mbox{Area}(\vec{\Phi}_k(\Sigma))^{1-1/2p}=o(1)\quad.
\end{array}
\]
This concludes the proof of lemma~\ref{lm-cv-app-sigma}.\hfill $\Box$

\begin{Lm}
\label{lm-vanishing-boundary}
{\bf[Asymptotic Vanishing of the Boundary of $\vec{T}^\rho_k$ in $B^4_\rho(\vec{\Phi}_\infty(x))$]}
Under the previous notations we have 
\be
\label{M-118}
\limsup_{k\rightarrow +\infty}\sup_{\mbox{supp}(\phi)\subset B^4_\rho(\vec{\Phi}_\infty(x))\ ;\ \|d\phi\|_\infty\le 1}\lf<\vec{T}^\rho_k,d\phi\rg>=o(\rho^2)\quad,
\ee
and for the two first directions $i=1,2$ we have
\be
\label{M-118-b}
\limsup_{k\rightarrow +\infty}\sup_{\mbox{supp}(\phi)\subset B^4_\rho(\vec{\Phi}_\infty(x))\ ;\ \|d\phi\|_\infty\le 1}\vec{\ep}_i\cdot\lf<\vec{T}^\rho_k,d\phi\rg>=O(\rho^4)\quad.
\ee
\hfill $\Box$
\end{Lm}
\noindent{\bf Proof of lemma~\ref{lm-vanishing-boundary}.} Because of (\ref{M-108-b}) it suffices to prove (\ref{M-118-b}). Because of the previous lemma it suffices to prove (\ref{M-118}) where $\vec{\ep}_i\cdot\vec{T}^\rho_k$ for $i=1,2$ is replaced by $\vec{\ep}_i\cdot\varphi_{\sigma_k}\star\vec{T}^\rho_k$. We assume $\phi(\vec{\Phi}_\infty(x))=0$ in such a way that $\|\phi\|_\infty\le \rho$. We have
\be
\label{M-119-0}
\lf<\varphi_{\sigma_k}\star\vec{T}^\rho_k,d\phi\rg>=\int_{B^4_\rho(\vec{\Phi}_\infty(x))\cap \vec{\Phi}_k(\Sigma)}\ d\lf(\varphi_{\sigma_k}\star\phi\rg)\wedge \ast_k d\vec{y}\quad.
\ee
Integrating by parts and using (\ref{III.15}) we have, omitting to write explicitly the subscript $k$,
\be
\label{M-119}
\begin{array}{l}
\ds\lf<\varphi_{\sigma}\star\vec{T}^\rho,d\phi\rg>=\int_{\vec{\Phi}_k^{-1}(B^4_\rho(\vec{\Phi}_\infty(x)))} \nabla\lf(\varphi_{\sigma}\star\phi(\vec{\Phi})\rg)\ \cdot\sigma^2\ f^p\ \nabla\vec{\Phi}\ dx^2\\[5mm]
\ds -\, 2\  p\ \sigma^2 \int_{\vec{\Phi}_k^{-1}(B^4_\rho(\vec{\Phi}_\infty(x)))} e^{-2\la}\ \nabla\lf(\varphi_{\sigma}\star\phi(\vec{\Phi})\rg)
\cdot \lf[\ov{\nabla}\lf[ f^{p-1}\, {\mathbb I}^0_{11}\rg] +  (\ov{\nabla})^\perp\lf[ f^{p-1}\, {\mathbb I}^0_{12}\rg]\ \rg]\vec{n}\ dx^2\\[5mm]
\ds-\,2\ p\ \sigma^2\ \int_{\vec{\Phi}_k^{-1}(B^4_\rho(\vec{\Phi}_\infty(x)))}\  \nabla\lf(\varphi_{\sigma}\star\phi(\vec{\Phi})\rg)\cdot\nabla\lf[f^{p-1}\ \vec{H}\rg]\ dx^2\\[5mm]
\ds+\,2\ p\ \sigma^2\ \int_{\vec{\Phi}_k^{-1}(B^4_\rho(\vec{\Phi}_\infty(x)))}\  \nabla\lf(\varphi_{\sigma}\star\phi(\vec{\Phi})\rg)\cdot
\lf[\,f^{p-1}\, H\ \nabla\vec{n}-\, e^{-2\la}\ f^{p-1}\, \lf<\nabla\vec{n}\dot{\otimes}\nabla\vec{n};\nabla\vec{\Phi}\rg>\rg]\ dx^2\\[5mm]
\ds-\int_{\vec{\Phi}_k^{-1}(B^4_\rho(\vec{\Phi}_\infty(x)))}\ \varphi_{\sigma}\star\phi(\vec{\Phi})\ \lf(\lf[1+\sigma^2\ (1-p)\, f^p+\, p\, \sigma^2\ f^{p-1}\rg]\ \vec{\Phi}\,|\nabla\vec{\Phi}|^2-4\, p\,\sigma^2\ f^{p-1}\ \vec{H}\rg)\ dx^2\quad.
\end{array}
\ee
Observe that $\|\p^2_{y_iy_j}(\varphi_{\sigma}\star\phi)\|_\infty\le \sigma^{-1/p}$ hence integrating by parts $\ov{\nabla}$ and $(\ov{\nabla})^\perp$ in the second  line of (\ref{M-119}) as well as integrating by parts
$\nabla$ in the fourth line of (\ref{M-119}) and using (\ref{M-5}) as in the proof of the monotonicity formula, we obtain that all the terms in the first, second, third and fourth lines of the r.h.s. of (\ref{M-119}) vanish as k goes to $+\infty$. In the fifth line only the term $\int_{\vec{\Phi}_k^{-1}(B^4_\rho(\vec{\Phi}_\infty(x)))}\varphi_{\sigma}\star\phi(\vec{\Phi})\ \vec{\Phi}\,|\nabla\vec{\Phi}|^2 dx^2$ is not necessarily converging towards 0. Since we are considering the first and second canonical directions and since $\Phi^1$ and $\Phi^2$ are $O(\rho)$ in $\vec{\Phi}_k^{-1}(B^4_\rho(\vec{\Phi}_\infty(x)))$ and since $\|\phi\|_\infty\le \rho$ we obtain (\ref{M-118-b}) and lemma~\ref{lm-vanishing-boundary} is proved. \hfill $\Box$

\medskip

\noindent{\bf Proof of lemma~\ref{lm-strong} continued.} Denote $\vec{\Phi}'_k:=(\Phi^3_k,\Phi^4_k)$. By taking $\phi(y):=h(y_1,y_2)\ \chi_\rho(y_3,y_4)$ where $\chi_\rho$ is identically equal to $\rho$ on $B^2_\rho(1,0)$, is non negative, supported in $B^2_{2\rho}(1,0)$, we have for
$i=1,2$
\be
\label{M-120}
\limsup_{k\rightarrow +\infty}\sup_{\mbox{supp}(h)\subset B^2_\rho(\vec{\Phi}_\infty(x))\ ;\ \|d h\|_\infty\le \rho^{-1}}\vec{\ep}_i\cdot\int_{B^4_{4\rho}(\vec{\Phi}_\infty(x))}\ast_kd\vec{y}\wedge (\chi_\rho\ dh+ h\, d\chi_\rho)=O(\rho^4)\quad.
\ee
Because of the existence of an approximate tangent plane at $\vec{\Phi}_\infty(x)$, which is equal to Span$\{\vec{\ep}_1,\vec{\ep}_2\}$, the asymptotic mass of the current in $B^4_{4\rho}(\vec{\Phi}_\infty(x))$ contained
in the support of $d\chi_\rho$ which is included in $B^2_{4\rho}(0,0)\times (B^2_{2\,\rho}(1,0)\setminus B^2_{\rho}(1,0))$ is a $o(\rho^2)$. Hence we deduce for $i=1,2$
\be
\label{M-121}
\limsup_{k\rightarrow +\infty}\sup_{\mbox{supp}(h)\subset B^2_\rho(0,0)\ ;\ \|d h\|_\infty\le \rho^{-1}}\int_{B^2_{\rho}(0,0)\times B^2_{\rho}(1,0)}\p_{y_i}h\ d{\mathcal H}^2\res\vec{\Phi}_k(\Sigma)=o(\rho^2)\quad.
\ee
This implies, using (\ref{M-108}),
\be
\label{M-122}
\limsup_{k\rightarrow +\infty}\sup_{\mbox{supp}(h)\subset B^2_\rho(0,0)\ ;\ \|d h\|_\infty\le \rho^{-1}}\int_{B^2_{\rho}(0,0)}N_k(y)\ \p_{y_i}h\ d{\mathcal L}^2=o(\rho^2)\quad,
\ee
where $N_k(y)$ is the number of pre-images of $y=(y_1,y_2)$ by $\vec{\zeta}_k$. Since $M(B^2_\rho(0,0)\cap\vec{\zeta}_k(\Sigma) )\simeq\rho^2$ we then have
\be
\label{M-123}
\limsup_{k\rightarrow +\infty}\ \sum_{i=1}^2\sup_{\mbox{supp}(h)\subset B^2_\rho(0,0)\ ;\ \|d h\|_\infty\le \rho^{-1}}\frac{\ds\int_{B^2_{\rho}(0,0)}N_k(y)\ \p_{y_i}h\ dy^2}{\ds\int_{B^2_{\rho}(0,0)}N_k(y)\ dy^2}=o_\rho(1)\quad.
\ee
The quantity on the numerator of (\ref{M-123}) is almost but not quite the {\it Flat Norm}\footnote{The flat norm would have been 
\[
\sup_{\mbox{supp}(X)\subset B^2_\rho(0,0)\ ;\ \|\mbox{div}X\|_\infty\le \rho^{-1}}\int_{B^2_{\rho}(0,0)}N_k(y)\ \mbox{div}(X)\ dy^2
\] and cannot a-priori be controlled by the numerator of (\ref{M-123}).}
 of the relative boundary in $B^2_\rho(0,0)$ of the 2 dimensional {\it integer rectifiable current} given by $C_k(\rho):=[N_k(y)\ dy^2]\res B^2_\rho(0,0)$  while the denominator equals it's total mass. 

In \cite{PiRi} the following inequality is proved. For any measurable function $f$ on the 2 dimensional unit ball $B_1(0)$ the following inequality holds
\be
\label{M-123-Rep}
\begin{array}{l}
\ds\lf\|f-\frac{1}{|B_{1/2}(0)|}\int_{B_{1/2}(0)} f(y)\ dy^2\rg\|_{L^{1,\infty}(B_{1/2}(0))}\\[5mm]
\ds\quad\quad\le C\ \sup\lf\{  \int_{B_{1}(0)} f(y)\ \nabla\phi(y)\ dy^2 \ ;\ \phi\in C^\infty_0(B_1(0))\quad \|\nabla\phi\|_\infty\le 1  \rg\}\quad.
\end{array}
\ee
Combining (\ref{M-123}) and (\ref{M-123-Rep}) gives that 
\be
\label{M-124-Rep}
\limsup_{k\rightarrow +\infty}\ \lf\|N_k(\rho\, x)-\frac{1}{|B_{1/2}(0)|}\int_{B_{1/2}(0)} N_k(\rho\, y)\ dy^2\rg\|_{L^{1,\infty}(B_{1/2}(0))}=o_\rho(1)\quad.
\ee
This shows that the average $\frac{1}{|B_{1/2}(0)|}\int_{B_{1/2}(0)} N_k(\rho\, y)\ dy^2$ is $o_\rho(1)$ close to an integer $n_k^\rho\in {\N}^\ast$ as $k$tends to infinity and that  
\be
\label{M-125-Rep}
\limsup_{k\rightarrow +\infty}\ \lf\|N_k(\rho\, x)-n_k^\rho\rg\|_{L^{1,\infty}(B_{1/2}(0))}=o_\rho(1)\quad.
\ee
Since this integer is bounded and bounded away from zero, modulo extraction of a subsequence we can assume that $n_k^\rho=n^\rho$ is independent of $k$ and, taking a sequence of radii
$\rho_j\rightarrow 0$ we can also assume that $n^{\rho_j}$ is independent of $j$ and we have the the existence of $n\in {\N}^\ast$ such that
\be
\label{M-126-Rep}
\lim_{j\rightarrow +\infty}\limsup_{k\rightarrow +\infty}\ \lf\|N_k(\rho_j\, x)-n\rg\|_{L^{1,\infty}(B_{1/2}(0))}=0\quad,
\ee
this proves (\ref{M-112}) and this concludes the proof of lemma~\ref{lm-strong}. \hfill $\Box$

\begin{Lm}
\label{lm-integ} {\bf[Convergence to a Bubble Tree of conformal ``integer target harmonic'' maps]}
Under the assumptions of theorem~\ref{th-limit}, we have that  one we can extract a subsequence such that the integer varifold
$|\vec{\Phi}_k(\Sigma)|$ converges to an integer rectifiable varifold supported by a finite union of the images by target harmonic conformal $W^{1,2}-$maps of Riemann surfaces .\hfill $\Box$
\end{Lm}

We adopt the same notations as in the proof of lemma~\ref{lm-strong} and assume to simplify the presentation that (\ref{m-112-b}) holds where we recall among other things that $x$ is chosen also to be a Lebesgue point for $\nabla\vec{\Phi}_\infty(x)$. One has
\be
\label{repa-1}
\lim_{\rho\rightarrow 0}\lim_{k\rightarrow +\infty}\frac{\ds\int_{\vec{\Phi}_\infty^{-1}(B^4_\rho(\vec{\Phi}_\infty(x)))} |\p_{x_1}\vec{\Phi}_k\wedge\p_{x_2}\vec{\Phi}_k|\ dx^2}{\ds\vec{\ep}_1\wedge\vec{\ep}_2\cdot\int_{\vec{\Phi}_\infty^{-1}(B^4_\rho(\vec{\Phi}_\infty(x)))} \p_{x_1}\vec{\Phi}_\infty\wedge\p_{x_2}\vec{\Phi}_\infty \ dx^2}=N_x\quad.
\ee
Observe also that\footnote{We recall among other things that $x$ is chosen also to be a Lebesgue point for $\nabla\vec{\Phi}_\infty$ and that $\nabla\vec{\Phi}_\infty(x)=\nabla\vec{\Xi}^\al(x)$.}
The lower semicontinuity of the norm gives
\be
\label{repa-2}
\begin{array}{l}
\ds\liminf_{k\rightarrow +\infty}\int_{\vec{\Phi}_\infty^{-1}(B^4_\rho(\vec{\Phi}_\infty(x)))} |\p_{x_1}\vec{\Phi}_k\wedge\p_{x_2}\vec{\Phi}_k|\ dx^2=\liminf_{k\rightarrow +\infty}\int_{\vec{\Phi}_\infty^{-1}(B^4_\rho(\vec{\Phi}_\infty(x)))} 2^{-1}\, |\nabla\vec{\Phi}_k|^2\ dx^2\\[5mm]
\ds\quad\ge\int_{\vec{\Phi}_\infty^{-1}(B^4_\rho(\vec{\Phi}_\infty(x)))} 2^{-1}\, |\nabla\vec{\Phi}_\infty|^2\ dx^2\quad.
\end{array}
\ee
Hence combining (\ref{repa-1}) and (\ref{repa-2}) one gets
\[
\ds\lim_{\rho\rightarrow 0}\frac{\ds\int_{\vec{\Phi}_\infty^{-1}(B^4_\rho(\vec{\Phi}_\infty(x)))} 2^{-1}\, |\nabla\vec{\Phi}_\infty|^2\ dx^2}{\ds\int_{\vec{\Phi}_\infty^{-1}(B^4_\rho(\vec{\Phi}_\infty(x)))} |\p_{x_1}\vec{\Phi}_\infty\wedge\p_{x_2}\vec{\Phi}_\infty| \ dx^2}\le N_x=\pi^{-1}\lim_{\rho\rightarrow 0}\rho^{-2}\,\mu_\infty(B^4_{\rho}(\vec{\Phi}_\infty(x)))\quad.
\]
This gives, using the {\it Monotonicity Formula}, we have 
\be
\label{repa-3}
\quad\mbox{ for }\quad{\nu_\infty}\mbox{ a. e. }\ x\in D^2\setminus{\mathcal B}\quad\quad1\le\frac{ |\nabla\vec{\Phi}_\infty|^2(x)}{2\, |\p_{x_1}\vec{\Phi}_\infty\wedge\p_{x_2}\vec{\Phi}_\infty|(x)}\le  \pi^{-1}\, e^{2\,C}\ \mu_\infty(S^3)=K\quad.
\ee
Take $g_{ij}:=\p_{x_i}\vec{\Phi}_\infty\cdot\p_{x_j}\vec{\Phi}_\infty$ and introduce 
\[
\mbox{for a. e.}\quad x\in D^2\setminus {\mathcal L}^0_{\nabla\vec{\Phi}_\infty}\quad\quad\mu(x):=\frac{g_{11}-g_{22}+2\, i\ g_{12}}{g_{11}+g_{22}+2\sqrt{g_{11}g_{22}-g_{12}^2}}
\]
on $D^2\setminus {\mathcal L}^0_{\nabla\vec{\Phi}_\infty}$, with the above notations (\ref{repa-3}) can be recast in the following way
\[
4\le\frac{(g_{11}+g_{22})^2}{g_{11}g_{22}-g_{12}^2}\le\frac{4}{\pi^2}\, e^{4\, C}\ \mu_\infty^2(S^3)=4\, K^2\quad.
\]   
 Extend $\mu$ by zero on the whole ${\C}$. Observe that we have
\[
\lf\|\mu\rg\|_\infty^2\le \lf\|\frac{(g_{11}+g_{22})^2-4(g_{11}g_{22}-g_{12}^2)}{(g_{11}+g_{22})^2+4(g_{11}g_{22}-g_{12}^2)}\rg\|_{L^\infty(D^2\setminus {\mathcal L}^0_{\nabla\vec{\Phi}_\infty})}\le \frac{K^2-1}{K^2+1}<1\quad.
\]
Hence $\mu $ defines a compactly supported {\it Beltrami coefficient}. Consider the {\it normal solution} of the {\it Beltrami equation} given by theorem 4.24 of \cite{ImTa}
\[
\p_{\ov{z}}\varphi=\mu\ \p_{{z}}\varphi\quad.
\]
The quasiconformal map $\varphi$ realizes in particular an homeomorphism whose inverse $\varphi^{-1}$ is also quasiconformal in $W^{1,p}_{loc}({\C})$ for some $p>2$ and one has
\[
\p_{\ov{w}}\varphi^{-1}=\om\ \p_{{w}}\varphi^{-1}\quad,
\]
where $\om=-(\mu\ \p_z\varphi\ /\ov{\p_z\varphi})\circ\varphi^{-1}$. Being an homeomorphic map of {\it bounded distortion} in $W^{1,2}(\varphi(D^2))$ it is quasi-regular, the chain rule applies with $\vec{\Phi}_\infty$ (see theorem 16.13.3 of \cite{IwMa}) and $\vec{\Phi}_\infty\circ \varphi^{-1}\in W^{1,2}(\varphi(D^2))$. A classical computation gives
\[
\p_{w}\lf(\vec{\Phi}_\infty\circ \varphi^{-1}\rg)\cdot\p_{w}\lf(\vec{\Phi}_\infty\circ \varphi^{-1}\rg)=0\quad\mbox{ a. e. on }\varphi(D^2)\quad.
\]
``Pasting'' together all these conformal charts gives a smooth conformal structure on $\Sigma$ and a global quasi-conformal homeomorphism $\psi$ of $\Sigma$ such that $\vec{\Phi}_\infty\circ\psi$ is weakly conformal. Moreover, the condition for the image of $\Sigma$ by $\vec{\Phi}_\infty$ equipped with the integer multiplicity $N$ to be stationary is equivalent to (\ref{targ}).
It remains to show that $(N,\vec{\Phi}_\infty\circ\psi)$ defines an integer target harmonic map.

\medskip

We omit to mention the composition by $\psi$ and we simply write $\vec{\Phi}_\infty$ for $\vec{\Phi}_\infty\circ\psi$. We can apply lemma~\ref{lm-VII.1} to $\Sigma\setminus \bigcup_{l=1}^n B_{r_k}(a_l)$ where $r_k$ are ``nice cuts'' taken between $\ep/2$ and $\ep$ on which $\vec{\Phi}_k$ converges in $C^0$ to deduce, using Because of (\ref{targ-ai}), that there exists $n$ points $\vec{q}_{l,\rho}$ such that
$$\lf|(\vec{\Phi}_\infty)_\ast(N\,[\Sigma])\res \lf({\R}^4\setminus\bigcup_{l=1}^n B^4_{s_\rho}(\vec{q}_{l,\rho})\rg)\rg|$$ realizes an integer rectifiable stationary varifold in $S^3\setminus\bigcup_{l=1}^n B^4_{s_\rho}(\vec{q}_{l,\rho})$. 
This is equivalent to
\be
\label{M-127}
\int_{\Sigma\setminus \bigcup_{l=1}^n B_r(a_l)}N\,\lf[\sum_{i=1}^4\lf<\p_{y_i}\vec{X}(\vec{\Phi}_\infty)\ \nabla\Phi^i_\infty;\nabla\vec{\Phi}_\infty\rg>-N\,\vec{X}(\vec{\Phi}_\infty)\cdot\vec{\Phi}_\infty\ |\nabla\vec{\Phi}_\infty|^2\rg]\ dx^2=0\quad.
\ee
 We chose a sequence of radii $\rho_k\rightarrow 0$ such that 
\[
\forall\ l=1\cdots n\quad\quad\quad\vec{q}_{l,\rho_k}\rightarrow \vec{q}_{l,0}\in S^3\quad.
\]
Since $s_{\rho_k}\rightarrow 0$, $(\vec{\Phi}_\infty)_\ast(N\,[\Sigma])$ is stationary in $S^3\setminus\{\vec{q}_{1,0}\cdots\vec{q}_{n,0}\}$. Let $\chi_\delta(t)=\chi(t/\delta)$ where $\chi\in C^\infty_0([0,2],{\R}_+)$, $\chi$ is identically equal to one on $[0,1]$. For any arbitrary smooth vector field $\vec{X}$ from $\Gamma(TS^3)$ we proceed to the following decomposition :
\[
\vec{X}(\vec{q})=\sum_{l=1}^n\chi_\delta(|\vec{q}-\vec{q}_{l,0}|)\ \vec{X}+\vec{X}_\delta(\vec{q})\quad\mbox{where }\quad \vec{X}_\delta(\vec{q}):=\lf[ 1-\sum_{l=1}^n\chi_\delta(|\vec{q}-\vec{q}_{l,0}|)\  \rg]\ \vec{X}
\]
Since Supp$(\vec{X}_\delta)\subset {\R}^4\setminus\bigcup_{l=1}^n B^4_\delta(\vec{q}_{l,0})$ we have
\be
\label{M-128}
\int_{\Sigma}N\,\lf[\sum_{i=1}^4\lf<\p_{y_i}\vec{X}_\delta(\vec{\Phi}_\infty)\ \nabla\Phi^i_\infty;\nabla\vec{\Phi}_\infty\rg>-\vec{X}_\delta(\vec{\Phi}_\infty)\cdot\vec{\Phi}_\infty\ |\nabla\vec{\Phi}_\infty|^2\rg]\ dx^2=0
\ee
and we have
\be
\label{M-129}
\begin{array}{l}
\ds\lf|\int_{\Sigma}N\lf[\sum_{i=1}^4\lf<\p_{y_i}(\vec{X}-\vec{X}_\delta)(\vec{\Phi}_\infty)\ \nabla\Phi^i_\infty;\nabla\vec{\Phi}_\infty\rg>-(\vec{X}-\vec{X}_\delta)(\vec{\Phi}_\infty)\cdot\vec{\Phi}_\infty\ |\nabla\vec{\Phi}_\infty|^2\rg]\ dx^2\rg|\\[5mm]
\ds\quad\quad\quad\le\|\vec{X}\|_\infty\ \frac{1}{\delta}\sum_{l=1}^n\mu_\infty(B^4_{2\,\delta}(\vec{q}_{l,0}))+\|\nabla\vec{X}\|_\infty\ \sum_{l=1}^n\mu_\infty(B^4_{2\,\delta}(\vec{q}_{l,0}))=O(\delta)
\end{array}
\ee
where we are using the monotonicity formula. Combining (\ref{M-128}) and (\ref{M-129}) with $\delta\rightarrow 0$ we obtain that 
\be
\label{M-130}
\int_{\Sigma}N\,\lf[\sum_{i=1}^4\lf<\p_{y_i}\vec{X}(\vec{\Phi}_\infty)\ \nabla\Phi^i_\infty;\nabla\vec{\Phi}_\infty\rg>-\vec{X}(\vec{\Phi}_\infty)\cdot\vec{\Phi}_\infty\ |\nabla\vec{\Phi}_\infty|^2\rg]\ dx^2=0\quad.
\ee
What we have done for the whole $\Sigma$ can be done for any subdomain $\Omega$ assuming that the support of $\vec{X}$ is contained in a complement of an open neighborhood of $\vec{\Phi}_\infty(\p \Om)$.
We deduce that $\vec{\Phi}_\infty$ is integer target harmonic from $\Sigma$ into $S^3$. This concludes the proof of the lemma~\ref{lm-integ}.\hfill $\Box$

\section{The proof of theorem~\ref{th-main}.}

We consider the general case where $(\Sigma,g_{\vec{\Phi}_k})$ possibly degenerate in the moduli space. Modulo extraction of a subsequence, following Deligne-Mumford compactification described in section II of \cite{Ri4} we have a ``splitting'' of the original surface into collars, called also ``thin parts'' and and a Nodal Riemann surface $\ti{\Sigma}$ called also ``thick part''. The parts of the collars that contain no bubbles can be treated exactly as the necks in lemma~\ref{lm-no-neck}, indeed a collar has the conformal type of a degenerating annulus and, if such a collar contains no bubble, by definition, it means that on each sub-annulus of controlled conformal type (in each dyadic annulus in particular) there is no concentration measure $\nu_\infty$. Hence in a collar region containing no bubble 
the statement of lemma~\ref{lm-no-neck} applies word by word. The ''thick parts'' as well as the ``bubbles'' formed either in the thick parts or in the collars can be treated exactly as the surface $\Sigma$ in the compact case presented in the previous section. 
So we deduce theorem~\ref{th-main}.

 \renewcommand{\theequation}{A.\arabic{equation}}
\renewcommand{\theTh}{A.\arabic{Th}}
\renewcommand{\theProp}{A.\arabic{Prop}}
\renewcommand{\theLma}{A.\arabic{Lma}}
\renewcommand{\theCo}{A.\arabic{Co}}
\renewcommand{\theRm}{A.\arabic{Rm}}
\renewcommand{\theequation}{A.\arabic{equation}}
\setcounter{equation}{0} 
\reset
\appendix
\section{Appendix}

\begin{Lma}
\label{lm-A.2}
There exists a universal number $\ep_0(m)>0$ such that, for any $\vec{\Phi}$  smooth immersion of $\Sigma$, a smooth surface with boundary,  into $B^m_2(0)\setminus B^m_1(0)$ and satisfying
\be
\label{A-8}
\mbox{Area}(\vec{\Phi}(\Sigma))<3\,\pi\quad,
\ee
and
\be
\label{A-9}
\forall \ r\in(1,2)\quad\vec{\Phi}(\Sigma)\cap\p B^m_r(0)\ne \emptyset\quad\mbox{ and }\quad\vec{\Phi}( \p\Sigma)\subset \p \lf(B^m_2(0)\setminus B^m_1(0)\rg)\quad,
\ee
then
\be
\label{A-10}
\int_\Sigma|d\vec{n}|^2_{g_{\vec{\Phi}}}\ dvol_{\vec{\Phi}}\ge \ep_0(m)\quad.
\ee
\hfill $\Box$
\end{Lma}
\noindent{\bf Proof of lemma~\ref{lm-A.2}.} We argue by contradiction. We consider a sequence $\Sigma_k$ and $\vec{\Phi}_k$ such that
\be
\label{A-10-a}
\mbox{Area}(\vec{\Phi}_k(\Sigma_k))<3\,\pi\quad,
\ee
such that
\be
\label{A-11}
\forall \ r\in(1,2)\quad\vec{\Phi}_k(\Sigma_k)\cap\p B^m_r(0)\ne \emptyset\quad\mbox{ and }\quad\vec{\Phi}_k( \p\Sigma_k)\subset \p \lf(B^m_2(0)\setminus B^m_1(0)\rg)\quad,
\ee
and
\be
\label{A-12}
\lim_{k\rightarrow+\infty} \int_{\Sigma_k}|d\vec{n}|^2_{g_{\vec{\Phi}_k}}\ dvol_{\vec{\Phi}_k}=0\quad.
\ee
 Let $V_k$ be the oriented varifold associated to the immersion of $\vec{\Phi}_k$ with $L^2-$bounded second fundamental form (see \cite{Hu}). Using theorem 3.1 and 5.3.2 of \cite{Hu}, modulo extraction of a subsequence
$V_k$ varifold converges to an integer oriented varifold $V_\infty$ with generalized second fundamental form equal to zero and without boundary in $B_2(0)\setminus B_1(0)$. $V_\infty$ is then
stationary and included in an at most countable union of 2-planes. Using the constancy theorem \cite{Si} we deduce that  $V_\infty$ is an oriented varifold given by 
 at most countably many intersections of 2-planes with the annulus $B_2(0)\setminus B_1(0)$ with locally constant integer multiplicities.
We claim that the intersection between the closed set given by the support of $V_\infty$ and $\p B_r(0)\times G_2({\R}^m)$ is non empty for any $r\in (1,2)$.
Indeed, from the assumption (\ref{A-11}), using Simon's monotonicity formula, for any $r\in(1,2)$  and $0<\rho<\min\{2-r,r-1\}$, there exists $x_k^r\in \p B_r(0)$ such that
\[
\frac{2\pi}{3}\rho^2\le M\lf(\vec{\Phi}_k(\Sigma_k)\cap B^m_\rho(x_k^r) \rg)+\frac{\rho^2}{2}\int_{\Sigma_k}
|\vec{H}_{\vec{\Phi}_k}|^2\ dvol_{g_{\vec{\Phi}_k}}\quad.
\]
Using (\ref{A-12}) we deduce that for any $\rho<\min\{2-r,r-1\}$
\[
\mu_{V_\infty}(B_{r+\rho}(0)\setminus B_{r-\rho}(0))\ge \frac{2\pi}{3}\rho^2\quad.
\]
Hence the support of $V_\infty$ intersects all the $\p B_r(0)\times G_2({\R}^m)$ for any $r\in (1,2)$. We consider a sequence of radii $r_i>1$ and converging to 1.
The 2-planes belonging to the support of $V_\infty$ and intersecting $\p B_{r_i}(0)\times G_2({\R}^m)$ has to be constant for $i$ large enough. This implies that the support of $V_\infty$ contains the intersection between the annulus $B_2(0)\setminus B_1(0)$ and a plane touching $\ov{B_1(0)}$. This imposes
\[
\mu_{V_\infty}(B_2(0)\setminus B_1(0))\ge\ 3\, \pi\quad.
\]
The later contradicts (\ref{A-10-a}) and lemma~\ref{lm-A.2} is proved.\hfill $\Box$

\end{document}